\documentclass[12pt]{article}
\usepackage{amssymb, amsmath}
\usepackage[dvips]{graphics}

\usepackage{bbm}
\usepackage{enumitem}
\usepackage{setspace}
\usepackage[colorlinks=true]{hyperref}
\hypersetup{colorlinks   = true}
\hypersetup{linkcolor=blue}
\begin{document}

\newtheorem{The}{Theorem}[section]
\newtheorem{Lem}[The]{Lemma}
\newtheorem{Prop}[The]{Proposition}
\newtheorem{Cor}[The]{Corollary}
\newtheorem{Rem}[The]{Remark}
\newtheorem{Obs}[The]{Observation}
\newtheorem{SConj}[The]{Standard Conjecture}
\newtheorem{Titre}[The]{\!\!\!\! }
\newtheorem{Conj}[The]{Conjecture}
\newtheorem{Question}[The]{Question}
\newtheorem{Prob}[The]{Problem}
\newtheorem{Def}[The]{Definition}
\newtheorem{Not}[The]{Notation}
\newtheorem{Claim}[The]{Claim}
\newtheorem{Conc}[The]{Conclusion}
\newtheorem{Ex}[The]{Example}
\newtheorem{Fact}[The]{Fact}
\newtheorem{Formula}[The]{Formula}
\newtheorem{Formulae}[The]{Formulae}
\newcommand{\C}{\mathbb{C}}
\newcommand{\R}{\mathbb{R}}
\newcommand{\N}{\mathbb{N}}
\newcommand{\Z}{\mathbb{Z}}
\newcommand{\Q}{\mathbb{Q}}
\newcommand{\Proj}{\mathbb{P}}

\begin{center}

{\Large\bf Degeneration at $E_2$ of Certain Spectral Sequences }

\end{center}

\begin{center}

{\large Dan Popovici}

\end{center}

\vspace{1ex}

\noindent {\small {\bf Abstract.} We propose a Hodge theory for the spaces $E_2^{p,\,q}$ featuring at the second step either in the Fr\"olicher spectral sequence of an arbitrary compact complex manifold $X$ or in the spectral sequence associated with a pair $(N,\,F)$ of complementary regular holomorphic foliations on such a manifold. The main idea is to introduce a Laplace-type operator associated with a given Hermitian metric on $X$ whose kernel in every bidegree $(p,\,q)$ is isomorphic to $E_2^{p,\,q}$ in either of the two situations discussed. The surprising aspect is that this operator is not a differential operator since it involves a harmonic projection, although it depends on certain differential operators. We then use this Hodge isomorphism for $E_2^{p,\,q}$ to give sufficient conditions for the degeneration at $E_2$ of the spectral sequence considered in each of the two cases in terms of the existence of certain metrics on $X$. For example, in the Fr\"olicher case we prove degeneration at $E_2$ if there exists an SKT metric $\omega$ (i.e. such that $\partial\bar\partial\omega=0$) whose torsion is small compared to the spectral gap of the elliptic operator $\Delta' + \Delta''$ defined by $\omega$. In the foliated case, we obtain degeneration at $E_2$ under a hypothesis involving the Laplacians $\Delta'_N$ and $\Delta'_F$ associated with the splitting $\partial = \partial_N + \partial_F$ induced by the foliated structure.}

\vspace{1ex}

\section{Introduction}\label{section:introduction1}

 This paper comprises two parts as it gives various sufficient conditions of a metric nature for the degeneration at $E_2$ of the Fr\"olicher spectral sequence of compact complex manifolds (Part I) and for the degeneration at $E_2$ of the spectral sequence associated with a complementary pair of foliations $(N,\,F)$ on a compact complex manifold (Part II). The first part is meant to serve as a blueprint for the second part. We briefly describe in this introduction the main results and ideas in each of the two parts. 

\vspace{2ex}

 {\bf Part I.} Let $X$ be a compact complex manifold of dimension $n$. For every Hermitian metric $\omega$ on $X$ and every bidegree $(p,\,q)$, we denote by $p''=p''_{p,\,q}:C^{\infty}_{p,\,q}(X,\,\C)\longrightarrow\ker\Delta''$ the orthogonal projection of the space of smooth $(p,\,q)$-forms onto the kernel of the $\bar\partial$-Laplacian $\Delta'':=\bar\partial\bar\partial^{\star} + \bar\partial^{\star}\bar\partial$ associated with $\omega$ and acting on $(p,\,q)$-forms. Clearly, $p''$ depends on the metric $\omega$.

 Our first result asserts that the commutation of $\partial$ with $p''$ suffices to ensure the degeneration at $E_2$ of the Fr\"olicher spectral sequence of $X$ (a property that will be denoted throughout by $E_2(X)=E_{\infty}(X)$).

\begin{The}\label{The:E_2-commutation} Let $X$ be a compact complex manifold. If $X$ carries a Hermitian metric $\omega$ such that $p''\partial = \partial p''$ on all $(p,\,q)$-forms for all bidegrees $(p,\,q)$, the Fr\"olicher spectral sequence of $X$ degenerates at $E_2$.

\end{The}

 We introduce and use as our main tool in Part I the Laplace-type operator 

$$\widetilde{\Delta} : = \partial p''\partial^{\star} + \partial^{\star}p''\partial + \bar\partial\bar\partial^{\star} + \bar\partial^{\star} \bar\partial:C^{\infty}_{p,\,q}(X,\,\C)\longrightarrow C^{\infty}_{p,\,q}(X,\,\C), \hspace{2ex} p,q=0,\dots , n,$$

\noindent which is not a differential operator since $p''$ isn't. However, $\widetilde{\Delta}$ is a pseudo-differential operator since the projector $p''$ is even a smoothing such operator of finite rank. (Indeed, $\ker\Delta''$ is finite-dimensional thanks to $\Delta''$ being elliptic.) We put $\Delta'_{p''}:= \partial p''\partial^{\star} + \partial^{\star}p''\partial \geq 0$. The introduction of $\widetilde{\Delta}$ is justified by a Hodge Isomorphism Theorem that we prove (cf. Theorem \ref{The:Hodge-isom_E2}, the starting point of this work) stating that each space $E_2^{p,\,q}$ featuring at the second step in the Fr\"olicher spectral sequence of $X$ is isomorphic to the kernel of $\widetilde{\Delta} : C^{\infty}_{p,\,q}(X,\,\C)\longrightarrow C^{\infty}_{p,\,q}(X,\,\C)$. In particular, we propose in this paper a Hodge theory for the second step of the Fr\"olicher spectral sequence by means of an operator belonging to a class of operators involving partial derivatives that is larger than the class of differential operators.

 On the other hand, with every Hermtian metric $\omega$ on $X$ we associate the following zero-order operators of type $(0,\,0)$ depending only on the torsion of $\omega$: $\bar{S}_{\omega}:= [\bar\partial\omega\wedge\cdot,\,(\bar\partial\omega\wedge\cdot)^{\star}]\geq 0$,

\begin{equation}\label{eqn:R-bar-omega_def}Z_{\omega}:= [\tau_{\omega},\,\tau_{\omega}^{\star}] + [\partial\omega\wedge\cdot,\,(\partial\omega\wedge\cdot)^{\star}]\geq 0  \hspace{2ex} \mbox{and}  \hspace{2ex}     \bar{R}_{\omega}:=[\bar\tau_{\omega},\,\bar\tau_{\omega}^{\star}] - [\bar\partial\omega\wedge\cdot,\,(\bar\partial\omega\wedge\cdot)^{\star}],\end{equation}

\noindent where the notation is the standard one: $[A,\,B]:=AB-(-1)^{a\,b}BA$ denotes the graded commutator of any pair of endomorphisms $A,B$ of respective degrees $a,b$ of the graded algebra $C^{\infty}_{\bullet,\,\bullet}(X,\,\C)$ of smooth differential forms on $X$, while $\tau = \tau_{\omega}:=[\Lambda,\,\partial\omega\wedge\cdot]$ is the torsion operator of order zero and bidegree $(1,\,0)$ associated with $\omega$ (cf. [Dem97, VII, $\S.1$.]) and $\Lambda=\Lambda_{\omega}$ is the formal adjoint of the Lefschetz operator $L:=\omega\wedge\cdot$ w.r.t. the $L^2$ inner product induced by $\omega$ on differential forms.

 The second sufficient condition of a metric nature for the degeneration at $E_2$ of the Fr\"olicher spectral sequence of $X$ that we give is the existence on $X$ of an SKT metric $\omega$ (i.e. a Hermitian metric $\omega$ such that $\partial\bar\partial\omega=0$) whose torsion is {\it small} in the sense that the upper bound of the torsion operator $Z_{\omega}$ (which is bounded) is dominated by a certain fixed multiple of the smallest positive eigenvalue of the non-negative self-adjoint elliptic operator $\Delta' + \Delta''$ in every bidegree $(p,\,q)$.

\begin{The}\label{The:E_2_small-torsionSKT_introd} Let $X$ be a compact complex $n$-dimensional manifold. If $X$ carries an SKT metric $\omega$ whose torsion satisfies the condition

\begin{equation}\label{eqn:torsion-smallness_spectral-gap_introd}\sup\limits_{u\in C^{\infty}_{p,\,q}(X,\,\C),\,\,||u||=1}\langle\langle Z_{\omega}u,\,u\rangle\rangle \leq \frac{1}{3}\,\mbox{min}\,\bigg(\mbox{Spec}\,(\Delta' + \Delta'')^{p,\,q}\cap(0,\,+\infty)\bigg)\end{equation}

\noindent for all $p,q\in\{0,\dots , n\}$, then the Fr\"olicher spectral sequence of $X$ degenerates at $E_2$.

\end{The}

 By $(\Delta' + \Delta'')^{p,\,q}$ we mean the operator $\Delta' + \Delta''$ acting on $(p,\,q)$-forms, while $\mbox{Spec}\,(\Delta' + \Delta'')^{p,\,q}$ stands for its spectrum and $||\cdot||, \langle\langle\cdot,\,\cdot\rangle\rangle$ denote the $L^2$-norm, resp. the $L^2$-inner product induced by $\omega$ on differential forms. Thus, the r.h.s. in (\ref{eqn:torsion-smallness_spectral-gap_introd}) is a third of the size of the {\it spectral gap} of $\Delta' + \Delta''$, an important quantity standardly associated with a given metric $\omega$. 

\vspace{2ex}

 We also give a third sufficient condition of a metric nature for the degeneration at $E_2$ of the Fr\"olicher spectral sequence of $X$ as follows. It, too, involves the torsion of the metric.

\begin{The}\label{The:E_2-Rbar} Let $X$ be a compact complex manifold. If $X$ carries an SKT metric $\omega$ (i.e. such that $\partial\bar\partial\omega=0$) with the property

\begin{equation}\label{eqn:tau-bar-S}\langle\langle\bar{R}_{\omega}\,u,\,u\rangle\rangle = 0  \hspace{2ex} \mbox{for all}\hspace{1ex} (p,\,q)\mbox{-forms} \hspace{1ex} u\in\ker\Delta'_{p''}\cap\ker\Delta'' \hspace{1ex} \mbox{of all bidegrees}\hspace{1ex} (p,\,q),\end{equation}

\noindent the Fr\"olicher spectral sequence of $X$ degenerates at $E_2$.

\end{The}

 We do not know at this point whether the mere existence of an SKT metric on $X$ suffices to ensure that $E_2(X)=E_{\infty}(X)$ in the Fr\"olicher spectral sequence, although this seems to be the case on all the examples of compact complex manifolds $X$ that we are aware of \footnote{The author is grateful to L. Ugarte for informing him that in the special class of nilmanifolds $X$ of complex dimension $3$, $E_2(X)=E_{\infty}(X)$ whenever an SKT metric exists on $X$.}. It is tempting to wonder whether the connection $\partial + \bar\partial$ can be changed in such a way that the spectral gap of $\Delta'+\Delta''$ becomes larger than the triple of the torsion operator bound (which remains unchanged if the metric does) and thus satisfies hypothesis (\ref{eqn:torsion-smallness_spectral-gap_introd}). It is with this question in mind that we computed a Bochner-Kodaira-Nakano-type formula for the Witten Laplacians as Theorem \ref{The:BKN_Witten2} in $\S.$\ref{subsection:BKN_Witten} of the Appendix. We hope this formula will be of independent interest as the Witten twisting of the connection only changes the resulting cohomology (hence the associated spectral sequence) by an isomorphism.

Meanwhile, we shall see that if $X$ carries an SKT metric $\omega$ for which $\bar{R}_{\omega}=0$ as an operator in all bidegrees, then the Fr\"olicher spectral sequence of $X$ degenerates even at $E_1$ (cf. Remark \ref{Rem:SKT-R=0_E1}).

\vspace{3ex}

 {\bf Part II.} Let $X$ be a compact complex manifold of dimension $n$. Suppose the holomorphic tangent bundle $T^{1,\,0}X$ splits into the direct sum of two holomorphic subbundles $N, F$ both of which are Frobenius integrable. This means that

\begin{equation}\label{eqn:product-structure_def}(i)\,\,T^{1,\,0}X = N\oplus F \hspace{2ex} \mbox{and}\hspace{2ex} (ii)\,\,[N,\,N]\subset N, \hspace{1ex}   [F,\,F]\subset F,\end{equation}

\noindent where $[\cdot,\,\cdot]$ stands for the Lie bracket. We may call such a splitting an {\it integrable holomorphic almost product structure} on $X$ by analogy with the real counterpart (see e.g. [Rei58]). By definition, it consists in a pair of complementary regular holomorphic foliations (N,\,F). 

 Let $r, n-r$ denote the ranks of $N$, resp. $F$. The splitting of $T^{1,\,0}X$ induces, for each $k\in\{0, 1, \dots , n\}$, a splitting of $\Lambda^{k,\,0}T^{\star}X:=\Lambda^k(T^{1,\,0}X)^{\star}$ as 

\begin{equation}\label{eqn:types}\Lambda^{k,\,0}T^{\star}X = \bigoplus\limits_{p+q=k}\Lambda^pN^{\star}\otimes\Lambda^qF^{\star}.\end{equation}

\noindent The sections of $\Lambda^pN^{\star}\otimes\Lambda^qF^{\star}$ will be called $(k,\,0)$-forms (or simply $k$-forms) of $(N,\,F)$-type $(p,\,q)$. The space of smooth global such forms on $X$ will be denoted by ${\cal E}^{p,\,q}(X) = {\cal E}^{p,\,q}_{N,\,F}(X)$. This is a subspace of the space ${\cal E}^k(X)$ of all smooth global $(k,\,0)$-forms (w.r.t. the complex structure of $X$).

 Thanks to the integrability assumption on both $N$ and $F$, the operator $\partial : {\cal E}^k(X) \longrightarrow {\cal E}^{k+1}(X)$ splits as $\partial = \partial_N + \partial_F$, where $\partial_N : {\cal E}^{p,\,q}(X) \longrightarrow {\cal E}^{p+1,\,q}(X)$ differentiates in the $N$-directions and $\partial_F : {\cal E}^{p,\,q}(X) \longrightarrow {\cal E}^{p,\,q+1}(X)$ differentiates in the $F$-directions while both $\partial_N$ and $\partial_F$ are integrable (i.e. $\partial_N^2=0$ and $\partial_F^2=0$). We are interested here in finding sufficient conditions under which the spectral sequence associated in the usual way with such an integrable holomorphic almost product structure $(N,\,F)$ on $X$ degenerates at $E_2$. The degeneration at $E_2$ is the best we can hope for since the spaces $E_1^{p,\,q}$ are not even finite-dimensional due to the lack of ellipticity of $\Delta'_F$, one of the two Laplace-type operators associated with $\partial_N$, resp. $\partial_F$ and with a given Hermitian metric $\omega$ on $X$:

$$\Delta'_N = \partial_N\partial_N^{\star} + \partial_N^{\star}\partial_N: {\cal E}^{p,\,q}(X) \longrightarrow {\cal E}^{p,\,q}(X) \hspace{2ex} \mbox{and} \hspace{2ex} \Delta'_F = \partial_F\partial_F^{\star} + \partial_F^{\star}\partial_F : {\cal E}^{p,\,q}(X) \longrightarrow {\cal E}^{p,\,q}(X).$$   

\noindent However, the sum of these two Laplacians $\Delta'_N + \Delta'_F : {\cal E}^{p,\,q}(X) \longrightarrow {\cal E}^{p,\,q}(X)$ is an elliptic operator and this feature will be exploited in Part II.

 By analogy with the discussion in Part I, for every $p, q$ we consider the orthogonal projection $p'_F : {\cal E}^{p,\,q}(X)\longrightarrow \ker\Delta'_F$ and introduce the Laplace-type operator

$$\widetilde{\Delta'}:=\partial_Np'_F \partial_N^{\star} + \partial_N^{\star} p'_F\partial_N + \partial_F\partial_F^{\star} + \partial_F^{\star}\partial_F : {\cal E}^{p,\,q}(X) \longrightarrow {\cal E}^{p,\,q}(X)$$  

\noindent that will be our main tool in Part II. Unlike its analogue $\widetilde{\Delta}$ of Part I, $\widetilde{\Delta'}$ is not a pseudo-differential operator since $p'_F$ is a projector onto an infinite-dimensional space. Our hope, to be investigated in future work, is that $p'_F$ will be, at least under certain conditions, a Fourier integral operator (FIO) with a complex phase in the sense of Melin-Sj\"ostrand [MS74] and Boutet de Monvel-Sj\"ostrand [BS76]. We put $\Delta'_{N,\,p'_F}:= \partial_Np'_F \partial_N^{\star} + \partial_N^{\star} p'_F\partial_N \geq 0$.

 Compared to the situation discussed in Part I, the main difference and new difficulty involved in Part II is that $\Delta'_F$ is not elliptic (unlike $\Delta''$), so the trivial inequality $\widetilde{\Delta'} \geq \Delta'_F$ does not reduce us to an elliptic operator whose G${\mathring a}$rding inequality can induce a similar inequality for $\widetilde{\Delta'}$. Unless a notion of {\it ellipticity} can be satisfactorily defined for $\widetilde{\Delta'}$ under certain conditions and subsequently used to infer that under those conditions $\ker\widetilde{\Delta'}$ is finite-dimensional and $\mbox{Im}\,\widetilde{\Delta'}$ is closed (our goal for future work), the only option available to us in this work is to ensure that under a certain hypothesis (see (\ref{eqn:Garding_sufficient-cond})) we have

$$\Delta'_{N,\,p'_F} \geq (1-\varepsilon)\,\Delta'_N, \hspace{2ex} \mbox{hence} \hspace{2ex} \widetilde{\Delta'} \geq (1-\varepsilon)\,(\Delta'_N + \Delta'_F)$$

\noindent for some constant $0<\varepsilon<1$. We then use the ellipticity of $\Delta'_N + \Delta'_F$ to deduce a G${\mathring a}$rding-type inequality for $\widetilde{\Delta'}$ from the G${\mathring a}$rding inequality satisfied by $\Delta'_N + \Delta'_F$. 

 The main results obtained in Part II can be summed up as follows.

\begin{The}\label{The:PartII_main} Let $X$ be a compact complex manifold endowed with a pair of complementary regular holomorphic foliations $(N,\,F)$. Suppose $X$ carries a Hermitian metric $\omega$ such that $[\partial_N,\,\partial_F^{\star}]=0$ and 

\begin{equation}\label{eqn:Garding_sufficient-cond_introd}\ker(\Delta'_N:{\cal E}^{p,\,q}(X)\rightarrow{\cal E}^{p,\,q}(X)) + \ker(\Delta'_F:{\cal E}^{p,\,q}(X)\rightarrow{\cal E}^{p,\,q}(X)) = {\cal E}^{p,\,q}(X)\end{equation}

\noindent for all $p,q$. Then 

\vspace{1ex}

\noindent $(a)$\, the following {\bf Hodge isomorphism} holds:

\begin{equation}\label{eqn:Hodge-isom_E2_NF_introd}{\cal H}_{\widetilde{\Delta'}}^{p,\,q}(N,\,F):=\ker\bigg(\widetilde{\Delta'}:{\cal E}^{p,\,q}(X)\longrightarrow{\cal E}^{p,\,q}(X)\bigg)\simeq E_2^{p,\,q}(N,\,F), \hspace{3ex} \alpha\longmapsto \bigg[[\alpha]_{\partial_F}\bigg]_{d_1},\end{equation}

\noindent where the $E_2^{p,\,q}(N,\,F)$ are the spaces of $(N,\,F)$-type $(p,\,q)$ featuring at step $2$ in the spectral sequence induced by $(N,\,F)$. Thus, every class $[[\alpha]_{\partial_F}]_{d_1}\in E_2^{p,\,q}(N,\,F)$ contains a unique $\widetilde{\Delta'}$-harmonic representative $\alpha$. In particular, $\mbox{dim}_{\C}E_2^{p,\,q}(N,\,F)<+\infty$ for all $p,q$.

\noindent $(b)$\, the spectral sequence induced by $(N,\,F)$ {\bf degenerates at $E_2$}.

\end{The}

 We end Part II with the study of a sufficient condition ensuring the anti-commutation identity $[\partial_N,\,\partial_F^{\star}]=0$ (i.e. $\partial_N\partial_F^{\star} = - \partial_F^{\star} \partial_N$) that is one of the two hypotheses we make in Theorem \ref{The:PartII_main}. Theorem \ref{The:delN-delFstar_anticomm} asserts that this anti-commutation holds if we work with a product metric $\omega = \omega_N + \omega_F$ that is {\it bundle-like} (cf. Definition \ref{Def:bundle-like-metric} borrowed from Reinhart's work [Rei59] and meant for our complex setting) and if $\partial_N\omega_N=0$ (i.e. $\omega$ is K\"ahler in the $N$-directions).  

\vspace{1ex}

 The sufficient degeneration conditions obtained in this paper can be looked at in the context of Deligne's classical work [Del68] giving degeneration criteria in terms of Lefschetz isomorphisms.

\vspace{2ex}

\noindent {\bf Acknowledgments.} This work was carried out at the UMI CNRS-CRM in Montr\'eal, specifically at the CIRGET and the Department of Mathematics of the UQ\`AM. The author is grateful to the CNRS for making his six-month stay here possible. He is also grateful to V. Apostolov and S. Lu for their hospitality and interest in many mathematical topics within and without the scope of this work, as well as to E. Giroux and C. Mourougane for stimulating discussions. Many thanks are also due to J. Sj\"ostrand for his interest in some of the matters discussed in this work. 

 This work has grown from a problem about the cohomology of nilmanifolds that was brought to the author's attention by S. Rollenske and the interest in it was later entertained by discussions with J. Ruppenthal. The author wishes to thank them both for stimulating discussions.

\vspace{2ex}

\begin{center} {\Large\bf Part I: $E_2$ degeneration of the Fr\"olicher spectral sequence}

\end{center}

\section{Review of standard material}\label{section:background} Let $X$ be a compact complex manifold of dimension $n$. Recall that the Fr\"olicher spectral sequence of $X$ is associated with the double complex $C^{\infty}_{\bullet,\,\bullet}(X,\,\C)$ defined by the total differential $d=\partial + \bar\partial$. This means that at step $0$ we put $E_0^{p,\,q}:=C^{\infty}_{p,\,q}(X,\,\C)$ and consider the differentials $d_0:=\bar\partial : E_0^{p,\,q}\rightarrow E_0^{p,\,q+1}$ for all $p,q\in\{0, 1, \dots , \,n\}$, so the groups $E_1^{p,\,q}$ at step $1$ in the spectral sequence are defined as the cohomology groups of the complex

\begin{equation}\label{eqn:E_0_complex}\cdots\stackrel{d_0}{\longrightarrow}E_0^{p,\,q-1}\stackrel{d_0}{\longrightarrow}E_0^{p,\,q}\stackrel{d_0}{\longrightarrow}E_0^{p,\,q+1}\stackrel{d_0}{\longrightarrow}\cdots,\end{equation} 

\noindent i.e. the $E_1^{p,\,q} : = H^q(E_0^{p,\,\bullet},\,d_0) = H^{p,\,q}_{\bar\partial}(X,\,\C)$ are the Dolbeault cohomology groups of $X$. The differentials $d_1$ are induced by $\partial$:

\begin{equation}\label{eqn:E_1_complex}\cdots\stackrel{d_1}{\longrightarrow}E_1^{p-1,\,q}\stackrel{d_1}{\longrightarrow}E_1^{p,\,q}\stackrel{d_1}{\longrightarrow}E_1^{p+1,\,q}\stackrel{d_1}{\longrightarrow}\cdots,\end{equation} 

\noindent i.e. for any form $\alpha\in E_0^{p,\,q}$ such that $\bar\partial\alpha = 0$, the class $[\alpha]_{\bar\partial}\in E_1^{p,\,q}$ is mapped by $d_1$ to the class $[\partial\alpha]_{\bar\partial}\in E_1^{p+1,\,q}$. Thus $d_1([\alpha]_{\bar\partial}) = [\partial\alpha]_{\bar\partial}$. This is meaningful since $\partial\bar\partial + \bar\partial\partial = 0$, so $\bar\partial(\partial\alpha) = -\partial(\bar\partial\alpha) = 0$ and thus $\partial\alpha$ defines indeed a cohomology class in $H^{p+1,\,q}_{\bar\partial}(X,\,\C) = E_1^{p+1,\,q}$. Moreover, the differential $d_1$ is well defined since $d_1([\alpha]_{\bar\partial})$ is independent of the choice of representative $\alpha$ of the $\bar\partial$-class $[\alpha]_{\bar\partial}$ as can be checked at once. Furthermore, $d_1^2=0$ (because $\partial^2=0$), so (\ref{eqn:E_1_complex}) is indeed a complex. The groups $E_2^{p,\,q}$ at step $2$ in the spectral sequence are defined as the cohomology groups of the complex (\ref{eqn:E_1_complex}), i.e.

\begin{equation}\label{eqn:E_2_def}E_2^{p,\,q}: = H^p(E_1^{\bullet,\,q},\,d_1) = \bigg\{\bigg[[\alpha]_{\bar\partial}\bigg]_{d_1}\bigg\slash\alpha\in C^{\infty}_{p,\,q}(X,\,\C)\cap\ker\bar\partial \hspace{1ex} \mbox{and} \hspace{1ex} \partial\alpha\in\mbox{Im}\,\bar\partial\bigg\}, \hspace{2ex} p,q\in\{0, 1, \dots , \,n\},    \end{equation}

\noindent so the elements of $E_2^{p,\,q}$ are $d_1$-classes of $\bar\partial$-classes. The process continues inductively by defining the groups $E_{r+1}^{p,\,q}$ at step $r+1$ as the cohomology groups of the complex $d_r: E_r^{p,\,q}\longrightarrow E_r^{p+r,\,q-r+1}$ already obtained at step $r$. At each step $r$ in the spectral sequence, the differentials $d_r$ are of type $(r,\,-r+1)$. We end up with $\C$-vector spaces $E_{\infty}^{p,\,q}$ and canonical isomorphisms

\begin{equation}\label{eqn:E_infty}H^k_{DR}(X,\,\C)\simeq\bigoplus_{p+q=k}E_{\infty}^{p,\,q}, \hspace{2ex} k=0,1, \dots , 2n,\end{equation}

\noindent relating the differential structure of $X$ encoded in the De Rham cohomology to its complex structure. The spectral sequence is said to degenerate at $E_r$ if $E_r^{p,\,q} = E_{r+1}^{p,\,q}$ for all $p,\, q$ (hence then also $E_r^{p,\,q} = E_{r+l}^{p,\,q} = E_{\infty}^{p,\,q}$ for all $l\geq 0$). This is a purely numerical property equivalent to the identities $\sum_{p+q=k}\mbox{dim}\,E_r^{p,\,q} = b_k:=\mbox{dim}\,H^k_{DR}(X,\,\C)$ for all $k\in\{0, \dots , 2n\}$ and also to the inequalities

\begin{equation}\label{eqn:E_r_degeneracy-condition}\sum\limits_{p+q=k}\mbox{dim}\,E_r^{p,\,q} \leq b_k  \hspace{2ex} \mbox{for all}\hspace{1ex} k\in\{0, \dots , 2n\}\end{equation}

\noindent since the reverse inequalities always hold thanks to (\ref{eqn:E_infty}) and to the obvious inequalities $\mbox{dim}\,E_1^{p,\,q} \geq \dots \geq \mbox{dim}\,E_r^{p,\,q} \geq \mbox{dim}\,E_{r+1}^{p,\,q} \geq\dots$. All these dimensions are, of course, always finite by compactness of $X$ and, for example, ellipticity of $\Delta''$. The degeneracy at $E_r$ of the spectral sequence will be denoted by $E_r(X) = E_{\infty}(X)$. For further details, see e.g. [Dem96].

\section{Pseudo-differential Laplacian and Hodge isomorphism for $E_2^{p,\,q}$}\label{section:nondiff-laplacian_Hodge}

 Let $\omega$ be an arbitrary Hermitian metric on $X$. Consider the formal adjoints $\partial^{\star}, \bar\partial^{\star}$ of  $\partial$, resp. $\bar\partial$ w.r.t. the $L^2$ inner product defined by $\omega$ and the usual Laplace-Beltrami operators $\Delta', \Delta'' : C^{\infty}_{p,\,q}(X,\,\C)\longrightarrow C^{\infty}_{p,\,q}(X,\,\C)$ defined as $\Delta'=\partial\partial^{\star} + \partial^{\star}\partial$ and $\Delta''=\bar\partial\bar\partial^{\star} + \bar\partial^{\star}\bar\partial$. It is standard that they are elliptic, self-adjoint and non-negative differential operators of order $2$ that induce $3$-space orthogonal decompositions

\begin{equation}\label{eqn:3-space-decomp}C^{\infty}_{p,\,q}(X,\,\C) = \ker\Delta' \oplus \mbox{Im}\,\partial \oplus \mbox{Im}\,\partial^{\star} \hspace{2ex} \mbox{and} \hspace{2ex} C^{\infty}_{p,\,q}(X,\,\C) = \ker\Delta'' \oplus \mbox{Im}\,\bar\partial \oplus \mbox{Im}\,\bar\partial^{\star}\end{equation}

\noindent where the harmonic spaces $\ker\Delta':={\cal H}^{p,\,q}_{\Delta'}(X,\,\C)$, $\ker\Delta'':={\cal H}^{p,\,q}_{\Delta''}(X,\,\C)$ are finite dimensional while

\begin{equation}\label{eqn:ker-del_ker_laplacian}\ker\partial =\ker\Delta' \oplus \mbox{Im}\,\partial \hspace{2ex} \mbox{and} \hspace{2ex} \ker\bar\partial = \ker\Delta'' \oplus \mbox{Im}\,\bar\partial.\end{equation}

\noindent We denote by 

\begin{equation}\label{eqn:p'p''_def}p'=p'_{p,\,q} : C^{\infty}_{p,\,q}(X,\,\C)\longrightarrow\ker\Delta' \hspace{2ex} \mbox{and} \hspace{2ex} p'' = p''_{p,\,q}: C^{\infty}_{p,\,q}(X,\,\C)\longrightarrow\ker\Delta''\end{equation}

\noindent the orthogonal projections defined by the orthogonal splittings (\ref{eqn:3-space-decomp}) onto the $\Delta'$-harmonic, resp. the $\Delta''$-harmonic
spaces in bidegree $(p,\,q)$. Similarly, let

\begin{equation}\label{eqn:p'p''_perp_def}p'_{\perp}: C^{\infty}_{p,\,q}(X,\,\C)\longrightarrow\mbox{Im}\,\Delta' = \mbox{Im}\,\partial \oplus \mbox{Im}\,\partial^{\star} \hspace{2ex} \mbox{and} \hspace{2ex} p''_{\perp}: C^{\infty}_{p,\,q}(X,\,\C)\longrightarrow\mbox{Im}\,\Delta'' = \mbox{Im}\,\bar\partial \oplus \mbox{Im}\,\bar\partial^{\star}\end{equation}

\noindent denote the orthogonal projections onto $(\ker\Delta')^{\perp} = \mbox{Im}\,\Delta'$, resp. onto $(\ker\Delta'')^{\perp} = \mbox{Im}\,\Delta''$. Note that the operators $p', p'', p'_{\perp}, p''_{\perp}$ are not differential operators and depend on the metric $\omega$. They clearly satisfy the properties:

\begin{equation}\label{eqn:orth-proj_prop}p' = (p')^{\star} = (p')^2, \hspace{2ex} p'' = (p'')^{\star} = (p'')^2, \hspace{2ex} p'_{\perp} = (p'_{\perp})^{\star} = (p'_{\perp})^2, \hspace{2ex} p''_{\perp} = (p''_{\perp})^{\star} = (p''_{\perp})^2.\end{equation}

\vspace{2ex}

 We start by giving a metric interpretation of the spaces $E_2^{p,\,q}$ in the Fr\"olicher spectral sequence of $X$.

\begin{Prop}\label{Prop:Hpq-tilde_def} For every $p,\,q=0,1, \dots ,\, n$, define the $\omega$-dependent $\C$-vector space 

\begin{equation}\label{eqn:Hpq-tilde_def}\widetilde{H}^{p,\,q}(X,\,\C):=\ker(p''\circ\partial)\cap\ker\bar\partial \bigg\slash \bigg(\mbox{Im}\,\bar\partial + \mbox{Im}\,(\partial_{|\ker\bar\partial})\bigg)\end{equation}

\noindent in which all the kernels and images involved are understood as subspaces of $C^{\infty}_{p,\,q}(X,\,\C)$. For every $C^{\infty}$ $(p,\,q)$-form $\alpha\in\ker(p''\circ\partial)\cap\ker\bar\partial$, let $\widetilde{[\alpha]}\in\widetilde{H}^{p,\,q}(X,\,\C)$ denote the class of $\alpha$ modulo $\mbox{Im}\,\bar\partial + \mbox{Im}\,(\partial_{|\ker\bar\partial})$. Then, for every $p,\,q$, the following linear map

\begin{eqnarray}\label{eqn:T_def} T = T^{p,\,q} : \widetilde{H}^{p,\,q}(X,\,\C) \longrightarrow E_2^{p,\,q}, \hspace{3ex} \widetilde{[\alpha]}\longmapsto\bigg[[\alpha]_{\bar\partial}\bigg]_{d_1},\end{eqnarray}

\noindent is {\bf well defined} and an {\bf isomorphism}.

\end{Prop}

\noindent {\it Proof.} First note that the inclusion $\mbox{Im}\,\bar\partial + \mbox{Im}\,(\partial_{|\ker\bar\partial})\subset\ker(p''\circ\partial)\cap\ker\bar\partial$ does hold, so the space $\widetilde{H}^{p,\,q}(X,\,\C)$ is meaningful. Indeed, $\mbox{Im}\,\bar\partial\subset\ker\bar\partial$ trivially and $\mbox{Im}\,\bar\partial\subset\ker(p''\circ\partial)$ because for every form $u$, $p''(\partial\bar\partial u) = -p''\bar\partial\partial u =0$ since $\mbox{Im}\,\bar\partial$ is orthogonal onto $\ker\Delta''$ (see (\ref{eqn:3-space-decomp})), so $p''\bar\partial = 0$. Thus $\mbox{Im}\,\bar\partial\subset\ker(p''\circ\partial)\cap\ker\bar\partial$. Moreover, $\mbox{Im}\,(\partial_{|\ker\bar\partial})\subset\ker(p''\circ\partial)$ because $\partial^2=0$ and $\mbox{Im}\,(\partial_{|\ker\bar\partial})\subset\ker\bar\partial$ because for any form $v\in\ker\bar\partial$, we have $\bar\partial(\partial v) = -\partial(\bar\partial v) =0$. Thus $\mbox{Im}\,(\partial_{|\ker\bar\partial})\subset\ker(p''\circ\partial)\cap\ker\bar\partial$.

  Then note that for any $\widetilde{[\alpha]}\in \widetilde{H}^{p,\,q}(X,\,\C)$, we do have $[\alpha]_{\bar\partial}\in\ker d_1$, so the $d_1$-class $[[\alpha]_{\bar\partial}]_{d_1}$ is a meaningful element of $E_2^{p,\,q}$. Indeed, $d_1([\alpha]_{\bar\partial}) = [\partial\alpha]_{\bar\partial}$, $\partial\alpha\in\ker\bar\partial = \ker\Delta''\oplus\mbox{Im}\,\bar\partial$ (because $\alpha\in\ker\bar\partial$ and (\ref{eqn:ker-del_ker_laplacian}) holds) and $p''(\partial\alpha) = 0$ (because $\alpha\in\ker(p''\circ\partial)$). The last two relations amount to $\partial\alpha\in\mbox{Im}\,\bar\partial$. This is equivalent to $[\partial\alpha]_{\bar\partial}=0$, i.e. to $d_1([\alpha]_{\bar\partial}) = 0$.

 To complete the proof of the well-definedness of $T$, it remains to show that $[[\alpha]_{\bar\partial}]_{d_1}$ does not depend on the choice of representative $\alpha$ of the class $\widetilde{[\alpha]}$, i.e. that the zero element of $\widetilde{H}^{p,\,q}(X,\,\C)$ is mapped by $T$ to the zero element of $E_2^{p,\,q}$. To prove this, let $\alpha\in\ker(p''\circ\partial)\cap\ker\bar\partial$ be a $(p,\,q)$-form such that $\alpha = \bar\partial u + \partial v$ with $v\in\ker\bar\partial$. We want to show that $[[\alpha]_{\bar\partial}]_{d_1} = 0\in E_2^{p,\,q}$, i.e. that $[\alpha]_{\bar\partial} = d_1([\beta]_{\bar\partial})$ or equivalently that $[\alpha]_{\bar\partial} = [\partial\beta]_{\bar\partial}$ for some $\beta\in C^{\infty}_{p-1,\,q}(X,\,\C)$ such that $\bar\partial\beta=0$. This is equivalent to showing that $\alpha = \partial\beta + \bar\partial\gamma$ for some $\beta\in C^{\infty}_{p-1,\,q}(X,\,\C)$ such that $\bar\partial\beta=0$ and some $\gamma\in C^{\infty}_{p,\,q-1}(X,\,\C)$. We can choose $\beta:=v$ and $\gamma:=u$. 

 To prove that $T$ is injective, let $\alpha\in\ker(p''\circ\partial)\cap\ker\bar\partial$ be a $(p,\,q)$-form s.t. $T(\widetilde{[\alpha]}) = [[\alpha]_{\bar\partial}]_{d_1} = 0$. Then $[\alpha]_{\bar\partial} = [\partial\beta]_{\bar\partial}$ for some $\beta\in C^{\infty}_{p-1,\,q}(X,\,\C)$ such that $\bar\partial\beta=0$. Hence $\alpha = \partial\beta + \bar\partial\gamma$ for some $\gamma\in C^{\infty}_{p,\,q-1}(X,\,\C)$. Thus, $\alpha\in\mbox{Im}\,\bar\partial + \mbox{Im}\,(\partial_{|\ker\bar\partial})$, so $\widetilde{[\alpha]} = 0$.

 To prove that $T$ is surjective, let $[[\alpha]_{\bar\partial}]_{d_1}\in E_2^{p,\,q}$. Then $\bar\partial\alpha=0$ (i.e. $\alpha\in\ker\bar\partial$) and $d_1([\alpha]_{\bar\partial}) = [\partial\alpha]_{\bar\partial}=0$ (i.e. $\partial\alpha\in\mbox{Im}\,\bar\partial$, which is equivalent, since we already have $\partial\alpha\in\ker\bar\partial = \ker\Delta''\oplus\mbox{Im}\,\bar\partial$, to $p''(\partial\alpha)=0$, i.e. to $\alpha\in\ker(p''\circ\partial_{|\ker\bar\partial})$). Thus, $\alpha\in\ker(p''\circ\partial)\cap\ker\bar\partial$. It is clear that $[[\alpha]_{\bar\partial}]_{d_1} = T(\widetilde{[\alpha]})$ by definition of $T$.   \hfill $\Box$

\vspace{2ex}

 The isomorphism (\ref{eqn:T_def}) naturally prompts the introduction of a Laplace-type operator which, surprisingly, is not a differential operator. It will be the main tool of investigation in this paper.

\begin{Def}\label{Def:Delta-tilde_def} Let $(X,\,\omega)$ be a compact Hermitian manifold with $\mbox{dim}_{\C}X=n$. For every $p, q\in\{0, 1, \dots , n\}$, we define the operator $\widetilde{\Delta} : C^{\infty}_{p,\,q}(X,\,\C)\longrightarrow C^{\infty}_{p,\,q}(X,\,\C)$ by

\begin{equation}\label{eqn:Delta-tilde_def}\widetilde{\Delta} : = \partial p''\partial^{\star} + \partial^{\star}p''\partial + \bar\partial\bar\partial^{\star} + \bar\partial^{\star} \bar\partial.\end{equation}

\noindent In other words, we have

\begin{equation}\label{eqn:Delta-tilde_def_bis}\widetilde{\Delta} = \Delta'_{p''} + \Delta'', \hspace{2ex} \mbox{where} \hspace{1ex} \Delta'_{p''}:= \partial p''\partial^{\star} + \partial^{\star}p''\partial : C^{\infty}_{p,\,q}(X,\,\C)\longrightarrow C^{\infty}_{p,\,q}(X,\,\C).\end{equation}

\noindent Thus $\widetilde{\Delta}$ is the sum of a pseudo-differential regularising operator ($\Delta'_{p''}$) and an elliptic differential operator of order two (the classical $\bar\partial$-Laplacian $\Delta''$).

\end{Def}

 Clearly, $\widetilde{\Delta}$ is a non-negative self-adjoint operator whose kernel is $\ker\widetilde{\Delta} = \ker\Delta'_{p''}\cap\ker\Delta''$ and

\begin{equation}\label{eqn:ker_Delta-tilde}\ker\Delta'_{p''} = \ker(p''\circ\partial)\cap\ker(p''\circ\partial^{\star}) \supset \ker\partial\cap\ker\partial^{\star} = \ker\Delta'\end{equation}

\noindent because $\langle\langle\Delta'_{p''}u,\,u\rangle\rangle = ||p''\partial u||^2 + ||p''\partial^{\star}u||^2$. Actually, if we put $\Delta'_{p''_{\perp}}: = \partial p''_{\perp}\partial^{\star} + \partial^{\star}p''_{\perp}\partial$, then $0\leq\Delta'_{p''}\leq\Delta' = \Delta'_{p''} + \Delta'_{p''_{\perp}}$ since

\begin{eqnarray}\label{eqn:3-laplacians-sum}\nonumber\langle\langle\Delta'u,\,u\rangle\rangle = ||\partial u||^2 + ||\partial^{\star}u||^2 & = & ||p''\partial u||^2 + ||p''\partial^{\star}u||^2 + ||p''_{\perp}\partial u||^2 + ||p''_{\perp}\partial^{\star}u||^2\\
   & = & \langle\langle\Delta'_{p''}u,\,u\rangle\rangle + \langle\langle\Delta'_{p''_{\perp}}u,\,u\rangle\rangle\end{eqnarray}

\noindent for any form $u$. Indeed, for example, $\partial u = p''\partial u + p''_{\perp}\partial u$ and $p''\partial u\perp p''_{\perp}\partial u$, while $\langle\langle \partial^{\star}p''\partial u,\, u\rangle\rangle = \langle\langle p''\partial u,\, \partial u\rangle\rangle = \langle\langle p''\partial u,\, p''\partial u\rangle\rangle = ||p''\partial u||^2$.

\vspace{2ex}

 We now pause briefly to notice some of the properties of $\widetilde{\Delta}$. 

\begin{Lem}\label{Lem:Delta-tilde_prop} $(i)$\, If the metric $\omega$ is K\"ahler, then $\Delta'_{p''}=0$, so $\widetilde{\Delta} = \Delta''$.

$(ii)$\, For every $p,q=0,1,\dots , n$, let $(\psi_j^{p,\,q})_{1\leq j\leq h^{p,\,q}}$ be an arbitrary orthonormal basis of the $\Delta''$-harmonic space ${\cal H}_{\Delta''}^{p,\,q}(X,\C)\subset C^{\infty}_{p,\,q}(X,\,\C)$. Then $\Delta'_{p''}$ is given by the formula

\begin{equation}\label{eqn:Delta-tilde_formula}\Delta'_{p''}u = \sum\limits_{j=1}^{h^{p-1,\,q}}\langle\langle u,\,\partial\psi_j^{p-1,\,q}\rangle\rangle\,\partial\psi_j^{p-1,\,q} + \sum\limits_{j=1}^{h^{p+1,\,q}}\langle\langle u,\,\partial^{\star}\psi_j^{p+1,\,q}\rangle\rangle\,\partial^{\star}\psi_j^{p+1,\,q}, \hspace{3ex} u\in C^{\infty}_{p,\,q}(X,\,\C).\end{equation}

\noindent $(iii)$\, For all $p,q$, $\widetilde{\Delta}: C^{\infty}_{p,\,q}(X,\,\C)\longrightarrow C^{\infty}_{p,\,q}(X,\,\C)$ behaves like an elliptic self-adjoint differential operator in the sense that $\ker\widetilde{\Delta}$ is finite-dimensional, $\mbox{Im}\,\widetilde{\Delta}$ is closed and finite codimensional in $C^{\infty}_{p,\,q}(X,\,\C)$, there is an orthogonal (for the $L^2$ inner product induced by $\omega$) $2$-space decomposition

\begin{equation}\label{eqn:2-space-decomp_Delta-tilde}C^{\infty}_{p,\,q}(X,\,\C) = \ker\widetilde{\Delta}\bigoplus\mbox{Im}\,\widetilde{\Delta}\end{equation}

\noindent giving rise to an orthogonal $3$-space decomposition

\begin{equation}\label{eqn:3-space-decomp_Delta-tilde}C^{\infty}_{p,\,q}(X,\,\C) = \ker\widetilde{\Delta}\bigoplus\bigg(\mbox{Im}\,\bar\partial + \mbox{Im}\,(\partial_{|\ker\bar\partial})\bigg)\bigoplus\bigg(\mbox{Im}\,(\partial^{\star}\circ p'') +  \mbox{Im}\,\bar\partial^{\star}\bigg)  \end{equation}

\noindent in which $\ker\widetilde{\Delta} \oplus (\mbox{Im}\,\bar\partial + \mbox{Im}\,(\partial_{|\ker\bar\partial})) = \ker(p''\circ\partial)\cap\ker\bar\partial$, $\ker\widetilde{\Delta} \oplus (\mbox{Im}\,(\partial^{\star}\circ p'') + \mbox{Im}\,\bar\partial^{\star}) = \ker\bar\partial^{\star}\cap\ker(p''\circ\partial^{\star})$ and $(\mbox{Im}\,\bar\partial + \mbox{Im}\,(\partial_{|\ker\bar\partial}))\oplus(\mbox{Im}\,(\partial^{\star}\circ p'') +  \mbox{Im}\,\bar\partial^{\star}) = \mbox{Im}\,\widetilde{\Delta}$.

\noindent Moreover, $\widetilde{\Delta}$ has a compact resolvent which is a pseudo-differential operator $G$ of order $-2$, the Green's operator of $\widetilde{\Delta}$, hence the spectrum of $\widetilde{\Delta}$ is discrete and consists of non-negative eigenvalues that tend to $+\infty$.

\end{Lem}

\noindent {\it Proof.} $(i)$\, If $\omega$ is K\"ahler, $\Delta'=\Delta''$, hence $p'=p''$. Since $\ker\Delta'$ is orthogonal to both $\mbox{Im}\,\partial$ and   $\mbox{Im}\,\partial^{\star}$, $p'\circ\partial=0$ and $p'\circ\partial^{\star}=0$. Thus $p''\circ\partial=0$ and $p''\circ\partial^{\star}=0$, so $\Delta'_{p''}=0$.

$(ii)$\, Since $\ker\Delta''$ is finite-dimensional, $p'': C^{\infty}_{p,\,q}(X,\,\C)\longrightarrow\ker\Delta''$ is a regularising operator of finite rank defined by the $C^{\infty}$ kernel $\sum\limits_{j=1}^{h^{p,\,q}}\psi_j^{p,\,q}(x)\otimes(\psi_j^{p,\,q})^{\star}(y)$. Consequently, for every $u\in C^{\infty}_{p,\,q}(X,\,\C)$, 

\begin{equation}\label{eqn:p''_integral-formula}(p''u)(x) = \int\limits_X\sum\limits_{j=1}^{h^{p,\,q}}\psi_j^{p,\,q}(x)\,\langle u(y),\,\psi_j^{p,\,q}(y)\rangle\,dV_{\omega}(y), \hspace{2ex} \mbox{i.e.} \hspace{2ex} p''u = \sum\limits_{j=1}^{h^{p,\,q}} \langle\langle u,\,\psi_j^{p,\,q}\rangle\rangle\,\psi_j^{p,\,q}.\end{equation}

\noindent Taking successively $u=\partial^{\star}v$ with $v\in C^{\infty}_{p+1,\,q}(X,\,\C)$ and $u=\partial w$ with $w\in C^{\infty}_{p-1,\,q}(X,\,\C)$, we get

$$p''\partial^{\star}v = \sum\limits_{j=1}^{h^{p,\,q}}\langle\langle v,\,\partial\psi_j^{p,\,q}\rangle\rangle\,\psi_j^{p,\,q} \hspace{2ex} \mbox{and} \hspace{2ex} p''\partial w = \sum\limits_{j=1}^{h^{p,\,q}} \langle\langle w,\,\partial^{\star}\psi_j^{p,\,q}\rangle\rangle\,\psi_j^{p,\,q}.$$

\noindent Formula (\ref{eqn:Delta-tilde_formula}) follows at once from these identities.

 $(iii)$\, Since $\ker\widetilde{\Delta}\subset\ker\Delta''$ and the latter kernel is finite-dimensional thanks to $\Delta''$ being elliptic, $\ker\widetilde{\Delta}$ is finite-dimensional.

 The operator $\widetilde{\Delta}$ is elliptic pseudo-differential as the sum of an elliptic differential operator and a regularising one, so the elliptic theory applies to it. But we can also argue starting from the obvious inequality $\widetilde{\Delta}\geq\Delta''\geq 0$ (which follows from $\langle\langle\Delta'_{p''}u,\,u\rangle\rangle \geq 0$ for all $u$) and combining it with the G$\mathring{a}$rding inequality for the elliptic differential operator $\Delta''$. We get constants $\delta_1,\delta_2>0$ such that

\begin{equation}\label{eqn:Garding}\delta_2\,||u||^2_1\leq \langle\langle\Delta''u,\,u\rangle\rangle + \delta_1\,||u||^2\leq  \langle\langle\widetilde{\Delta}u,\,u\rangle\rangle + \delta_1\,||u||^2, \hspace{2ex} u\in C^{\infty}_{p,\,q}(X,\,\C),\end{equation}  

\noindent where $||\,\,||_1$ denotes the Sobolev norm $W^1$ and $||\,\,||$ denotes the $L^2=W^0$ norm. Since $\langle\langle\widetilde{\Delta}u,\,u\rangle\rangle \leq \frac{1}{2}\,||\widetilde{\Delta}u||^2 + \frac{1}{2}\,||u||^2$, we get

\begin{equation}\label{eqn:Garding_consequence}\delta_2\,||u||^2_1\leq \frac{1}{2}\,||\widetilde{\Delta}u||^2 + (\delta_1 + \frac{1}{2})\,||u||^2, \hspace{2ex} u\in C^{\infty}_{p,\,q}(X,\,\C).\end{equation}

\noindent This suffices to prove that $\mbox{Im}\,\widetilde{\Delta}$ is closed by the usual method using the Rellich Lemma (see e.g. [Dem96, 3.10, p. 18-19]). From closedness of $\mbox{Im}\,\widetilde{\Delta}$ and self-adjointness of $\widetilde{\Delta}$ we get (\ref{eqn:2-space-decomp_Delta-tilde}).

 Now (\ref{eqn:3-space-decomp_Delta-tilde}) is easily deduced from (\ref{eqn:2-space-decomp_Delta-tilde}) as follows. It is clear that

$$\mbox{Im}\,\widetilde{\Delta}\subset \mbox{Im}\,(\partial\circ p'') + \mbox{Im}\,(\partial^{\star}\circ p'') + \mbox{Im}\,\bar\partial + \mbox{Im}\,\bar\partial^{\star}.$$

\noindent Since $\mbox{Im}\,(\partial\circ p'') = \mbox{Im}\,(\partial_{|\ker\Delta''})$ and $\ker\Delta''\subset\ker\bar\partial$, we get $\mbox{Im}\,(\partial\circ p'')\subset \mbox{Im}\,(\partial_{|\ker\bar\partial})$, hence

\begin{equation}\label{eqn:image-delta-tilde_inclusion}\mbox{Im}\,\widetilde{\Delta}\subset \bigg(\mbox{Im}\,\bar\partial + \mbox{Im}\,(\partial_{|\ker\bar\partial})\bigg )\oplus \bigg(\mbox{Im}\,(\partial^{\star}\circ p'') +  \mbox{Im}\,\bar\partial^{\star}\bigg).\end{equation}

\noindent Indeed, we can easily check that the middle sum on the r.h.s. of (\ref{eqn:image-delta-tilde_inclusion}) is orthogonal. We have $\mbox{Im}\,\bar\partial\perp\mbox{Im}\,\bar\partial^{\star}$ since $\bar\partial^2=0$ and $\mbox{Im}\,\bar\partial\perp\mbox{Im}\,(\partial^{\star}\circ p'')$ since $\langle\langle\bar\partial u,\,\partial^{\star}p''v\rangle\rangle = \langle\langle\partial\bar\partial u,\,p''v\rangle\rangle =0$ for all $u,v$ because $\partial\bar\partial u\in\mbox{Im}\,\bar\partial\perp\ker\Delta''\ni p''v$. Similarly, $\mbox{Im}\,(\partial_{|\ker\bar\partial})\perp\mbox{Im}\,(\partial^{\star}\circ p'')$ since $\partial^2=0$ and $\mbox{Im}\,(\partial_{|\ker\bar\partial})\perp\mbox{Im}\,\bar\partial^{\star}$ since $\langle\langle\partial u,\,\bar\partial^{\star}v\rangle\rangle = \langle\langle\bar\partial\partial u,\,v\rangle\rangle = 0$ for all $u\in\ker\bar\partial$ and all $v$.

 Now, putting together (\ref{eqn:2-space-decomp_Delta-tilde}) and (\ref{eqn:image-delta-tilde_inclusion}), we get

$$C^{\infty}_{p,\,q}(X,\,\C) \subset \ker\widetilde{\Delta}\bigoplus\bigg(\mbox{Im}\,\bar\partial + \mbox{Im}\,(\partial_{|\ker\bar\partial})\bigg)\bigoplus\bigg(\mbox{Im}\,(\partial^{\star}\circ p'') +  \mbox{Im}\,\bar\partial^{\star}\bigg)$$

\noindent in which the inclusion must be an equality because all the three mutually orthogonal spaces on the r.h.s. are contained in $C^{\infty}_{p,\,q}(X,\,\C)$. This proves (\ref{eqn:3-space-decomp_Delta-tilde}) and also that the inclusion in (\ref{eqn:image-delta-tilde_inclusion}) is an equality.

 The first of the three $2$-space decompositions stated after (\ref{eqn:3-space-decomp_Delta-tilde}) will be proved as (\ref{eqn:orthogonal_decomp_tilde}) in the proof of the next Theorem \ref{The:Hodge-isom_E2}, while the second one can be proved analogously. The third one is (\ref{eqn:image-delta-tilde_inclusion}) that was seen above to be an equality.

 The last two statements about the Green's operator and the spectrum are proved in the usual way using the elliptic theory.  \hfill $\Box$

\vspace{2ex}

We can now state our Hodge isomorphism for the spaces $E_2^{p,\,q}$ in the Fr\"olicher spectral sequence.

\begin{The}\label{The:Hodge-isom_E2} Let $(X,\,\omega)$ be a compact Hermitian manifold with $\mbox{dim}_{\C}X=n$. For every $p, q\in\{0, 1, \dots , n\}$, let $\widetilde{\cal H}^{p,\,q}_{\widetilde\Delta}(X,\,\C)$ stand for the kernel of $\widetilde\Delta$ acting on $(p,\,q)$-forms. Then the map

\begin{equation}\label{eqn:S_def}S=S^{p,\,q} : \widetilde{\cal H}^{p,\,q}_{\widetilde\Delta}(X,\,\C) \longrightarrow \widetilde{H}^{p,\,q}(X,\,\C), \hspace{3ex} \alpha\longmapsto \widetilde{[\alpha]},\end{equation}

\noindent is an {\bf isomorphism}. In particular, its composition with the isomorphism $T:\widetilde{H}^{p,\,q} (X,\,\C)\longrightarrow E_2^{p,\,q}$ defined in (\ref{eqn:T_def}) yields the {\bf Hodge isomorphism}

\begin{equation}\label{eqn:Hodge-isom_E2}T\circ S = T^{p,\,q}\circ S^{p,\,q} : \widetilde{\cal H}^{p,\,q}_{\widetilde\Delta}(X,\,\C) \longrightarrow  E_2^{p,\,q}, \hspace{3ex} \alpha\longmapsto \bigg[[\alpha]_{\bar\partial}\bigg]_{d_1}.\end{equation}

\noindent Thus, every class $[[\alpha]_{\bar\partial}]_{d_1}\in E_2^{p,\,q}$ contains a unique $\widetilde\Delta$-harmonic representative $\alpha$.

\end{The}

\noindent {\it Proof.} Thanks to (\ref{eqn:ker_Delta-tilde}), we have

\begin{equation}\label{eqn:Delta-tilde-harmonic-space}\widetilde{\cal H}^{p,\,q}_{\widetilde\Delta}(X,\,\C) = \ker(p''\circ\partial)\cap\ker(p''\circ\partial^{\star})\cap\ker\bar\partial\cap\ker\bar\partial^{\star} \subset \ker(p''\circ\partial)\cap\ker\bar\partial.\end{equation}

\noindent In particular, every form $\alpha\in\widetilde{\cal H}^{p,\,q}_{\widetilde\Delta}(X,\,\C)$ defines a class $\widetilde{[\alpha]}\in \widetilde{H}^{p,\,q} (X,\,\C)$, so the map $S^{p,\,q}$ is well defined. We now prove the following orthogonal decomposition

\begin{equation}\label{eqn:orthogonal_decomp_tilde}\ker(p''\circ\partial)\cap\ker\bar\partial = \ker\widetilde{\Delta}\bigoplus\bigg(\mbox{Im}\,\bar\partial + \mbox{Im}\,(\partial_{|\ker\bar\partial})\bigg),\end{equation}

\noindent where $\ker\widetilde{\Delta} = \widetilde{\cal H}^{p,\,q}_{\widetilde\Delta}(X,\,\C)$ is given by (\ref{eqn:Delta-tilde-harmonic-space}). It is clear that (\ref{eqn:orthogonal_decomp_tilde}) implies that $S$ is an isomorphism. 

 Thanks to the $3$-space orthogonal decomposition (\ref{eqn:3-space-decomp_Delta-tilde}), proving (\ref{eqn:orthogonal_decomp_tilde}) is equivalent to proving

\begin{equation}\label{eqn:orthogonal_decomp_tilde1}\ker(p''\circ\partial)\cap\ker\bar\partial = \bigg(\mbox{Im}\,(\partial^{\star}\circ p'') + \mbox{Im}\,\bar\partial^{\star}\bigg)^{\perp}.\end{equation}

\noindent Now, the r.h.s. term in (\ref{eqn:orthogonal_decomp_tilde1}) is the intersection of $(\mbox{Im}\,(\partial^{\star}\circ p''))^{\perp} = \ker\,(\partial^{\star}\circ p'')^{\star} = \ker\,(p''\circ\partial)$ with $(\mbox{Im}\,\bar\partial^{\star})^{\perp} = \ker\bar\partial$. This proves (\ref{eqn:orthogonal_decomp_tilde1}), hence also (\ref{eqn:orthogonal_decomp_tilde}).   \hfill $\Box$.

\section{Harmonic metrics for the pseudo-differential Laplacian}\label{section:harmonic-metrics}

 It is well known that only K\"ahler metrics are harmonic for the most commonly used Laplace-type operators they induce. We shall now point out which Hermitian metrics lie in the kernel of $\Delta'_{p''}$.

\begin{Prop}\label{Prop:harmonic-metrics} Let $(X,\,\omega)$ be a compact Hermitian manifold of dimension $n$. Consider the operator $\Delta'_{p''}$ induced by $\omega$.

\noindent $(i)$\, The following equivalences hold: \begin{eqnarray}\label{eqn:harmonic-metrics_11}\Delta'_{p''}\omega=0 & \iff & \partial\omega\in\mbox{Im}\,\bar\partial \oplus \mbox{Im}\,\bar\partial^{\star} \hspace{1ex} \mbox{and} \hspace{1ex} \partial^{\star}\omega\in\mbox{Im}\,\bar\partial \oplus \mbox{Im}\,\bar\partial^{\star},\\
\label{eqn:harmonic-metrics_n-1n-1}\Delta'_{p''}\omega^{n-1}=0 & \iff & \partial\omega^{n-1}\in\mbox{Im}\,\bar\partial \oplus \mbox{Im}\,\bar\partial^{\star} \hspace{1ex} \mbox{and} \hspace{1ex} \partial^{\star}\omega^{n-1}\in\mbox{Im}\,\bar\partial \oplus \mbox{Im}\,\bar\partial^{\star}.\end{eqnarray}

\noindent $(ii)$\, In particular, if $\partial\omega\in\mbox{Im}\,\bar\partial$ and $\partial\omega^{n-1}\in\mbox{Im}\,\bar\partial$, then $\Delta'_{p''}\omega=0$ and $\Delta'_{p''}\omega^{n-1}=0$.

\noindent $(iii)$\, If $\omega$ is both SKT (i.e.\!\! $\partial\bar\partial\omega=0$) and Gauduchon (i.e.\!\! $\partial\bar\partial\omega^{n-1}=0$), the next equivalences hold:

\begin{equation}\label{eqn:SKT-Gauduchon_equiv}\Delta'_{p''}\omega=0 \iff  \Delta'_{p''}\omega^{n-1}=0 \iff \partial\omega\in\mbox{Im}\,\bar\partial \hspace{1ex} \mbox{and} \hspace{1ex} \partial\omega^{n-1}\in\mbox{Im}\,\bar\partial.\end{equation}

\end{Prop}

 We shall term the metrics $\omega$ with the property $\partial\omega\in\mbox{Im}\,\bar\partial$ {\it super SKT}, while those satisfying $ \partial\omega^{n-1}\in\mbox{Im}\,\bar\partial$ are the strongly Gauduchon (sG) metrics of [Pop13].

\vspace{2ex}

\noindent {\it Proof.} $(i)$\, The condition $\Delta'_{p''}\omega=0$ is equivalent to $p''(\partial\omega) = 0$ and $p''(\partial^{\star}\omega) = 0$ (cf. (\ref{eqn:ker_Delta-tilde})) which amount to $\partial\omega\perp\ker\Delta''$ and $\partial^{\star}\omega\perp\ker\Delta''$. These conditions are, in turn, equivalent to the conditions stated on the r.h.s. of (\ref{eqn:harmonic-metrics_11}) thanks to the orthogonal $3$-space decomposition $C^{\infty}_{p,\,q}(X,\,\C) = \ker\Delta''\oplus\mbox{Im}\,\bar\partial\oplus\mbox{Im}\,\bar\partial^{\star}$ (cf. (\ref{eqn:3-space-decomp})) for $(p,\,q)=(2,\,1)$, resp. $(p,\,q)=(0,\,1)$. This proves (\ref{eqn:harmonic-metrics_11}) and (\ref{eqn:harmonic-metrics_n-1n-1}) is proved similarly.

 $(ii)$\, It suffices to show that the sG condition $\partial\omega^{n-1}\in\mbox{Im}\,\bar\partial$ is equivalent to $\partial^{\star}\omega\in\mbox{Im}\,\bar\partial^{\star}$ and that the super SKT condition $\partial\omega\in\mbox{Im}\,\bar\partial$ is equivalent to $\partial^{\star}\omega^{n-1}\in\mbox{Im}\,\bar\partial^{\star}$. This follows immediately from the Hodge star operator $\star = \star_{\omega} : \Lambda^{p,\,q}T^{\star}X\longrightarrow\Lambda^{n-q,\,n-p}T^{\star}X$ being an isomorphism and from the well-known formulae

$$\partial^{\star} = -\star\bar\partial\star, \hspace{2ex} \bar\partial^{\star} = -\star\partial\star,   \hspace{2ex}  \star\star = \pm\, 1   \hspace{2ex} \mbox{and} \hspace{2ex} \star\omega = \frac{\omega^{n-1}}{(n-1)!}.$$  

\noindent For example, we have the equivalences:

\vspace{1ex}

\hspace{3ex} $\partial^{\star}\omega\in\mbox{Im}\,\bar\partial^{\star} \iff -\star\bar\partial\star\omega \in\mbox{Im}\,(-\star\partial\star) \iff \bar\partial\omega^{n-1}\in\mbox{Im}\,\partial \iff \partial\omega^{n-1}\in\mbox{Im}\,\bar\partial,$

\vspace{1ex}

\noindent where the last equivalence follows by conjugation.

 $(iii)$\, The SKT condition $\partial\bar\partial\omega=0$ is equivalent to $\partial\omega\in\ker\bar\partial$ (and also to $\partial^{\star}\omega^{n-1}\in\ker\bar\partial^{\star}$), while the Gauduchon condition $\partial\bar\partial\omega^{n-1}=0$ is equivalent to $\partial\omega^{n-1}\in\ker\bar\partial$ (and also to $\partial^{\star}\omega\in\ker\bar\partial^{\star}$). Thus, the equivalences (\ref{eqn:SKT-Gauduchon_equiv}) follow from $(i)$, $(ii)$ and from the orthogonal splittings $\ker\bar\partial = \ker\Delta''\oplus\mbox{Im}\,\bar\partial$ and $\ker\bar\partial^{\star} = \ker\Delta''\oplus\mbox{Im}\,\bar\partial^{\star}$ (cf. (\ref{eqn:ker-del_ker_laplacian})).  \hfill $\Box$

\vspace{3ex}

 We do not know at this point whether non-K\"ahler Hermitian metrics that are both super SKT and strongly Gauduchon exist. The author is grateful to L. Ugarte for informing him that 

\vspace{1ex}

 $(1)$\, no such metrics exist on nilmanifolds of complex dimension $3$;

 $(2)$\, on nilmanifolds of complex dimension $3$, there exist non-K\"ahler Hermitian metrics that are both SKT and strongly Gauduchon.

\section{Sufficient metric conditions for the $E_2$ degeneration of the Fr\"olicher spectral sequence}\label{section:sufficient-Froelicher}

 Throughout this section, $(X,\,\omega)$ will be a compact Hermitian manifold with $\mbox{dim}_{\C}X=n$. Recall that for every $k\in\{0,\dots , 2n\}$, the $d$-Laplacian $\Delta : C^{\infty}_k(X,\,\C)\longrightarrow C^{\infty}_k(X,\,\C)$ is defined by $\Delta = dd^{\star} + d^{\star}d$. If we denote by ${\cal H}_{\Delta}^k(X,\,\C)\subset C^{\infty}_k(X,\,\C)$ the kernel of $\Delta$ acting on smooth forms of degree $k$, we have the Hodge isomorphism ${\cal H}_{\Delta}^k(X,\,\C)\simeq H^k_{DR}(X,\,\C)$ with the De Rham cohomology group of degree $k$.

 We start with the following very simple observation.

\begin{Lem}\label{Lem:sufficient-cond_E2-deg} $(a)$\,\, If for every $p, q\in\{0, 1, \dots , n\}$ the following map induced by the identity

\begin{equation}\label{eqn:Jpq_def}J^{p,\,q} : \widetilde{\cal H}^{p,\,q}_{\widetilde\Delta}(X,\,\C) \longrightarrow {\cal H}_{\Delta}^{p+q}(X,\,\C), \hspace{3ex} \gamma\longmapsto\gamma,\end{equation}

\noindent is well defined, then the Fr\"olicher spectral sequence of $X$ degenerates at $E_2$.

 $(b)$\, A sufficient condition for the map $J^{p,\,q}$ to be well defined is that the following inequality hold

\begin{equation}\label{eqn:Jpq_well-definedness_suff}\Delta'-\Delta'_{p''}\leq \Delta'' + (C \Delta'' + (1-\varepsilon) \Delta')  \hspace{3ex} \mbox{on}\hspace{1ex}(p,\,q)\mbox{-forms},\end{equation}

\noindent for some constants $C\geq 0$ and $0<\varepsilon \leq 1$ depending only on $X$, $\omega$ and $(p,\,q)$. (Recall that $\Delta'-\Delta'_{p''} = \Delta'_{p''_{\perp}}\geq 0$.)

 Thus, (\ref{eqn:Jpq_well-definedness_suff}) implies the degeneracy at $E_2$ of the Fr\"olicher spectral sequence of $X$.

\end{Lem}

\noindent {\it Proof.} $(a)$ Well-definedness for $J^{p,\,q}$ means that for every smooth $(p,\,q)$-form $\gamma$ we have $\Delta\gamma=0$ whenever $\widetilde{\Delta}\gamma = 0$. It is clear that $J^{p,\,q}$ is automatically injective if it is well defined, hence in that case $\mbox{dim}\,\widetilde{\cal H}^{p,\,q}_{\widetilde\Delta}(X,\,\C)\leq\mbox{dim}\,{\cal H}_{\Delta}^{p+q}(X,\,\C)$. Therefore, if all the maps $J^{p,\,q}$ are well defined, then

\begin{equation}\label{eqn:dim-E2pqleq-b_k}\sum\limits_{p+q=k}\mbox{dim}\,E_2^{p,\,q}\leq b_k:=\mbox{dim}\,H^k_{DR}(X,\,\C) \hspace{2ex} \mbox{for all}\hspace{1ex} k\in\{0,\dots , 2n\}\end{equation} 

\noindent since $\mbox{dim}\,E_2^{p,\,q} = \mbox{dim}\,\widetilde{\cal H}^{p,\,q}_{\widetilde\Delta}(X,\,\C)$ by the Hodge isomorphism (\ref{eqn:Hodge-isom_E2}) and the images $J^{p,\,q}(\widetilde{\cal H}^{p,\,q}_{\widetilde\Delta}(X,\,\C))$ in ${\cal H}_{\Delta}^k(X,\,\C)$ have pairwise intersections reduced to zero for $p+q=k$ for bidegree reasons. Inequality (\ref{eqn:dim-E2pqleq-b_k}) is precisely the degeneracy condition (\ref{eqn:E_r_degeneracy-condition}) at $E_2$.

$(b)$ Clearly, a sufficient condition for $J^{p,\,q}$ to be well defined is that the following inequality hold

\begin{equation}\label{eqn:Delta_Delta-tilde_ineq}\langle\langle\Delta\gamma,\,\gamma\rangle\rangle \leq C\,\langle\langle\widetilde\Delta\gamma,\,\gamma\rangle\rangle \hspace{3ex}\mbox{for all}\hspace{1ex}\gamma\in C^{\infty}_{p,\,q}(X,\,\C)\end{equation}

\noindent since $ \Delta, \widetilde\Delta\geq 0$. Now, by definition of $\widetilde\Delta$ (cf. (\ref{eqn:Delta-tilde_def_bis})), $\langle\langle\widetilde\Delta\gamma,\,\gamma\rangle\rangle = \langle\langle\Delta'_{p''}\gamma,\,\gamma\rangle\rangle + \langle\langle\Delta''\gamma,\,\gamma\rangle\rangle$. Meanwhile, for every $(p,\,q)$-form $\gamma$, we have $\langle\langle\Delta\gamma,\,\gamma\rangle\rangle = ||\partial\gamma + \bar\partial\gamma||^2 + ||\partial^{\star}\gamma + \bar\partial^{\star}\gamma||^2 = ||\partial\gamma||^2 + ||\partial^{\star}\gamma||^2 + ||\bar\partial\gamma||^2 + ||\bar\partial^{\star}\gamma||^2 = \langle\langle\Delta'\gamma,\,\gamma\rangle\rangle + \langle\langle\Delta''\gamma,\,\gamma\rangle\rangle$ since $\partial\gamma$ is orthogonal to $\bar\partial\gamma$ and $\partial^{\star}\gamma$ is orthogonal to $\bar\partial^{\star}\gamma$ for bidegree reasons. (This argument breaks down if $\gamma$ is not of pure type.) Thus

 \begin{equation}\label{eqn:Delta_gamma_gamma_pre-type}\langle\langle\Delta\gamma,\,\gamma\rangle\rangle = \langle\langle\Delta'\gamma,\,\gamma\rangle\rangle + \langle\langle\Delta''\gamma,\,\gamma\rangle\rangle   \hspace{3ex}\mbox{for all}\hspace{1ex}\gamma\in C^{\infty}_{p,\,q}(X,\,\C).\end{equation}

\noindent It is now clear that (\ref{eqn:Jpq_well-definedness_suff}) implies (\ref{eqn:Delta_Delta-tilde_ineq}) with a possibly different constant $C$, so (\ref{eqn:Jpq_well-definedness_suff}) implies the well-definedness of $J^{p,\,q}$.  \hfill $\Box$

\vspace{2ex}

 Concerning inequality (\ref{eqn:Jpq_well-definedness_suff}), note that the stronger inequality $\langle\langle\Delta'\gamma,\,\gamma\rangle\rangle \leq C\,\langle\langle\Delta''\gamma,\,\gamma\rangle\rangle$ for all $(p,\,q)$-forms $\gamma$ and all bidegrees $(p,\,q)$ implies the degeneracy at $E_1$ of the Fr\"olicher spectral sequence, but we shall not pursue this here.

\subsection{Use of $(b)$ of Lemma \ref{Lem:sufficient-cond_E2-deg}}\label{subsection:use-of-b}

 We shall now concentrate on proving inequality (\ref{eqn:Jpq_well-definedness_suff}) under the SKT assumption (i.e. $\partial\bar\partial\omega = 0$) coupled with a torsion assumption on the metric $\omega$.

\begin{Lem}\label{Lem:reduction-ineq-starstar} A sufficient condition for (\ref{eqn:Jpq_well-definedness_suff}) to hold (hence for $E_2(X)=E_{\infty}(X)$) is that there exist constants $0<\delta<1-\varepsilon<1$ and $C\geq 0$ such that the following inequality holds

\begin{eqnarray}\label{eqn:reduction-ineq-starstar}\nonumber  (1-\varepsilon-\delta)\,\bigg(||p''_{\perp}\partial u||^2 + ||p''_{\perp}\partial^{\star} u||^2\bigg) + (1-\varepsilon)\,(||p''\partial u||^2 + ||p''\partial^{\star} u||^2) + C\,\langle\langle\Delta''u,\,u\rangle\rangle \geq & & \\
     \bigg(\frac{1}{\delta} - 1\bigg)\,\bigg(||p''_{\perp}\tau u||^2 + ||p''_{\perp}\tau^{\star} u||^2\bigg) + \langle\langle[\partial\omega\wedge\cdot,\,(\partial\omega\wedge\cdot)^{\star}]\,u,\,u\rangle\rangle - \langle\langle[\Lambda,\,[\Lambda,\,\frac{i}{2}\,\partial\bar\partial\omega]]\,u,\,u\rangle\rangle  & & \end{eqnarray}

\noindent for every form $u\in C^{\infty}_{p,\,q}(X,\,\C)$ and every bidegree $(p,\,q)$. (Note that all the terms on the r.h.s. of (\ref{eqn:reduction-ineq-starstar}) are of order zero, hence bounded, while the last and only signless term vanishes if $\omega$ is SKT.)

\end{Lem}

\noindent {\it Proof.} By Demailly's non-K\"ahler Bochner-Kodaira-Nakano identity $\Delta'' = \Delta'_{\tau} + T_{\omega}$ (cf. (\ref{eqn:BKN_demailly2})), inequality (\ref{eqn:Jpq_well-definedness_suff}) is equivalent to each of the following inequalities:

\begin{eqnarray}\label{eqn:Jpq_well-definedness_suff1}\nonumber \Delta' - \Delta'_{p''} & \leq & \Delta' + [\tau,\,\partial^{\star}] + [\partial,\,\tau^{\star}] + [\tau,\,\tau^{\star}] + C\,\Delta'' + (1-\varepsilon)\,\Delta' + T_{\omega} \iff \\
     0 & \leq & \bigg(\Delta'_{p''} + (\tau p''\partial^{\star} + \partial^{\star}p''\tau) + (\partial p''\tau^{\star} + \tau^{\star}p''\partial) + (\tau p''\tau^{\star} + \tau^{\star}p''\tau)\bigg) \\
\nonumber & + & (1-\varepsilon)\,\Delta' + C\,\Delta'' + (\tau p''_{\perp}\partial^{\star} + \partial^{\star}p''_{\perp}\tau) + (\partial p''_{\perp}\tau^{\star} + \tau^{\star}p''_{\perp}\partial) + (\tau p''_{\perp}\tau^{\star} + \tau^{\star}p''_{\perp}\tau) + T_{\omega}.\end{eqnarray}

\noindent Since $\Delta'_{p''} + (\tau p''\partial^{\star} + \partial^{\star}p''\tau) + (\partial p''\tau^{\star} + \tau^{\star}p''\partial) + (\tau p''\tau^{\star} + \tau^{\star}p''\tau) = (\partial + \tau)p''(\partial^{\star} + \tau^{\star}) + (\partial^{\star} + \tau^{\star})p''(\partial + \tau) \geq 0$, inequality (\ref{eqn:Jpq_well-definedness_suff1}) holds if the following inequality holds

\begin{eqnarray}\label{eqn:Jpq_well-definedness_suff2}\nonumber (1-\varepsilon)\,\langle\langle\Delta'u,\,u\rangle\rangle & + & C\,\langle\langle\Delta''u,\,u\rangle\rangle + ||p''_{\perp}\tau u||^2 + ||p''_{\perp}\tau^{\star}u||^2 \\
 & \geq & - 2\,\mbox{Re}\,\langle\langle p''_{\perp}\partial^{\star}u,\,p''_{\perp}\tau^{\star}u\rangle\rangle - 2\,\mbox{Re}\,\langle\langle p''_{\perp}\partial u,\,p''_{\perp}\tau u\rangle\rangle - \langle\langle T_{\omega}u,\,u\rangle\rangle.\end{eqnarray}

\noindent Now, suppose that $0<\varepsilon<1$ and choose any $0<\delta<1-\varepsilon$. The Cauchy-Schwarz inequality gives

\begin{equation}\label{eqn:C-S}\nonumber\bigg| 2\,\mbox{Re}\,\langle\langle p''_{\perp}\partial u,\,p''_{\perp}\tau u\rangle\rangle\bigg| \leq \delta\,||p''_{\perp}\partial u||^2 + \frac{1}{\delta}\,||p''_{\perp}\tau u||^2,\, \bigg| 2\,\mbox{Re}\,\langle\langle p''_{\perp}\partial^{\star} u,\,p''_{\perp}\tau^{\star} u\rangle\rangle\bigg| \leq \delta\,||p''_{\perp}\partial^{\star} u||^2 + \frac{1}{\delta}\,||p''_{\perp}\tau^{\star} u||^2.\end{equation}

\noindent Thus, for (\ref{eqn:Jpq_well-definedness_suff2}) to hold, it suffices that the following inequality hold:

\begin{eqnarray}\label{eqn:Jpq_well-definedness_suff3}\nonumber (1-\varepsilon)\,\langle\langle\Delta'u,\,u\rangle\rangle + C\,\langle\langle\Delta''u,\,u\rangle\rangle & \geq & \delta\,(||p''_{\perp}\partial u||^2 + ||p''_{\perp}\partial^{\star} u||^2) + \bigg(\frac{1}{\delta} - 1\bigg)\,(||p''_{\perp}\tau u||^2 + ||p''_{\perp}\tau^{\star} u||^2) \\
   & + & \langle\langle[\partial\omega\wedge\cdot,\,(\partial\omega\wedge\cdot)^{\star}]\,u,\,u\rangle\rangle - \langle\langle[\Lambda,\,[\Lambda,\,\frac{i}{2}\,\partial\bar\partial\omega]]\,u,\,u\rangle\rangle.\end{eqnarray}

\noindent This is equivalent to (\ref{eqn:reduction-ineq-starstar}) since 

\vspace{1ex}

\noindent $\langle\langle\Delta'u,\,u\rangle\rangle = ||p''\partial u + p''_{\perp}\partial u||^2 + ||p''\partial^{\star} u + p''_{\perp}\partial^{\star} u||^2 = (||p''\partial u||^2 + ||p''\partial^{\star} u||^2) + (||p''_{\perp}\partial u||^2 + ||p''_{\perp}\partial^{\star} u||^2)$

\vspace{1ex}

\noindent thanks to the obvious orthogonality relations $p''\partial u \perp p''_{\perp}\partial u$ and $p''\partial^{\star} u \perp p''_{\perp}\partial^{\star} u$. \hfill $\Box$

\vspace{3ex}

 To apply Lemma \ref{Lem:reduction-ineq-starstar}, we start with a very simple elementary observation.

\begin{Lem}\label{Lem:Hilbert-space-obs} Let ${\cal H}$ be a Hilbert space and let $A,B:{\cal H}\rightarrow{\cal H}$ be closed linear operators such that $A, B\geq 0$, $A=A^{\star}$ and $B=B^{\star}$.

 If $\ker A\subset\ker B$ and if $B\leq A$ on $(\ker A)^{\perp}$, then $B\leq A$.

\end{Lem}

\noindent {\it Proof.} We have to prove that $\langle Bu,\,u\rangle \leq \langle Au,\,u\rangle$ for all $u$. Since $A$ is closed, $\ker A$ is closed in ${\cal H}$, so every $u\in{\cal H}$ splits uniquely as $u = u_A + u_A^{\perp}$ with $u_A\in\ker A$ and $u_A^{\perp}\in(\ker A)^{\perp}$. Moreover,

\begin{equation}\label{eqn:A_kerA-orthogonal_inclusion}A((\ker A)^{\perp}) \subset (\ker A)^{\perp}.\end{equation}

\noindent Indeed, for every $u_A^{\perp}\in(\ker A)^{\perp}$ and every $v\in\ker A$, we have: $\langle A(u_A^{\perp}),\,v\rangle = \langle u_A^{\perp},\,Av\rangle = 0$ since $A^{\star}v = Av = 0$. Therefore, for every $u$, we get:

\vspace{1ex}

\hspace{6ex} $\langle Au,\,u\rangle = \langle Au_A^{\perp},\,u_A + u_A^{\perp}\rangle = \langle Au_A^{\perp},\, u_A^{\perp}\rangle \geq \langle Bu_A^{\perp},\, u_A^{\perp}\rangle = \langle Bu_A^{\perp},\, u\rangle = \langle Bu,\, u\rangle.$

\vspace{1ex}

\noindent The second identity above followed from (\ref{eqn:A_kerA-orthogonal_inclusion}), the inequality followed from the hypothesis and the last two identities followed from the next relations:

\vspace{1ex}

\hspace{20ex} $(i)\,\, B((\ker A)^{\perp}) \subset (\ker A)^{\perp} \hspace{2ex} \mbox{and} \hspace{2ex} (ii)\,\, B(\ker A)=0.$

\noindent To prove $(i)$, let $u_A^{\perp}\in(\ker A)^{\perp}$ and $v\in\ker A\subset\ker B$. We have: $\langle B(u_A^{\perp}),\,v\rangle = \langle u_A^{\perp},\,Bv\rangle = 0$ since $B^{\star}v = Bv = 0$. Identity $(ii)$ follows from the hypothesis $\ker A\subset\ker B$.   \hfill $\Box$

\vspace{2ex}

 We shall now apply Lemma \ref{Lem:Hilbert-space-obs} to the non-negative self-adjoint operators 

$$B:=\Delta'_{p''_{\perp}} = \Delta'- \Delta'_{p''}\geq 0  \hspace{2ex} \mbox{and} \hspace{2ex} A:= (C+1)\,\Delta'' + (1-\varepsilon)\,\Delta'\geq 0$$

\noindent for which we obviously have $\ker B = \ker\Delta'_{p''_{\perp}} \supset \ker\Delta' \supset \ker\Delta'\cap\ker\Delta'' = \ker A$. The choice of constants
 $C>0$ and $0<\varepsilon<1$ will be specified later on.

 We know from $(b)$ of Lemma \ref{Lem:sufficient-cond_E2-deg} that a sufficient condition for $E_2(X)=E_{\infty}(X)$ in the Fr\"olicher spectral sequence is the validity of inequality (\ref{eqn:Jpq_well-definedness_suff}), i.e. of the inequality $B\leq A$. By Lemma \ref{Lem:Hilbert-space-obs}, this is equivalent to having $B\leq A$ on $(\ker A)^{\perp} = (\ker\Delta'\cap\ker\Delta'')^{\perp}$. Now, the proof of Lemma \ref{Lem:reduction-ineq-starstar} shows that for this to hold, it suffices for the inequality (\ref{eqn:Jpq_well-definedness_suff3}) to hold on $(\ker\Delta'\cap\ker\Delta'')^{\perp}$. If we assume $\partial\bar\partial\omega = 0$, after bounding above $||p''_{\perp}v||$ by $||v||$ for $v\in\{\partial u, \partial^{\star}u, \tau u, \tau^{\star}u\}$ in the r.h.s. of (\ref{eqn:Jpq_well-definedness_suff3}), we see that it suffices to have

\begin{eqnarray}\label{eqn:E2_suff_SKT1} (1-\varepsilon-\delta)\,\langle\langle\Delta'u,\,u\rangle\rangle + C\,\langle\langle\Delta''u,\,u\rangle\rangle \geq \bigg(\frac{1}{\delta} - 1\bigg)\,\langle\langle[\tau,\,\tau^{\star}]\,u,\,u\rangle\rangle +  \langle\langle[\partial\omega\wedge\cdot,\,(\partial\omega\wedge\cdot)^{\star}]\,u,\,u\rangle\rangle .\end{eqnarray}

\noindent for all $u\in (\ker\Delta'\cap\ker\Delta'')^{\perp}$ and some fixed constants $C>0$, $0<\delta < 1-\varepsilon < 1$.

 Now, we choose the constants such that $\delta=1-2\varepsilon >0$ (so $0<\varepsilon < \frac{1}{2}$) and $C=1-\varepsilon-\delta = \varepsilon$. Thus, $(1/\delta)-1 = 2\varepsilon/(1-2\varepsilon)$. If, moreover, we choose $\varepsilon$ such that $2/(1-2\varepsilon) < 3$ (i.e. such that $0<\varepsilon<1/6$), (\ref{eqn:E2_suff_SKT1}) holds with these choices of constants whenever the following inequality holds:

\begin{eqnarray}\label{eqn:E2_suff_SKT2} \langle\langle(\Delta' + \Delta'')\,u,\,u\rangle\rangle \geq 3\, \langle\langle([\tau,\,\tau^{\star}] + [\partial\omega\wedge\cdot,\,(\partial\omega\wedge\cdot)^{\star}])\,u,\,u\rangle\rangle \hspace{2ex} \mbox{for all}\hspace{1ex} u\in(\ker\Delta'\cap\ker\Delta'')^{\perp}.\end{eqnarray}

 For all $p,q\in\{0,\dots , n\}$, the non-negative self-adjoint differential operator $\Delta' + \Delta'' : C^{\infty}_{p,\,q}(X,\,\C)\longrightarrow C^{\infty}_{p,\,q}(X,\,\C)$ is elliptic. Therefore, since $X$ is compact, it has a discrete spectrum contained in $[0,\, +\infty)$ with $+\infty$ as its only accumulation point. In particular, it has a smallest positive eigenvalue that we denote by

\begin{equation}\label{eqn:lambda_0mu_0_def}\rho_{\omega}^{p,\,q}:=\mbox{min}\,\bigg(\mbox{Spec}\,(\Delta' + \Delta'')^{p,\,q}\cap(0,\,+\infty)\bigg)>0.\end{equation}

\noindent Thus, $\rho_{\omega}^{p,\,q}$ is the size of the spectral gap of $\Delta' + \Delta''$ acting on $(p,\,q)$-forms. We get

\begin{equation}\label{eqn:lower-bound_Delta'+Delta''_orthogonal}\langle\langle(\Delta' + \Delta'')\,u,\,u\rangle\rangle \geq \rho_{\omega}^{p,\,q}\,||u||^2   \hspace{2ex} \mbox{for all}\hspace{1ex} u\in C^{\infty}_{p,\,q}(X,\,\C)\cap(\ker\Delta'\cap\ker\Delta'')^{\perp},\end{equation}

\noindent since $\ker(\Delta' + \Delta'') = \ker\Delta'\cap\ker\Delta''$. On the other hand, the non-negative torsion operator $[\tau,\,\tau^{\star}] + [\partial\omega\wedge\cdot,\,(\partial\omega\wedge\cdot)^{\star}]$ is of order zero, hence bounded, hence

\begin{equation}\label{eqn:nonnegative-torsion_upper-bound} \langle\langle([\tau,\,\tau^{\star}] + [\partial\omega\wedge\cdot,\,(\partial\omega\wedge\cdot)^{\star}])\,u,\,u\rangle\rangle \leq C_{\omega}^{p,\,q}\,||u||^2 \hspace{2ex} \mbox{for all}\hspace{1ex} u\in C^{\infty}_{p,\,q}(X,\,\C),\end{equation}

\noindent where $C_{\omega}^{p,\,q}:= \sup\limits_{u\in C^{\infty}_{p,\,q}(X,\,\C),\,\,||u||=1}\langle\langle([\tau,\,\tau^{\star}] + [\partial\omega\wedge\cdot,\,(\partial\omega\wedge\cdot)^{\star}])\,u,\,u\rangle\rangle$.

 We conclude from (\ref{eqn:lower-bound_Delta'+Delta''_orthogonal}) and (\ref{eqn:nonnegative-torsion_upper-bound}) that (\ref{eqn:E2_suff_SKT2}) holds if $\rho_{\omega}^{p,\,q} \geq 3\,C_{\omega}^{p,\,q}$. We have thus proved the following statement which is nothing but Theorem \ref{The:E_2_small-torsionSKT_introd}.

\begin{The}\label{The:E_2_small-torsionSKT} Let $X$ be a compact complex $n$-dimensional manifold. If $X$ carries an SKT metric $\omega$ whose torsion satisfies the condition

\begin{equation}\label{eqn:torsion-smallness_spectral-gap}C_{\omega}^{p,\,q} \leq \frac{1}{3}\,\rho_{\omega}^{p,\,q}\end{equation}

\noindent for all $p,q\in\{0,\dots , n\}$, then the Fr\"olicher spectral sequence of $X$ degenerates at $E_2$.

\end{The}

\subsection{Use of $(a)$ of Lemma \ref{Lem:sufficient-cond_E2-deg}}\label{subsection:use-of-a}

 We shall now give a different kind of metric condition ensuring that $E_2(X) = E_{\infty}(X)$ in the Fr\"olicher spectral sequence. 

\begin{Lem}\label{Lem:kernels-intersection_E2} Let $X$ be a compact complex manifold with $\mbox{dim}_{\C}X=n$. If $X$ admits a Hermitian metric $\omega$ whose induced operators $\Delta', \Delta'', \Delta'_{p''} : C^{\infty}_{p,\,q}(X,\,\C)\longrightarrow C^{\infty}_{p,\,q}(X,\,\C)$ satisfy the condition

\begin{equation}\label{eqn:kernels-intersection_E2}\ker\Delta'_{p''}\cap\ker\Delta''\subset\ker\Delta'  \hspace{2ex}\mbox{in every bidegree}\hspace{1ex}(p,\,q),\end{equation}

\noindent the Fr\"olicher spectral sequence of $X$ degenerates at $E_2$.

\end{Lem}

\noindent{\it Proof.} As noticed in (\ref{eqn:ker_Delta-tilde}), we always have $\ker\Delta'_{p''}\supset\ker\Delta'$. Recall that $\ker\Delta'_{p''}\cap\ker\Delta'' = \ker\widetilde\Delta$ and that this space is denoted by $\widetilde{\cal H}^{p,\,q}_{\widetilde\Delta}(X,\,\C)$ in bidegree $(p,\,q)$. For every $u\in\widetilde{\cal H}^{p,\,q}_{\widetilde\Delta}(X,\,\C)$, we have $\Delta'u=0$ thanks to (\ref{eqn:kernels-intersection_E2}), hence from (\ref{eqn:Delta_gamma_gamma_pre-type}) we get

$$\langle\langle\Delta u,\,u\rangle\rangle = \langle\langle\Delta'u,\,u\rangle\rangle + \langle\langle\Delta''u,\,u\rangle\rangle = 0 + 0 =0.$$

\noindent This shows that the identity map induces a well-defined linear map $ \widetilde{\cal H}^{p,\,q}_{\widetilde\Delta}(X,\,\C) \longrightarrow {\cal H}_{\Delta}^{p+q}(X,\,\C)$ for all $(p,\,q)$, hence $E_2(X) = E_{\infty}(X)$ by $(a)$ of Lemma \ref{Lem:sufficient-cond_E2-deg}.  \hfill $\Box$

\vspace{2ex}

 We now use Lemma \ref{Lem:kernels-intersection_E2} to give two sufficient metric conditions ensuring that $E_2(X) = E_{\infty}(X)$ in the Fr\"olicher spectral sequence. The following theorem subsumes Theorems \ref{The:E_2-commutation} and \ref{The:E_2-Rbar}.

\begin{The}\label{The:commutation_E2-commutation_Rbar} Let $X$ be a compact complex manifold with $\mbox{dim}_{\C}X=n$. 

\noindent $(i)$\, For any Hermitian metric $\omega$ on $X$, the following three conditions are equivalent:

\vspace{1ex}

$(a)$\, $p''\partial = \partial p''$ on all $(p,\,q)$-forms for all bidegrees $(p,\,q)$;

$(b)$\, $[\partial,\,\bar\partial^{\star}](\ker\Delta'')=0$ and $[\partial,\,\bar\partial^{\star}](\mbox{Im}\,\bar\partial\oplus\mbox{Im}\,\bar\partial^{\star})\subset\mbox{Im}\,\bar\partial\oplus\mbox{Im}\,\bar\partial^{\star}$;

$(c)$\, $[\partial,\,\bar\tau^{\star}](\ker\Delta'')=0$ and $[\partial,\,\bar\tau^{\star}](\mbox{Im}\,\bar\partial\oplus\mbox{Im}\,\bar\partial^{\star})\subset\mbox{Im}\,\bar\partial\oplus\mbox{Im}\,\bar\partial^{\star}$.

\vspace{1ex}

\noindent Moreover, if $X$ carries a Hermitian metric $\omega$ satisfying one of the equivalent conditions $(a),(b), (c)$, the Fr\"olicher spectral sequence of $X$ degenerates at $E_2$.

\vspace{1ex}

\noindent $(ii)$\, If $X$ carries an SKT metric $\omega$ (i.e. such that $\partial\bar\partial\omega=0$) which moreover satisfies the identity

\begin{equation}\label{eqn:tau-bar-S}\langle\langle[\bar\tau,\,\bar\tau^{\star}]\,u,\,u\rangle\rangle = \langle\langle[\bar\partial\omega\wedge\cdot,\,(\bar\partial\omega\wedge\cdot)^{\star}]\,u,\,u\rangle\rangle \hspace{2ex} \mbox{for all}\hspace{1ex} u\in\ker\Delta'_{p''}\cap\ker\Delta'',\end{equation}

\noindent the Fr\"olicher spectral sequence of $X$ degenerates at $E_2$.

\end{The}

\noindent {\it Proof.} $(i)$\, Let $\omega$ be any Hermitian metric on $X$ and let $u$ be any smooth $(p,\,q)$-form. Then $u = u_0 + \bar\partial v + \bar\partial^{\star}w$ with $u_0\in\ker\Delta''$ and $v,w$ smooth forms of bidegrees $(p,\,q-1)$, resp. $(p,\,q+1)$. (Note that we can choose $v\in\mbox{Im}\,\bar\partial^{\star}$ and $w\in\mbox{Im}\,\bar\partial$ if these forms are chosen to have minimal $L^2$ norms.) Thus $p''u=u_0$, so the following equivalences hold:

\begin{eqnarray}\label{eqn:commutation_E2_1}\nonumber p''\partial u = \partial p''u & \iff & p''\partial u_0 + p''\partial\bar\partial v + p''\partial\bar\partial^{\star}w = \partial u_0 \iff p''\partial\bar\partial^{\star}w = p''_{\perp}\partial u_0\\
\nonumber & \iff & p''\partial\bar\partial^{\star}w = 0 \hspace{1ex} \mbox{and} \hspace{1ex} p''_{\perp}\partial u_0 = 0 \iff \partial u_0\in\ker\Delta'' \hspace{1ex} \mbox{and} \hspace{1ex} \partial\bar\partial^{\star}w \in \mbox{Im}\,\bar\partial\oplus\mbox{Im}\,\bar\partial^{\star}\\
        & \iff & \partial u_0\in\ker\bar\partial^{\star} \hspace{1ex} \mbox{and} \hspace{1ex} \partial\bar\partial^{\star}w \in \mbox{Im}\,\bar\partial\oplus\mbox{Im}\,\bar\partial^{\star}.\end{eqnarray}

\noindent We have successively used the following facts: $ p''\partial\bar\partial v = - p''\bar\partial\partial v = 0 $ because $\ker\Delta''\perp\mbox{Im}\,\bar\partial$, $1-p'' = p''_{\perp}$, $\mbox{Im}\,p'' = \ker\Delta'' \perp \mbox{Im}\,\bar\partial\oplus\mbox{Im}\,\bar\partial^{\star} = \mbox{Im}\,p''_{\perp}$, $\partial u_0 \in\ker\bar\partial = \ker\Delta'' \oplus\mbox{Im}\,\bar\partial$ (because $u_0\in\ker\Delta''\subset\ker\bar\partial$), hence the equivalence $\partial u_0\in\ker\Delta''\iff\partial u_0\in\ker\bar\partial^{\star}$.  

 Now, $\bar\partial^{\star}\partial w\in\mbox{Im}\,\bar\partial^{\star}$ and $\partial\bar\partial^{\star}u_0=0$ because $u_0\in\ker\Delta''=\ker\bar\partial\cap\ker\bar\partial^{\star}\subset\ker\bar\partial^{\star}$, so (\ref{eqn:commutation_E2_1}) is equivalent to

$$[\partial,\,\bar\partial^{\star}]\,u_0 =0 \hspace{1ex} \mbox{and} \hspace{1ex} [\partial,\,\bar\partial^{\star}]\,w\in\mbox{Im}\,\bar\partial\oplus\mbox{Im}\,\bar\partial^{\star}.$$

\noindent On the other hand, $[\partial,\,\bar\partial^{\star}] = - [\partial,\,\bar\tau^{\star}]$ by (\ref{eqn:basic-anticommutation}). Since $u_0\in\ker\Delta''$ and $w\in\mbox{Im}\,\bar\partial\oplus\mbox{Im}\,\bar\partial^{\star}$ are arbitrary, the equivalences stated under $(i)$ are proved. 

 To prove the last statement of $(i)$, let $\omega$ be a metric satisfying condition $(a)$. We are going to show that the inclusion (\ref{eqn:kernels-intersection_E2}) holds, hence by Lemma \ref{Lem:kernels-intersection_E2} we shall have $E_2(X)=E_{\infty}(X)$ in the Fr\"olicher spectral sequence of $X$. Let $u\in\ker\Delta'_{p''}\cap\ker\Delta''$ of an arbitrary bidegree $(p,\,q)$. Then

$$0 = \langle\langle\Delta'_{p''}\,u,\,u\rangle\rangle = \langle\langle\Delta'(p''u),\,u\rangle\rangle = \langle\langle\Delta'\,u,\,u\rangle\rangle,$$

\noindent where the second identity followed from $p''\partial = \partial p''$ (which also implies $p''\partial^{\star} = \partial^{\star} p''$) and the last identity followed from $u\in\ker\Delta''$ (which amounts to $p''u = u$). Thus $\Delta'u=0$, i.e. $u\in\ker\Delta'$. This proves (\ref{eqn:kernels-intersection_E2}), so  Lemma \ref{Lem:kernels-intersection_E2} applies.

 $(ii)$\, We prove that inclusion (\ref{eqn:kernels-intersection_E2}) holds under the assumptions made. Let $u\in\ker\Delta'_{p''}\cap\ker\Delta''$. 

Note that the conjugate of Demailly's non-K\"ahler Bochner-Kodaira-Nakano identity $\Delta'' = \Delta'_{\tau} + T_{\omega}$ (cf. (\ref{eqn:BKN_demailly2})) is

\begin{equation}\label{eqn:BKN_conjugate_demailly2}\Delta' = \Delta''_{\tau} + \overline{T}_{\omega},\end{equation}

\noindent where $\Delta''_{\tau}:=[\bar\partial + \bar\tau,\, \bar\partial^{\star} + \bar\tau^{\star}]$ and $\overline{T}_{\omega} = [\Lambda,\,[\Lambda,\,\frac{i}{2}\,\partial\bar\partial\omega\wedge\cdot]] - [\bar\partial\omega\wedge\cdot,\,(\bar\partial\omega\wedge\cdot)^{\star}]$. Thanks to formula (\ref{eqn:BKN_conjugate_demailly2}), we have \begin{eqnarray}\label{eqn:Delta'-SKT-Delta''-tau1}\nonumber\langle\langle\Delta'\,u,\,u\rangle\rangle & = & \langle\langle(\Delta'' + [\bar\partial,\,\bar\tau^{\star}] + [\bar\tau,\,\bar\partial^{\star}])\,u,\,u\rangle\rangle + \langle\langle[\bar\tau,\,\bar\tau^{\star}]\,u,\,u\rangle\rangle - \langle\langle[\bar\partial\omega\wedge\cdot,\,(\bar\partial\omega\wedge\cdot)^{\star}]\,u,\,u\rangle\rangle\\  
    & = & \langle\langle[\bar\tau,\,\bar\tau^{\star}]\,u,\,u\rangle\rangle - \langle\langle[\bar\partial\omega\wedge\cdot,\,(\bar\partial\omega\wedge\cdot)^{\star}]\,u,\,u\rangle\rangle,\end{eqnarray}

\noindent where we have used the SKT assumption on $\omega$ to have $\bar{T}_{\omega}$ reduced to $-[\bar\partial\omega\wedge\cdot,\,(\bar\partial\omega\wedge\cdot)^{\star}]$ in formula (\ref{eqn:BKN_conjugate_demailly2}) and the argument below to infer that $ \langle\langle[\bar\partial,\,\bar\tau^{\star}]\,u,\,u\rangle\rangle = \langle\langle[\bar\tau,\,\bar\partial^{\star}]\,u,\,u\rangle\rangle = 0$ from the assumption $u\in\ker\Delta''= \ker\bar\partial\cap\ker\bar\partial^{\star}$: \begin{eqnarray}\label{eqn:Delta''-torsion-terms_vanishing}\nonumber \langle\langle[\bar\partial,\,\bar\tau^{\star}]\,u,\,u\rangle\rangle & = & \langle\langle\bar\tau^{\star}\,u,\,\bar\partial^{\star}u\rangle\rangle + \langle\langle\bar\partial u,\,\bar\tau u\rangle\rangle = 0 + 0 = 0, \\ 
    \langle\langle[\bar\tau,\,\bar\partial^{\star}]\,u,\,u\rangle\rangle & = & \langle\langle\bar\partial^{\star}\,u,\,\bar\tau^{\star}u\rangle\rangle + \langle\langle\bar\tau u,\,\bar\partial u\rangle\rangle = 0 + 0 = 0.\end{eqnarray}

Now, $\Delta' = \Delta'_{p''} + \Delta'_{p''_{\perp}}$, so the assumption $u\in\ker\Delta'_{p''}$ reduces (\ref{eqn:Delta'-SKT-Delta''-tau1}) to

\begin{eqnarray}\label{eqn:Delta'-SKT-Delta''-tau2}\langle\langle\Delta'_{p''_{\perp}}\,u,\,u\rangle\rangle =  \langle\langle[\bar\tau,\,\bar\tau^{\star}]\,u,\,u\rangle\rangle - \langle\langle[\bar\partial\omega\wedge\cdot,\,(\bar\partial\omega\wedge\cdot)^{\star}]\,u,\,u\rangle\rangle.\end{eqnarray}

\noindent The r.h.s. of (\ref{eqn:Delta'-SKT-Delta''-tau2}) vanishes thanks to the hypothesis (\ref{eqn:tau-bar-S}), so $\Delta'_{p''_{\perp}}\,u = 0$, hence also $\Delta'u=0$.    \hfill $\Box$

\begin{Rem}\label{Rem:SKT-R=0_E1} The proof of $(ii)$ of the above Theorem \ref{The:commutation_E2-commutation_Rbar} shows that if $X$ carries an SKT metric $\omega$ whose torsion satisfies the condition $[\tau,\,\tau^{\star}] = [\partial\omega\wedge\cdot,\,(\partial\omega\wedge\cdot)^{\star}]$, then the Fr\"olicher spectral sequence of $X$ degenerates at $E_1$.

\end{Rem}

\noindent {\it Proof.} To get $E_1$ degeneration, it suffices for the inclusion ${\cal H}_{\Delta''}^{p,\,q}(X,\,\C)\subset {\cal H}_{\Delta}^{p+q}(X,\,\C)$ of $\Delta''$-, resp. $\Delta$-harmonic spaces to hold for all $p,q$. (The argument is analogous to the one for $(a)$ of Lemma \ref{Lem:sufficient-cond_E2-deg}.) Now, (\ref{eqn:Delta'-SKT-Delta''-tau1}) holds for all $u\in{\cal H}_{\Delta''}^{p,\,q}(X,\,\C)$ if $\omega$ is SKT, hence $\Delta'u=0$ whenever $\Delta''u=0$ under the present assumptions. Then, by (\ref{eqn:Delta_gamma_gamma_pre-type}), we get $\Delta u=0$ for all $(p,\,q)$-forms $u$ satisfying $\Delta''u=0$ and for all $p,q$. This proves the above inclusion, hence the contention. \hfill $\Box$

\vspace{1ex}

\subsubsection{Alternative expression for the torsion operator $\bar{R}_{\omega}$}

 We shall now compute the operator $\bar{R}_{\omega}:=[\bar\tau,\,\bar\tau^{\star}] - [\bar\partial\omega\wedge\cdot,\,(\bar\partial\omega\wedge\cdot)^{\star}]$ featuring in $(ii)$ of Theorem \ref{The:commutation_E2-commutation_Rbar} in terms of the non-negative operator $\bar{S}_{\omega}$ (cf. (\ref{eqn:R-bar-omega_def})).

\begin{Lem}\label{Lem:tau-tau-star-S_computation} Let $(X,\,\omega)$ be an arbitrary compact Hermitian manifold of dimension $n$. Put $\bar{S}_{\omega}:=[\bar\partial\omega\wedge\cdot,\,(\bar\partial\omega\wedge\cdot)^{\star}] \geq 0$. The following formula holds:

\begin{equation}\label{eqn:tau-tau-star-S_computation}[\bar\tau,\,\bar\tau^{\star}] - \bar{S}_{\omega} = 2\bar{S}_{\omega} + [[\Lambda,\,\bar{S}_{\omega}],\,L],\end{equation}

\noindent where, as usual, $L=L_{\omega}:=\omega\wedge\cdot$. Moreover, for any bidegree $(p,\,q)$, $[[\Lambda,\,\bar{S}_{\omega}],\,L]$ is given by

\begin{eqnarray} \label{eqn:Lambda-S-L_formula}\nonumber \langle\langle[[\Lambda,\,\bar{S}_{\omega}],\,L]\,u,\,u\rangle\rangle = \langle\langle \bar{S}_{\omega}(\omega\wedge u),\,\omega\wedge u\rangle\rangle & + & \langle\langle \bar{S}_{\omega}(\Lambda\, u),\,\Lambda\, u\rangle\rangle + (p+q-n)\,\langle\langle \bar{S}_{\omega}u,\,u\rangle\rangle \\
     & - & 2\mbox{Re}\,\langle\langle\Lambda(\bar{S}_{\omega}u),\,\Lambda\,u\rangle\rangle, \hspace{3ex} u\in C^{\infty}_{p,\,q}(X,\,\C).\end{eqnarray}

\end{Lem}

\noindent {\it Proof.} Since $\tau = [\Lambda,\,\partial\omega\wedge\cdot]$, we get

\begin{eqnarray}\label{eqn:tau-bar_tau-bar-star}[\bar\tau,\,\bar\tau^{\star}] = \bigg[[\Lambda,\,\bar\partial\omega\wedge\cdot],\,[(\bar\partial\omega\wedge\cdot)^{\star},\,L]\bigg] = \bigg[[[(\bar\partial\omega\wedge\cdot)^{\star},\,L],\,\Lambda],\,\bar\partial\omega\wedge\cdot\bigg] - \bigg[[\bar\partial\omega\wedge\cdot,\,[(\bar\partial\omega\wedge\cdot)^{\star},\,L]],\,\Lambda\bigg],\end{eqnarray}

\noindent where the last identity followed from Jacobi's identity applied to the operators $[(\bar\partial\omega\wedge\cdot)^{\star},\,L]$, $\Lambda$ and $\bar\partial\omega\wedge\cdot$. 

 To compute the first factor in the first term on the r.h.s. of (\ref{eqn:tau-bar_tau-bar-star}), we apply again Jacobi's identity:

\begin{eqnarray}\label{eqn:1st-1st}\bigg[[(\bar\partial\omega\wedge\cdot)^{\star},\,L],\,\Lambda\bigg] = -\bigg[[L,\,\Lambda],\,(\bar\partial\omega\wedge\cdot)^{\star}\bigg] - \bigg[[\Lambda,\,(\bar\partial\omega\wedge\cdot)^{\star}],\,L\bigg].\end{eqnarray}

\noindent Using the standard fact that $[L,\,\Lambda] = (p+q-n)\,\mbox{Id}$ on $(p,\,q)$-forms, for any $(p,\,q)$-form $u$ we get

\begin{eqnarray}\nonumber \bigg[[L,\,\Lambda],\,(\bar\partial\omega\wedge\cdot)^{\star}\bigg]\,u & = & [L,\,\Lambda]\bigg((\bar\partial\omega\wedge\cdot)^{\star}u\bigg) - (\bar\partial\omega\wedge\cdot)^{\star}\bigg([L,\,\Lambda],\,u\bigg)\\
\nonumber & = & (p+q-3-n)\,(\bar\partial\omega\wedge\cdot)^{\star}u - (\bar\partial\omega\wedge\cdot)^{\star}((p+q-n)\,u) = -3(\bar\partial\omega\wedge\cdot)^{\star}u.\end{eqnarray}

\noindent Thus $\bigg[[L,\,\Lambda],\,(\bar\partial\omega\wedge\cdot)^{\star}\bigg] = -3(\bar\partial\omega\wedge\cdot)^{\star}$. On the other hand, $[\Lambda,\,(\bar\partial\omega\wedge\cdot)^{\star}] = [\bar\partial\omega\wedge\cdot,\,L]^{\star} =0$ since, clearly, $[\bar\partial\omega\wedge\cdot,\,L]\,u = \bar\partial\omega\wedge\omega\wedge u - \omega\wedge\bar\partial\omega\wedge u = 0$ for any $u$. Therefore, (\ref{eqn:1st-1st}) reduces to

\begin{eqnarray}\label{eqn:1st-1st2}\bigg[[(\bar\partial\omega\wedge\cdot)^{\star},\,L],\,\Lambda\bigg] = 3(\bar\partial\omega\wedge\cdot)^{\star}.\end{eqnarray}

 Similarly, to compute the first factor in the second term on the r.h.s. of (\ref{eqn:tau-bar_tau-bar-star}), we start by applying Jacobi's identity:

\begin{eqnarray}\label{eqn:1st-2nd}\nonumber\bigg[\bar\partial\omega\wedge\cdot,\,[(\bar\partial\omega\wedge\cdot)^{\star},\,L]\bigg] & = & \bigg[(\bar\partial\omega\wedge\cdot)^{\star},\,[L,\,\bar\partial\omega\wedge\cdot]\bigg] -\bigg[L,\,[\bar\partial\omega\wedge\cdot,\,(\bar\partial\omega\wedge\cdot)^{\star}]\bigg]\\
      & = & \bigg[[\bar\partial\omega\wedge\cdot,\,(\bar\partial\omega\wedge\cdot)^{\star}],\,L\bigg] = [\bar{S}_{\omega},\,L],\end{eqnarray}

\noindent where the last but one identity followed from $[L,\,\bar\partial\omega\wedge\cdot] = 0$ seen above.

 Putting together (\ref{eqn:tau-bar_tau-bar-star}), (\ref{eqn:1st-1st2}) and (\ref{eqn:1st-2nd}), we get:

\begin{eqnarray}\label{eqn:tau-bar_tau-bar-star2}[\bar\tau,\,\bar\tau^{\star}] = 3\bar{S}_{\omega} - [[\bar{S}_{\omega},\,L],\,\Lambda].\end{eqnarray}

\noindent A new application of Jacobi's identity spells

\begin{equation}\label{eqn:Jacobi_S-L-Lambda}[[\bar{S}_{\omega},\,L],\,\Lambda] + [[L,\,\Lambda],\,\bar{S}_{\omega}] + [[\Lambda,\,\bar{S}_{\omega}],\,L] =0,    \hspace{2ex} \mbox{which gives} \hspace{2ex} -[[\bar{S}_{\omega},\,L],\,\Lambda] = [[\Lambda,\,\bar{S}_{\omega}],\,L].\end{equation}

\noindent Indeed, since $[L,\,\Lambda] = (p+q-n)\,\mbox{Id}$ on $(p,\,q)$-forms and $\bar{S}_{\omega}$ is an operator of type $(0,\,0)$, we get $[[L,\,\Lambda],\,\bar{S}_{\omega}] = 0$ which accounts for the last statement in (\ref{eqn:Jacobi_S-L-Lambda}).

 It is now clear that the combined (\ref{eqn:tau-bar_tau-bar-star2}) and (\ref{eqn:Jacobi_S-L-Lambda}) prove (\ref{eqn:tau-tau-star-S_computation}).

 To prove (\ref{eqn:Lambda-S-L_formula}), we start by computing \begin{eqnarray}\label{eqn:Lambda-S-L__computation}\nonumber \langle\langle[[\Lambda,\,\bar{S}_{\omega}],\,L]\,u,\,u\rangle\rangle = \langle\langle[\Lambda,\,\bar{S}_{\omega}]\,(\omega\wedge u),\,u\rangle\rangle & - & \langle\langle\omega\wedge[\Lambda,\,\bar{S}_{\omega}]\, u,\,u\rangle\rangle \\
  \nonumber  =  \langle\langle \bar{S}_{\omega}(\omega\wedge u),\,\omega\wedge u\rangle\rangle & - & \langle\langle \omega\wedge u,\, \omega\wedge \bar{S}_{\omega}u\rangle\rangle + \langle\langle \bar{S}_{\omega}(\Lambda\,u),\, \Lambda\,u\rangle\rangle\\
 & - & \langle\langle \Lambda(\bar{S}_{\omega} u),\, \Lambda\,u\rangle\rangle.\end{eqnarray}

 Then we notice the general fact that for every $(p,\,q)$-forms $u,v$ we have:

\begin{equation}\label{eqn:omega-Lambda-wedge-product}\langle\langle\omega\wedge u,\,\omega\wedge v\rangle\rangle = \langle\langle\Lambda\, u,\,\Lambda\, v\rangle\rangle - (p+q-n)\, \langle\langle u,\,v\rangle\rangle.\end{equation}

\noindent Indeed, $\langle\langle\omega\wedge u,\,\omega\wedge v\rangle\rangle = \langle\langle\Lambda(\omega\wedge u),\, v\rangle\rangle$ and $\Lambda(\omega\wedge u) = \omega\wedge\Lambda\,u - (p+q-n)\,u$. Now, applying (\ref{eqn:omega-Lambda-wedge-product}), we get

\begin{equation}\label{eqn:omega-Lambda-wedge-product_applied}\langle\langle \omega\wedge u,\, \omega\wedge \bar{S}_{\omega}u\rangle\rangle = \langle\langle \Lambda\, u,\, \Lambda(\bar{S}_{\omega}u)\rangle\rangle - (p+q-n)\,\langle\langle u,\,\bar{S}_{\omega}u\rangle\rangle.\end{equation}

 It is now clear that the combined (\ref{eqn:Lambda-S-L__computation}) and (\ref{eqn:omega-Lambda-wedge-product_applied}) prove (\ref{eqn:Lambda-S-L_formula}) because $\langle\langle \Lambda\, u,\, \Lambda(\bar{S}_{\omega}u)\rangle\rangle$ is the conjugate of $\langle\langle \Lambda(\bar{S}_{\omega} u),\, \Lambda\,u\rangle\rangle$ and $\langle\langle u,\,\bar{S}_{\omega}u\rangle\rangle = \langle\langle \bar{S}_{\omega}u,\,u\rangle\rangle$.  \hfill $\Box$

\vspace{1ex}

\subsubsection{Putting the hypothesis $\partial p'' = p''\partial$ in context}

 We now reinterpret the commutation of $\partial$ with $p''$ (the simplest sufficient condition for $E_2(X)=E_{\infty}(X)$ found so far, cf. Theorem \ref{The:commutation_E2-commutation_Rbar}).

\begin{Lem}\label{Lem:reinterpret_commutation_E2} Let $(X,\,\omega)$ be a compact Hermitian manifold. The following implication and equivalence hold:

\begin{equation}\label{eqn:reinterpret_commutation_E2}\partial\Delta'' = \Delta''\partial \implies \partial p'' = p''\partial  \iff \partial(\ker\Delta'')\subset\ker\Delta'' \hspace{1ex} \mbox{and} \hspace{1ex} \partial^{\star}(\ker\Delta'')\subset\ker\Delta''.\end{equation}

\end{Lem}

\noindent {\it Proof.} Suppose that $\partial\Delta'' = \Delta''\partial$. Then, taking adjoints, we also have $\Delta''\partial^{\star} = \partial^{\star}\Delta''$. These identities immediately imply

\begin{equation}\label{eqn:del-ker-Delta''_inclusion}\partial(\ker\Delta'')\subset\ker\Delta'' \hspace{1ex} \mbox{and} \hspace{1ex} \partial^{\star}(\ker\Delta'')\subset\ker\Delta''.\end{equation}

 Now suppose that (\ref{eqn:del-ker-Delta''_inclusion}) holds. We shall prove that $\partial p'' = p''\partial$. Let $u$ be an arbitrary smooth form. Then $u$ splits as $u = u_0 + \bar\partial v + \bar\partial^{\star}w$ with $u_0\in\ker\Delta''$. Thus $\partial p''u=\partial u_0$ and 

$$p''\partial u = p''\partial u_0 + p''\partial\bar\partial v + p''\partial\bar\partial^{\star}w = \partial u_0 + p''\partial\bar\partial^{\star}w$$

\noindent because $\partial u_0\in\ker\Delta''$ by (\ref{eqn:del-ker-Delta''_inclusion}) and $p''\partial\bar\partial v = -p''\bar\partial\partial v =0$ since $\mbox{Im}\,\bar\partial\perp\ker\Delta''$. We now prove that $p''\partial\bar\partial^{\star}w = 0$ and this will show that $\partial p''u = p''\partial u$, as desired. Proving that $p''\partial\bar\partial^{\star}w = 0$ is equivalent to proving that $\partial\bar\partial^{\star}w\in(\ker\Delta'')^{\perp}$. Let $\zeta\in\ker\Delta''$, arbitrary. We have

$$\langle\langle\zeta,\,\partial\bar\partial^{\star}w\rangle\rangle = \langle\langle\partial^{\star}\zeta,\,\bar\partial^{\star}w\rangle\rangle =0$$

\noindent because $\partial^{\star}\zeta\in\ker\Delta''$ thanks to (\ref{eqn:del-ker-Delta''_inclusion}), $\bar\partial^{\star}w\in\mbox{Im}\,\bar\partial^{\star}$ and $\ker\Delta''\perp\mbox{Im}\,\bar\partial^{\star}$.

 It remains to prove that if $\partial p'' = p''\partial$, then (\ref{eqn:del-ker-Delta''_inclusion}) holds. Note the general fact that for any form $u$, $u\in\ker\Delta''$ iff $p''u=u$. Let us now suppose that $\partial p'' = p''\partial$. Then, taking adjoints, we also have $\partial^{\star} p'' = p''\partial^{\star}$, so (\ref{eqn:del-ker-Delta''_inclusion}) holds.  \hfill $\Box$

\vspace{2ex}

\begin{center} {\Large\bf Part II: $E_2$ degeneration of the spectral sequence associated with a pair of foliations}

\end{center}

\section{Review of standard material}\label{section:introduction2}

Let $X$ be a compact complex manifold of dimension $n$ endowed with a complementary pair of regular foliations $(N,\,F)$ (cf. (\ref{eqn:product-structure_def})). We keep the notation used in the Introduction.

\vspace{2ex}

\noindent{\bf Global picture.} The integrability of both $F$ and $N$ induces a splitting 

\begin{equation}\label{eqn:del-splitting}\partial = \partial_N + \partial_F\end{equation}

\noindent of $\partial$ into components of $(N,\, F)$-types $(1,\,0)$, resp. $(0,\,1)$, that are intrinsically defined by the following formulae (see e.g. [Raw77]) modelled on the classical Cartan formula. For any $\alpha\in C^{\infty}(X,\,\Lambda^qF^{\star})$ (i.e. for any global smooth $(q,\,0)$-form of $(N,\,F)$-type $(0,\,q)$) and any $\xi_0, \dots , \xi_q\in C^{\infty}(X,\,F)$ (i.e. for any global vector fields of $(N,\, F)$-type $(0,\,1)$), we set

\begin{eqnarray}\label{eqn:del_F-def}\nonumber (\partial_F\alpha)(\xi_0,\dots , \xi_q) & := & \sum\limits_{j=0}^q(-1)^j\,\xi_j\cdot\alpha(\xi_0,\dots , \widehat{\xi_j},\dots , \xi_q) \\
      & + & \sum\limits_{0\leq j<k\leq q}(-1)^{j+k}\,\alpha([\xi_j,\,\xi_k],\,\xi_0,\dots , \widehat{\xi_j},\dots , \widehat{\xi_k},\dots , \xi_q).\end{eqnarray}

\noindent This defines a form $\partial_F\alpha\in C^{\infty}(X,\,\Lambda^{q+1}F^{\star})$. More generally, for any $\alpha\in C^{\infty}(X,\,\Lambda^pN^{\star}\otimes\Lambda^qF^{\star})$ and any $\xi_0, \dots , \xi_{p-1}\in C^{\infty}(X,\,N)$, $\xi_p, \dots , \xi_{p+q}\in C^{\infty}(X,\,F)$, we set

\noindent $\displaystyle (\partial_F\alpha)(\xi_0,\dots , \xi_{p-1},\xi_p,\dots , \xi_{p+q}) := (-1)^p\,\sum\limits_{j=0}^q(-1)^j\,\xi_j\cdot\alpha(\xi_0,\dots , \xi_{p-1}, \dots , \widehat{\xi_{p+j}},\dots , \xi_{p+q})$

\begin{eqnarray}\label{eqn:del_F-def_full} \nonumber   & + & \sum\limits_{0\leq j<k\leq q}(-1)^{j+k}\,\alpha(\xi_0,\dots , \xi_{p-1},\,[\xi_{p+j},\,\xi_{p+k}],\dots , \widehat{\xi_{p+j}},\dots , \widehat{\xi_{p+k}},\dots , \xi_{p+q}) \\
           & + & (-1)^p\,\sum\limits_{j=0}^{p-1}\sum\limits_{k=0}^q(-1)^{j+k}\,\alpha([\xi_j,\,\xi_{p+k}],\,\xi_0,\dots , \widehat{\xi_j},\dots , \xi_{p-1}, \dots, \widehat{\xi_{p+k}}, \dots , \xi_{p+q}).\end{eqnarray}

\noindent This defines a form $\partial_F\alpha\in C^{\infty}(X,\,\Lambda^pN^{\star}\otimes\Lambda^{q+1}F^{\star})$. Again by the integrability of $F$, we have $\partial_F^2=0$. Analogous formulae define $\partial_N$ and the integrability of $N$ implies $\partial_N^2=0$. Thus $\partial_N$ and $\partial_F$ are exterior differentials in the holomorphic $N$-directions, resp. the holomorphic $F$-directions. They are of the respective $(N,\,F)$-types $(1,\,0)$ and $(0,\,1)$ and anti-commute, i.e. $\partial_N\partial_F + \partial_F\partial_N =0$.

\vspace{2ex}

\noindent{\bf Local picture.} The integrability of both $F$ and $N$ implies that in a neighbourhood $U$ of any point $x\in X$ there exist local holomorphic coordinates $z_1, \dots z_r, z_{r+1}, \dots , z_n$ centred at $x$ such that $\{\frac{\partial}{\partial z_1}, \dots , \frac{\partial}{\partial z_r}\}$ is a frame for $N$ over $U$ and $\{\frac{\partial}{\partial z_{r+1}}, \dots , \frac{\partial}{\partial z_n}\}$ is a frame for $F$ over $U$. Such coordinates will be called {\it local product coordinates} by analogy with the terminology of [Rei58]. Put $z'=(z_1, \dots , z_r)$ (the $N$-tangent coordinates) and $z''=(z_{r+1}, \dots , z_n)$ (the $F$-tangent coordinates). Thus, the leaves of $F$ are locally defined by the equations 

\vspace{1ex}

\hspace{25ex} $\{z_1=c_1,\dots , z_r=c_r\}$,  \hspace{3ex} where $c_1,\dots , c_r$ are constants,

\vspace{1ex}

\noindent while the operators $\partial_N$ and $\partial_F$ are locally defined by

\begin{equation}\label{eqn:del-coord-def}\partial_N = \sum\limits_{l=1}^rdz_l\wedge\frac{\partial}{\partial z_l} \hspace{2ex} \mbox{and} \hspace{2ex} \partial_F = \sum\limits_{l=r+1}^ndz_l\wedge\frac{\partial}{\partial z_l}.\end{equation}

 There always exists (e.g. by the analogue of the argument for the real case in [Rei58, p.245-246]) a Hermitian metric $\omega$ on $X$, called henceforth a {\it product metric}, which in every system of local product coordinates $z_1, \dots z_r, z_{r+1}, \dots , z_n$ has the shape 

\begin{equation}\label{eqn:almost-product-metric_def}\omega = \sum\limits_{j,\,k=1}^r\omega_{j\bar{k}}(z',\,z'')\,idz_j\wedge d\bar{z}_k + \sum\limits_{j,\,k=r+1}^n\omega_{j\bar{k}}(z',\,z'')\,idz_j\wedge d\bar{z}_k.\end{equation}

\noindent In general, the coefficients $\omega_{j\bar{k}}$ are functions on $U$ depending on both groups of coordinates $z'$ and $z''$. We put

\begin{equation}\label{eqn:omega_NF_def}\omega_N:=\sum\limits_{j,\,k=1}^r\omega_{j\bar{k}}(z',\,z'')\,idz_j\wedge d\bar{z}_k  \hspace{2ex} \mbox{and} \hspace{2ex} \omega_F:=\sum\limits_{j,\,k=r+1}^n\omega_{j\bar{k}}(z',\,z'')\,idz_j\wedge d\bar{z}_k.\end{equation}

\noindent For some of our results that follow, we shall need to assume the existence of a special kind of Hermitian metric on $X$ adapted to the $(N,\,F)$ structure. The name is borrowed from [Rei59] (where the real case was studied) whose analogue to the complex Hermitian context we consider.

\begin{Def}\label{Def:bundle-like-metric} A {\bf bundle-like Hermitian metric} on a compact complex manifold $X$ endowed with an integrable holomorphic almost product structure $(N,\,F)$ is a Hermitian metric $\omega$ which in every system of local product coordinates $z_1, \dots z_r, z_{r+1}, \dots , z_n$ has the shape

\begin{equation}\label{eqn:bundle-like-metric_def}\omega = \sum\limits_{j,\,k=1}^r\omega_{j\bar{k}}(z')\,idz_j\wedge d\bar{z}_k + \sum\limits_{j,\,k=r+1}^n\omega_{j\bar{k}}(z'')\,idz_j\wedge d\bar{z}_k.\end{equation}

\end{Def}

\vspace{2ex}

\noindent{\bf The spectral sequence induced by $(N,\,F)$.} Put ${\cal E}^k(X):=C^{\infty}(X,\,\Lambda^{k,\,0}T^{\star}X)$. For every $p,q\in\{0, \dots , k\}$, let $E^{p,\,q}$ denote the holomorphic vector subbundle $\Lambda^pN^{\star}\otimes\Lambda^qF^{\star}$ of $\Lambda^{k,\,0}T^{\star}X$, and let ${\cal E}^{p,\,q}(X) = {\cal E}^{p,\,q}_{N,\,F}(X):=C^{\infty}(X,\,E^{p,\,q})$ stand for the space of its global smooth sections. \noindent The integrable operators

\begin{equation}\label{eqn:double-complex}\partial_N:{\cal E}^{p,\,q}(X)\to{\cal E}^{p+1,\,q}(X) \hspace{2ex} \mbox{and} \hspace{2ex} \partial_F:{\cal E}^{p,\,q}(X)\to{\cal E}^{p,\,q+1}(X)\end{equation}

\noindent define a double complex ${\cal E}^{\bullet,\,\bullet}(X)$ with the total differential $\partial=\partial_N + \partial_F$. We consider the spectral sequence associated in the standard way with this double complex (see e.g. [Dem 96, $\S.9$] or [Voi02, $\S. 8.3.2$]). As usual, the first two steps in the spectral sequence are defined by the second differential $\partial_F$, resp. the first differential $\partial_N$. Indeed, we put $E_0^{p,\,q}:={\cal E}^{p,\,q}(X)$ and $d_0:=\partial_F$, so the groups $E_1 = E_1(N,\,F)$ are defined as the cohomology of the complex 

\begin{equation}\label{eqn:del_F_complex} \cdots\stackrel{\partial_F}{\longrightarrow}{\cal E}^{p,\,q-1}(X)\stackrel{\partial_F}{\longrightarrow}{\cal E}^{p,\,q}(X)\stackrel{\partial_F}{\longrightarrow}{\cal E}^{p,\,q+1}(X)\stackrel{\partial_F}{\longrightarrow}\cdots,\end{equation}

\noindent i.e. $E_1^{p,\,q}= E_1^{p,\,q}(N,\,F)=H^q({\cal E}^{p,\,\bullet}(X),\,\partial_F)$, while the differentials $d_1$ are induced by $\partial_N$:

\begin{equation}\label{eqn:d_1_complex}\cdots\stackrel{d_1}{\longrightarrow} E^{p-1,\,q}_1\stackrel{d_1}{\longrightarrow}E_1^{p,\,q}\stackrel{d_1}{\longrightarrow}E_1^{p+1,\,q}\stackrel{d_1}{\longrightarrow}\cdots.\end{equation}

\noindent This means that for any form $\alpha\in {\cal E}^{p,\,q}(X)$ such that $\partial_F\alpha=0$, the class $[\alpha]_{\partial_F}\in E_1^{p,\,q}$ is mapped by $d_1$ to the class $[\partial_N\alpha]_{\partial_F}\in E_1^{p+1,\,q}$. This is meaningful since from $\partial = \partial_N + \partial_F$ and from $\partial_N^2=0$, $\partial_F^2=0$ (the integrability of $N$ and $F$), we get $\partial_N\partial_F + \partial_F\partial_N=0$, so $\partial_F(\partial_N\alpha)=0$, which means that the class $[\partial_N\alpha]_{\partial_F}\in E_1^{p+1,\,q}$ is well defined. The map $d_1$ is also well defined since $d_1([\alpha]_{\partial_F}) = [\partial_N\alpha]_{\partial_F}$ is independent of the choice of representative $\alpha$ of the $\partial_F$-class $[\alpha]_{\partial_F}$. (Indeed, if $\alpha=\partial_F\beta$, then $\partial_N\alpha = \partial_F(-\partial_N\beta)\in\mbox{Im}\,\partial_F$.) Moreover, $d_1^2=0$ since $\partial_N^2=0$, so (\ref{eqn:d_1_complex}) is indeed a complex. The groups $E_2 = E_2(N,\,F)$ are defined as the cohomology of this complex, i.e. for all $p,q$ we have

\begin{equation}\label{eqn:E_2_def_NF}E_2^{p,\,q} = E_2^{p,\,q}(N,\,F): = H^p(E_1^{\bullet,\,q},\,d_1) = \bigg\{\bigg[[\alpha]_{\partial_F}\bigg]_{d_1}\bigg\slash\alpha\in {\cal E}^{p,\,q}(X)\cap\ker\partial_F \hspace{1ex} \mbox{and} \hspace{1ex} \partial_N\alpha\in\mbox{Im}\,\partial_F\bigg\},    \end{equation}

\noindent so the elements of $E_2^{p,\,q}$ are $d_1$-classes of $\partial_F$-classes.

 This spectral sequence converges to the $\partial$-cohomology in bidegrees $(k,\,0)$ of the manifold $X$, i.e. if we let $E_{\infty}^{p,\,q}:=\lim_{s\to +\infty}E_s^{p,\,q}$, we have canonical isomorphisms

\begin{equation}\label{eqn:Einfty_sum}H^{k,\,0}_{\partial}(X,\,\C)\simeq\bigoplus\limits_{p+q=k}E_{\infty}^{p,\,q},   \hspace{3ex} k=0, 1, \dots , n.\end{equation}

 Note that the vector spaces $E_s^{p,\,q}$ need not be finite-dimensional since the system $(\partial_N,\,\partial_F)$ is not elliptic in general. Indeed, if we have fixed a Hermitian metric $\omega$ on $X$ and denote by $\partial_N^{\star},\partial_F^{\star}$ the adjoints of $\partial_N,\partial_F$ w.r.t. the $L^2$ inner product induced by $\omega$, the associated Laplace-Beltrami operators

$$\Delta'_N=[\partial_N,\,\partial_N^{\star}]=\partial_N\partial_N^{\star} + \partial_N^{\star}\partial_N \hspace{2ex}\mbox{and}\hspace{2ex} \Delta'_F=[\partial_F,\,\partial_F^{\star}]=\partial_F\partial_F^{\star} + \partial_F^{\star}\partial_F$$

\noindent are not elliptic in general (since each ``misses'' the complementary directions). However, a fact that will be important for us is that the sum of these Laplacians

$$\Delta'_N + \Delta'_F$$

\noindent is an elliptic operator. So is also the usual $\partial$-Laplacian

$$\Delta' = [\partial,\,\partial^{\star}] = \Delta'_N + \Delta'_F + [\partial_N,\,\partial_F^{\star}] + [\partial_F,\,\partial_N^{\star}].$$

\noindent Since $X$ is compact, $\ker\Delta'_N =\ker\partial_N\cap\ker\partial_N^{\star}$ and $\ker\Delta'_F =\ker\partial_F\cap\ker\partial_F^{\star}$.

\section{The non-differential Laplacian and $E_2$ degeneration}\label{section:FI-degenration}

 Let $X$ be an $n$-dimensional compact complex manifold equipped with an integrable holomorphic almost product structure $(N,\,F)$ and with a product Hermitian metric $\omega$. By analogy with the Fr\"olicher case described in Part I, we consider the following operators

\begin{equation}\label{eqn:pNpF_def}p'_N: {\cal E}^{p,\,q}(X)\longrightarrow\ker\Delta'_N \hspace{2ex} \mbox{and} \hspace{2ex} p'_F: {\cal E}^{p,\,q}(X)\longrightarrow\ker\Delta'_F\end{equation}

\noindent the orthogonal projections onto the $\Delta'_N$-harmonic, resp. the $\Delta'_F$-harmonic spaces. Similarly, let

\begin{equation}\label{eqn:pNpF_perp_def}p_N^{'\perp}: {\cal E}^{p,\,q}(X)\longrightarrow(\ker\Delta'_N)^{\perp}  \hspace{2ex} \mbox{and} \hspace{2ex} p_F^{'\perp}: {\cal E}^{p,\,q}(X)\longrightarrow(\ker\Delta'_F)^{\perp}\end{equation}

\noindent denote the orthogonal projections onto the orthogonal complements of the respective harmonic spaces. The operators $p'_N, p'_F, p_N^{'\perp}, p_F^{'\perp}$ depend on the metric $\omega$ and are no longer pseudo-differential operators, let alone regularising or of finite rank, since the kernels of $\Delta'_N$ and $\Delta'_F$ are no longer finite-dimensional due to the lack of ellipticity of $\Delta'_N$ and $\Delta'_F$. They clearly satisfy the properties:

\begin{equation}\label{eqn:orth-proj_prop_NF}p'_N = (p'_N)^{\star} = (p'_N)^2, \hspace{2ex} p'_F = (p_F)^{'\star} = (p'_F)^2, \hspace{2ex} p_N^{'\perp} = (p_N^{'\perp})^{\star} = (p_N^{'\perp})^2, \hspace{2ex} p_F^{'\perp} = (p_F^{'\perp})^{\star} = (p_F^{'\perp})^2.\end{equation}

 Again by analogy with the Fr\"olicher case, we define our main object of study in this second part.

\begin{Def}\label{Def:Delta-tilde_def_NF} Let $X$ be a compact complex manifold with $\mbox{dim}_{\C}X=n$ equipped with an integrable holomorphic almost product structure $(N,\,F)$ and with a product Hermitian metric $\omega$. For every $p, q$, we define the operator $\widetilde{\Delta'} : {\cal E}^{p,\,q}(X)\longrightarrow {\cal E}^{p,\,q}(X)$ by

\begin{equation}\label{eqn:Delta-tilde_def_NF}\widetilde{\Delta'} : = \partial_N p'_F\partial_N^{\star} + \partial_N^{\star}p'_F\partial_N + \partial_F\partial_F^{\star} + \partial_F^{\star}\partial_F.\end{equation}

\noindent In other words, we have

\begin{equation}\label{eqn:Delta-tilde_def_bis_NF}\widetilde{\Delta'} = \Delta'_{N,\,p'_F} + \Delta'_F, \hspace{2ex} \mbox{where} \hspace{1ex} \Delta'_{N,\,p'_F}:= \partial_N p'_F\partial_N^{\star} + \partial_N^{\star}p'_F\partial_N : {\cal E}^{p,\,q}(X)\longrightarrow {\cal E}^{p,\,q}(X).\end{equation}

\noindent Thus $\widetilde{\Delta'}$ is the sum of a Fourier integral operator ($\Delta'_{N,\,p'_F}$) and a non-elliptic differential operator of order two (the $\partial_F$-Laplacian $\Delta'_F$).

\end{Def}

 The idea we shall now be pursuing is to find a hypothesis ensuring that $\widetilde{\Delta'}$ satisfies G$\mathring{a}$rding's inequality by ensuring that $\Delta'_{N,\,p'_F}$ dominates a constant multiple of $\Delta'_N$ and then using the ellipticity of $\Delta'_N + \Delta'_F$ and G$\mathring{a}$rding's inequality it satisfies.

\begin{The}\label{The:Garding_sufficient-cond} Suppose that for every $p,q$ the following identity holds:

\begin{equation}\label{eqn:Garding_sufficient-cond}\ker(\Delta'_N:{\cal E}^{p,\,q}(X)\rightarrow{\cal E}^{p,\,q}(X)) + \ker(\Delta'_F:{\cal E}^{p,\,q}(X)\rightarrow{\cal E}^{p,\,q}(X)) = {\cal E}^{p,\,q}(X).\end{equation}

\noindent $(i)$\, There exists a constant $0<\varepsilon<1$ such that for all $p,q$ we have

\begin{equation}\label{eqn:Delta'Np_Delta'N}\langle\langle\Delta'_{N,\,p'_F}u,\,u\rangle\rangle \geq (1-\varepsilon)\,\langle\langle\Delta'_Nu,\,u\rangle\rangle, \hspace{3ex} u\in{\cal E}^{p,\,q}(X).\end{equation}

\noindent $(ii)$\, There exist constants $\delta_1,\delta_2>0$ such that for all $p,q$, {\bf G$\mathring{a}$rding's inequality} holds for $\widetilde{\Delta'}$:

\begin{equation}\label{eqn:Garding_NF}\langle\langle\widetilde{\Delta'}u,\,u\rangle\rangle + \delta_1\,||u||^2 \geq \delta_2\,||u||_1^2,  \hspace{3ex} u\in{\cal E}^{p,\,q}(X),\end{equation} 

\noindent where $||\,\,||_1$ stands for the Sobolev norm $W^1$ and $||\,\,||$ stands for the $L^2$ norm.

\noindent $(iii)$\, The above G$\mathring{a}$rding's inequality implies in turn that $\ker\widetilde{\Delta'}$ is finite-dimensional, that the image $\mbox{Im}\,\widetilde{\Delta'}$ is closed in ${\cal E}^{p,\,q}(X)$ and that the following $3$-space orthogonal decomposition holds:

\begin{equation}\label{eqn:3-space-decomp_Delta-tilde_NF}{\cal E}^{p,\,q}(X) = \ker\widetilde{\Delta'}\bigoplus\bigg(\mbox{Im}\,\partial_F + \mbox{Im}\,(\partial_{N|\ker\Delta'_F})\bigg)\bigoplus\bigg(\mbox{Im}\,(\partial^{\star}_N\circ p'_F) +  \mbox{Im}\,\partial_F^{\star}\bigg).\end{equation}

\noindent Moreover, the decomposition (\ref{eqn:3-space-decomp_Delta-tilde_NF}) also holds when $\mbox{Im}\,(\partial_{N|\ker\Delta'_F})$ is replaced with $\mbox{Im}\,(\partial_{N|\ker\partial_F})$.

\noindent $(iv)$\, If, moreover, $[\partial_N,\,\partial_F^{\star}]=0$, then $\mbox{Im}\,\partial_F$ is closed in ${\cal E}^{p,\,q}(X)$ and the following {\bf Hodge isomorphism} holds:

\begin{equation}\label{eqn:Hodge-isom_E2_NF}{\cal H}_{\widetilde{\Delta'}}^{p,\,q}(N,\,F):=\ker\bigg(\widetilde{\Delta'}:{\cal E}^{p,\,q}(X)\longrightarrow{\cal E}^{p,\,q}(X)\bigg)\simeq E_2^{p,\,q}(N,\,F), \hspace{3ex} \alpha\longmapsto \bigg[[\alpha]_{\partial_F}\bigg]_{d_1}.\end{equation}

\noindent Thus, every class $[[\alpha]_{\partial_F}]_{d_1}\in E_2^{p,\,q}(N,\,F)$ contains a unique $\widetilde{\Delta'}$-harmonic representative $\alpha$. In particular, $\mbox{dim}_{\C}E_2^{p,\,q}(N,\,F)<+\infty$ for all $p,q$.

\noindent $(v)$\, Much more holds under the above assumption $[\partial_N,\,\partial_F^{\star}]=0$: the spectral sequence induced by $(N,\,F)$ {\bf degenerates at $E_2$} (i.e. $E_2(N,\,F)=E_{\infty}(N,\,F)$.)

\end{The}

\noindent {\it Proof.} $(i)$\, The hypothesis $\ker\Delta'_N + \ker\Delta'_F = {\cal E}^{p,\,q}(X)$ is equivalent (using the Open Mapping Theorem in Fr\'echet spaces) to each of the following equivalent conditions: \\

\noindent $(\ker\Delta'_N)^{\perp}\cap(\ker\Delta'_F)^{\perp} = \{0\} \iff \mbox{the map}\hspace{1ex} p'_F:(\ker\Delta'_N)^{\perp}\longrightarrow\ker\Delta'_F \hspace{2ex} \mbox{is injective}$ \\

 \hspace{20ex}  $ \iff  \exists 0<\varepsilon<1 \hspace{1ex} \mbox{such that} \hspace{1ex} ||p'_Fv||^2\geq (1-\varepsilon)\,||v||^2 \hspace{2ex} \mbox{for all}\hspace{1ex} v\in(\ker\Delta'_N)^{\perp}.$

\vspace{2ex}

\noindent The constant in the last inequality is necessarily in the interval $(0,\,1)$ since $||p'_Fv||\leq||v||$ for all forms $v$. The choices $v:=\partial_Nu$ and $v:=\partial_N^{\star}u$ are allowed for all the forms $u$ since $\mbox{Im}\,\partial_N\perp\ker\Delta'_N$ and $\mbox{Im}\,\partial_N^{\star}\perp\ker\Delta'_N$. Thus, we obtain

\vspace{1ex}

\hspace{3ex} $\langle\langle\Delta'_{N,\,p'_F}u,\,u\rangle\rangle = ||p'_F\partial_Nu||^2 + ||p'_F\partial_N^{\star}u||^2 \geq (1-\varepsilon)\,(||\partial_Nu||^2 + ||\partial_N^{\star}u||^2) \geq (1-\varepsilon)\,\langle\langle\Delta'_Nu,\,u\rangle\rangle$

\vspace{1ex}

\noindent for all $u\in{\cal E}^{p,\,q}(X)$. This proves (\ref{eqn:Delta'Np_Delta'N}).

 $(ii)$\, Thanks to (\ref{eqn:Delta'Np_Delta'N}) and to G$\mathring{a}$rding's inequality satisfied by the elliptic differential operator $\Delta'_N + \Delta'_F$, there exist constants $\delta'_1,\delta'_2>0$ such that

$$\langle\langle\widetilde{\Delta'}u,\,u\rangle\rangle + \delta'_1\,||u||^2 \geq (1-\varepsilon)\,\langle\langle(\Delta'_N + \Delta'_F)\,u,\,u\rangle\rangle + \delta'_1\,||u||^2 \geq (1-\varepsilon)\delta'_2\,||u||^2_1  $$

\noindent for all $u\in{\cal E}^{p,\,q}(X)$. This proves (\ref{eqn:Garding_NF}) after putting $\delta_1:=\delta_1'$ and $\delta_2:= (1-\varepsilon)\delta'_2$.

$(iii)$\, G$\mathring{a}$rding's inequality (\ref{eqn:Garding_NF}) implies the finite dimensionality of $\ker\widetilde{\Delta'}$ and the closedness of $\mbox{Im}\,\widetilde{\Delta'}$ by standard arguments (see e.g. [Dem96, 3.10. p. 18-19]). Since $\widetilde{\Delta'}$ is self-adjoint, (\ref{eqn:Garding_NF}) also implies the following $2$-space orthogonal decomposition:

\begin{equation}\label{eqn:2-space-decomp_Delta-tilde_NF}{\cal E}^{p,\,q}(X) = \ker\widetilde{\Delta'}\bigoplus\mbox{Im}\,\widetilde{\Delta'}.\end{equation}

\noindent Now, it is clear that $\mbox{Im}\,\widetilde{\Delta'}$ is contained in \begin{equation}\label{eqn:image_inclusions}(\mbox{Im}\,\partial_F + \mbox{Im}\,(\partial_{N|\ker\Delta'_F}))\oplus(\mbox{Im}\,(\partial^{\star}_N\circ p'_F) +  \mbox{Im}\,\partial_F^{\star}) \subset (\mbox{Im}\,\partial_F + \mbox{Im}\,(\partial_{N|\ker\partial_F}))\oplus(\mbox{Im}\,(\partial^{\star}_N\circ p'_F) +  \mbox{Im}\,\partial_F^{\star}),\end{equation}

\noindent where the direct sums are orthogonal because $\partial_N$ (resp. $\partial_N^{\star}$) anti-commutes with $\partial_F$ (resp. $\partial_F^{\star}$) and $\ker\Delta'_F = \ker\partial_F\cap\ker\partial_F^{\star}$. Since the direct sums in (\ref{eqn:image_inclusions}) are, in turn, contained in ${\cal E}^{p,\,q}(X)$, the conjunction of (\ref{eqn:2-space-decomp_Delta-tilde_NF}) and (\ref{eqn:image_inclusions}) implies (\ref{eqn:3-space-decomp_Delta-tilde_NF}) and also that the inclusion (\ref{eqn:image_inclusions}) is actually an equality.

$(iv)$\, If $\partial_N$ and $\partial_F^{\star}$ anti-commute, then the following orthogonality relations hold:

\begin{equation}\label{eqn:orthogonality1-2} \mbox{Im}\,\partial_F \perp \mbox{Im}\,(\partial_{N|\ker\Delta'_F})  \hspace{2ex} \mbox{and} \hspace{2ex} \mbox{Im}\,(\partial^{\star}_N\circ p'_F) \perp \mbox{Im}\,\partial_F^{\star}.\end{equation}     

\noindent Indeed, for $u$ arbitrary and $v\in\ker\Delta'_F = \ker\partial_F\cap\ker\partial_F^{\star}$, we have 

\vspace{1ex}

$\langle\langle\partial_Fu,\,\partial_Nv\rangle\rangle = \langle\langle u,\,\partial_F^{\star}\partial_Nv\rangle\rangle = -\langle\langle u,\,\partial_N(\partial_F^{\star}v)\rangle\rangle = 0$ because $\partial_F^{\star}v=0$.

\vspace{1ex}

\noindent Similarly, for $u,v$ arbitrary, we have

\vspace{1ex}

$\langle\langle\partial_N^{\star}(p_F'u),\,\partial_F^{\star}v\rangle\rangle =  \langle\langle p_F'u,\,\partial_N\partial_F^{\star}v\rangle\rangle = - \langle\langle p_F'u,\,\partial_F^{\star}\partial_Nv\rangle\rangle = 0$ 

\vspace{1ex}

\noindent because $p_F'u\in\ker\Delta'_F \perp \mbox{Im}\,\partial_F^{\star}$.

 On the other hand, it is clear that in any pre-Hilbert space, whenever the sum of two mutually orthogonal subspaces is closed, each of the two subspaces is closed. From $(iii)$ we know that $\mbox{Im}\,\widetilde{\Delta'}$ is closed in ${\cal E}^{p,\,q}(X)$ and that it splits orthogonally as $\mbox{Im}\,\widetilde{\Delta'} = (\mbox{Im}\,\partial_F + \mbox{Im}\,(\partial_{N|\ker\Delta'_F}))\oplus(\mbox{Im}\,(\partial^{\star}_N\circ p'_F) +  \mbox{Im}\,\partial_F^{\star})$, hence $\mbox{Im}\,\partial_F + \mbox{Im}\,(\partial_{N|\ker\Delta'_F})$ and $\mbox{Im}\,(\partial^{\star}_N\circ p'_F) +  \mbox{Im}\,\partial_F^{\star}$ are closed in ${\cal E}^{p,\,q}(X)$. Thanks to the orthogonality relations (\ref{eqn:orthogonality1-2}), we infer that $\mbox{Im}\,\partial_F$, $\mbox{Im}\,(\partial_{N|\ker\Delta'_F}$, $\mbox{Im}\,(\partial^{\star}_N\circ p'_F)$ and $\mbox{Im}\,\partial_F^{\star}$ are closed in ${\cal E}^{p,\,q}(X)$ under the assumption $[\partial_N,\,\partial_F^{\star}]=0$.

 The proof of the Hodge isomorphism statement (\ref{eqn:Hodge-isom_E2_NF}) uses crucially the closedness of $\mbox{Im}\,\partial_F$ and runs along the lines of the proof of the analogous Theorem \ref{The:Hodge-isom_E2} of the Fr\"olicher case. We shall therefore limit ourselves to pointing out the main steps. 

 The sum of the first two main terms on the r.h.s. of (\ref{eqn:3-space-decomp_Delta-tilde_NF}), after replacing $\mbox{Im}\,(\partial_{N|\ker\Delta'_F})$ with $\mbox{Im}\,(\partial_{N|\ker\partial_F})$, is given by

\begin{equation}\label{eqn:pre_Hodge-isom_E2_NF}\ker(p'_F\circ\partial_N)\cap\ker\partial_F = \ker\widetilde{\Delta'}\bigoplus\bigg(\mbox{Im}\,\partial_F + \mbox{Im}\,(\partial_{N|\ker\partial_F})\bigg)\end{equation}

\noindent since $\ker(p'_F\circ\partial_N)\cap\ker\partial_F$ is easily seen to be the orthogonal complement of $\mbox{Im}\,(\partial_N^{\star}\circ p'_F) + \mbox{Im}\,\partial_F^{\star}$. Therefore, we immediately obtain the isomorphism (cf. (\ref{eqn:S_def})):

\begin{eqnarray}\label{eqn:S_NF_def}\nonumber S=S^{p,\,q}_{N,\,F} : {\cal H}^{p,\,q}_{\widetilde\Delta'}(N,\,F) & \longrightarrow & \widetilde{H}^{p,\,q}_{N,\,F}(X,\,\C):=\ker(p'_F\circ\partial_N)\cap\ker\partial_F\bigg\slash \bigg(\mbox{Im}\,\partial_F + \mbox{Im}\,(\partial_{N|\ker\partial_F})\bigg),\\
  \alpha & \longmapsto & \widetilde{[\alpha]},\end{eqnarray}

\noindent where $\widetilde{[\alpha]}$ denotes the class in $\widetilde{H}^{p,\,q}_{N,\,F}(X,\,\C)$ of $\alpha\in\ker(p'_F\circ\partial_N)\cap\ker\partial_F$. 

 On the other hand, the linear map

\begin{eqnarray}\label{eqn:T_NF_def} T = T^{p,\,q}_{N,\,F} : \widetilde{H}^{p,\,q}_{N,\,F}(X,\,\C) \longrightarrow E_2^{p,\,q}(N,\,F), \hspace{3ex} \widetilde{[\alpha]}\longmapsto\bigg[[\alpha]_{\partial_F}\bigg]_{d_1},\end{eqnarray}

\noindent is seen to be well defined and an isomorphism as in the proof of the analogous Proposition \ref{Prop:Hpq-tilde_def}. The closedness of $\mbox{Im}\,\partial_F$ is a key ingredient here. For example, for $T$ to be well defined, we need every $\alpha\in\ker(p'_F\circ\partial_N)\cap\ker\partial_F$ to induce a unique class $[[\alpha]_{\partial_F}]_{d_1}$. Now, the class $[\alpha]_{\partial_F}$ is well-defined since $\alpha\in\ker\partial_F$, but we also need to have $d_1([\alpha]_{\partial_F})=0$ in $E_1^{p,\,q}(N,\,F)$, i.e. we need $\partial_N\alpha\in\mbox{Im}\,\partial_F$. However, $\partial_N\alpha\in\ker\partial_F$ and $p'_F(\partial_N\alpha)=0$ (i.e. $\partial_N\alpha\perp\ker\Delta'_F$), which amounts precisely to $\partial_N\alpha\in\mbox{Im}\,\partial_F$. If $\mbox{Im}\,\partial_F$ were not closed, this would only amount to the weaker property $\partial_N\alpha\in\overline{\mbox{Im}\,\partial_F}$.

 It is now clear that the composition of the isomorphisms (\ref{eqn:S_NF_def}) and (\ref{eqn:T_NF_def}) provides the Hodge isomorphism (\ref{eqn:Hodge-isom_E2_NF}). 

 $(v)$\, If $\Delta'=\partial\partial^{\star} + \partial^{\star}\partial : {\cal E}^k(X)\longrightarrow{\cal E}^k(X)$ is the standard $\partial$-Laplacian induced by the metric $\omega$ (where ${\cal E}^k(X)$ is the space of smooth $(k,\,0)$-forms on $X$), the usual Hodge isomorphism theorem for $\partial$ gives the isomorphism

$${\cal H}^{p+q,\,0}_{\Delta'}(X,\,\C)\longrightarrow H^{p+q,\,0}_{\partial}(X,\,\C), \hspace{3ex}  \alpha\mapsto[\alpha]_{\partial}.$$

\noindent Coupled with the Hodge isomorphism (\ref{eqn:Hodge-isom_E2_NF}), this shows (as in the proof of Lemma \ref{Lem:sufficient-cond_E2-deg}) that if the identity map induces a well-defined linear map

\begin{equation}\label{eqn:harmonic-space-map_E2}{\cal H}_{\widetilde{\Delta'}}^{p,\,q}(N,\,F) \longrightarrow {\cal H}^{p+q,\,0}_{\Delta'}(X,\,\C), \hspace{3ex} \gamma\mapsto\gamma,\end{equation}

\noindent then $E_2(N,\,F)=E_{\infty}(N,\,F)$. Indeed, if well defined, the map (\ref{eqn:harmonic-space-map_E2}) is necessarily injective, so we get an injection $E_2^{p,\,q}(N,\,F)\hookrightarrow H^{p+q,\,0}_{\partial}(X,\,\C)$ for all $p,q$. Now, since $\omega$ is a product metric, one easily checks that for any $(p,\,q)\neq(r,\,s)$ such that $p+q=r+s=k$, the images in $H^{k,\,0}_{\partial}(X,\,\C)$ of $E_2^{p,\,q}(N,\,F)$ and $E_2^{r,\,s}(N,\,F)$ intersect only at zero because ${\cal E}^{p,\,q}(X)$ is orthogonal to ${\cal E}^{r,\,s}(X)$. Thus, for every $k$, there is a linear injection

$$\bigoplus\limits_{p+q=k}E_2^{p,\,q}(N,\,F)\hookrightarrow H^{k,\,0}_{\partial}(X,\,\C).$$

\noindent This implies the degeneration at $E_2$ of the spectral sequence.

 On the other hand, $\ker\widetilde{\Delta'} = \ker\Delta'_{N,\,p'_F}\cap\ker\Delta'_F$ (cf. (\ref{eqn:Delta-tilde_def_bis_NF})). Meanwhile, the following analogue of (\ref{eqn:Delta_gamma_gamma_pre-type}) holds:

$$\langle\langle\Delta'u,\,u\rangle\rangle = \langle\langle\Delta'_Nu,\,u\rangle\rangle + \langle\langle\Delta'_Fu,\,u\rangle\rangle \hspace{3ex} \mbox{for all}\hspace{1ex} u\in{\cal E}^{p,\,q}(X),$$

\noindent so for $u\in{\cal E}^{p,\,q}(X)$ we have the equivalence: $u\in\ker\Delta' \Leftrightarrow u\in\ker\Delta'_N\cap\ker\Delta'_F$. Putting these facts together, we see that the map (\ref{eqn:harmonic-space-map_E2}) is well defined whenever the following inclusion holds

\begin{equation}\label{eqn:kernels-intersection_E2_NF}\ker\Delta'_{N,\,p'_F}\cap\ker\Delta'_F\subset\ker\Delta'_N  \hspace{3ex} \mbox{in}\hspace{1ex} {\cal E}^{p,\,q}(X).\end{equation}

\noindent This is the analogue of Lemma \ref{Lem:kernels-intersection_E2}.

  Summing up, we have just argued that if (\ref{eqn:kernels-intersection_E2_NF}) holds for all $p,q$, then $E_2(N,\,F)=E_{\infty}(N,\,F)$. Now, we claim that the following implications hold:

\begin{eqnarray}\label{eqn:implications-E_2-degeneration_NF}\nonumber[\partial_N,\,\partial_F^{\star}] = 0 & \stackrel{(a)}{\implies} & [\partial_N,\,\Delta'_F] = 0 \hspace{2ex} \mbox{and} \hspace{2ex} [\partial_N,\,p_F] = 0 \stackrel{(b)}{\implies} \\
 & & \ker\Delta'_{N,\,p'_F}\cap\ker\Delta'_F\subset\ker\Delta'_N  \hspace{3ex} \mbox{in}\hspace{1ex} {\cal E}^{p,\,q}(X) \hspace{2ex} \mbox{for all}\hspace{1ex} p,q.\end{eqnarray}

\noindent In view of the above arguments, the implications (\ref{eqn:implications-E_2-degeneration_NF}) prove $(v)$.

 To prove the first part of implication $(a)$ of (\ref{eqn:implications-E_2-degeneration_NF}), recall that $\partial_N$ and $\partial_F$ anti-commute, so if $\partial_N$ and $\partial_F^{\star}$ anti-commute as well, then $\partial_N$ commutes with $\Delta'_F$.

 To prove the second part of implication $(a)$ of (\ref{eqn:implications-E_2-degeneration_NF}), recall that $\mbox{Im}\,\partial_F$ and $\mbox{Im}\,\partial_F^{\star}$ have been proved to be {\it closed} in ${\cal E}^{p,\,q}(X)$ under the assumption $[\partial_N,\,\partial_F^{\star}] = 0$ (cf. proof of $(iv)$). This implies the orthogonal $3$-space decompostion

\begin{equation}\label{eqn:3-space_decomp_Delta'_NF}{\cal E}^{p,\,q}(X) = \ker\Delta'_F \oplus \mbox{Im}\,\partial_F \oplus \mbox{Im}\,\partial_F^{\star}.\end{equation}

\noindent Now, let $u\in{\cal E}^{p,\,q}(X)$. By (\ref{eqn:3-space_decomp_Delta'_NF}), $u$ splits uniquely and orthogonally as $u = u_0 + \partial_Fv + \partial_F^{\star}w$, with $u_0\in\ker\Delta'_F$. Thus, $p_Fu=u_0$. We get $(\partial_Np_F)\,u = \partial_Nu_0$ and, using $[\partial_N,\,\partial_F]=0$ and $[\partial_N,\,\partial_F^{\star}]=0$, we also get

\begin{equation}\label{eqn:p_Fdel_N_commutation1}(p_F\partial_N)\,u = p_F(\partial_Nu_0) - p_F(\partial_F\partial_Nv) - p_F(\partial_F^{\star}\partial_Nw) = p_F(\partial_Nu_0),\end{equation}

\noindent since $p_F\circ\partial_F=0$ and $p_F\circ\partial_F^{\star}=0$ thanks to the orthogonal splitting (\ref{eqn:3-space_decomp_Delta'_NF}). Now, since $\Delta'_Fu_0=0$ and $\Delta'_F$ commutes with $\partial_N$, we get $\Delta'_F(\partial_Nu_0) = 0$, i.e. $p_F(\partial_Nu_0) = \partial_Nu_0$. From this and (\ref{eqn:p_Fdel_N_commutation1}), we get $(p_F\partial_N)\,u = \partial_Nu_0$. Therefore, $(\partial_Np_F)\,u =(p_F\partial_N)\,u$.

 To prove implication $(b)$ of (\ref{eqn:implications-E_2-degeneration_NF}), we run the analogue of the proof of the last statement in $(i)$ of Theorem \ref{The:commutation_E2-commutation_Rbar}. Let $u\in\ker\Delta'_{N,\,p'_F}\cap\ker\Delta'_F$. Since the assumption $\partial_Np_F = p_F\partial_N$ implies $\Delta'_{N,\,p'_F} = \Delta'_N\circ p_F$, we get the second identity below:

$$0 = \langle\langle\Delta'_{N,\,p'_F}u,\,u\rangle\rangle = \langle\langle\Delta'_N(p'_Fu),\,u\rangle\rangle = \langle\langle\Delta'_Nu,\,u\rangle\rangle,$$

\noindent where we have used the assumption $u\in\ker\Delta'_{N,\,p'_F}$ to get the first identity and the assumption $u\in\ker\Delta'_F$ to get the last identity. Consequently, $\Delta'_Nu = 0$, i.e. $u\in\ker\Delta'_N$.  \hfill $\Box$

\vspace{2ex}

\begin{Rem}\label{Rem:1st-hypothesis_enough} The hypothesis (\ref{eqn:Garding_sufficient-cond}) suffices to prove part $(v)$ of Theorem \ref{The:Garding_sufficient-cond}, but the anti-commutation hypothesis $[\partial_N,\,\partial_F^{\star}]=0$ is needed to get closedness of $\mbox{Im}\,\partial_F$ and the Hodge isomorphism (\ref{eqn:Hodge-isom_E2_NF}) in part $(iv)$.

\end{Rem}

\noindent {\it Proof.} From $(i)$, we get $(1-\varepsilon)\,(\Delta'_N + \Delta'_F) \leq \widetilde{\Delta} = \Delta'_{N,\,p'_F} + \Delta'_F \leq \Delta'_N + \Delta'_F,$ hence $\ker(\Delta'_{N,\,p'_F} + \Delta'_F) = \ker(\Delta'_N + \Delta'_F)$. This means that   

 $$\ker\Delta'_{N,\,p'_F}\cap\Delta'_F = \ker\Delta'_N\cap\ker\Delta'_F \subset \ker\Delta'_N$$

\noindent which proves (\ref{eqn:kernels-intersection_E2_NF}) in every ${\cal E}^{p,\,q}(X$.  \hfill $\Box$

\section{The anti-commutation of $\partial_N$ and $\partial_F^{\star}$}\label{section:anti-commutation}

 We now give a sufficient condition for the crucial hypothesis $[\partial_N,\,\partial_F^{\star}]=0$ of Theorem \ref{The:Garding_sufficient-cond} to hold.

\begin{The}\label{The:delN-delFstar_anticomm} Let $X$ be a compact complex manifold with $\mbox{dim}_{\C}X=n$ equipped with an integrable holomorphic almost product structure $(N,\,F)$ and with a product Hermitian metric $\omega$. 

 If $\omega$ is a bundle-like metric (cf. Definition \ref{Def:bundle-like-metric}) and if $\partial_N\omega_N=0$ (i.e. $\omega_N$ is K\"ahler in the $N$-directions), then

\begin{equation}\label{eqn:delN-delFstar_anticomm} [\partial_N,\,\partial_F^{\star}]=0.   \end{equation}

\end{The}

 The proof will proceed in local coordinates along the lines of the proof of the standard Hermitian commutation relations (\ref{eqn:standard-comm-rel}) (cf. [Dem97, VII, $\S.1$]) adapted to the foliated situation. We start with the following calculation.

\begin{Lem}\label{Lem:dz_l-adjoint} Let $z_1, \dots , z_r, z_{r+1}, \dots , z_n$ be local product coordinates and let $\omega$ be a product metric for $(N,\,F)$ (cf. $\S.$\ref{section:introduction2} for the terminology). Then, the formal adjoint of $dz_l\wedge\cdot$ is given by

\begin{equation}\label{eqn:dz_l-adjoint} (dz_l\wedge\cdot)^{\star} = \sum\limits_{k=1}^r\omega^{l\bar{k}}\,\frac{\partial}{\partial z_k}\lrcorner\cdot \hspace{2ex} \mbox{if}\hspace{1ex} l\in\{1,\dots , r\} \hspace{2ex} \mbox{and} \hspace{2ex}  (dz_l\wedge\cdot)^{\star} = \sum\limits_{k=r+1}^n\omega^{l\bar{k}}\,\frac{\partial}{\partial z_k}\lrcorner\cdot \hspace{2ex} \mbox{if}\hspace{1ex} l\in\{r+1,\dots , r\}.\end{equation}

\end{Lem}

\noindent {\it Proof.} For any vector field $\xi$ of type $(1,\,0)$, the formal adjoint w.r.t. $\omega$ of the contraction by $\xi$ is easily checked to be given by the formula

\begin{equation}\label{eqn:xi-contraction_adjoint} (\xi\lrcorner\,\,\cdot)^{\star} = i\,\overline{\alpha}_{\xi}\wedge\cdot \end{equation}

\noindent where $\alpha_{\xi}$ is the smooth $(0,\,1)$-form defined by $\alpha_{\xi}:=\xi\lrcorner\omega$. Indeed, (\ref{eqn:xi-contraction_adjoint}) is a pointwise formula since it involves only operators of order zero. Fixing an arbitrary point $x$ and choosing local coordinates about $x$ in which $\omega(x)$ is given by the identity matrix, we have $(\partial/\partial z_l\lrcorner\cdot)^{\star} = dz_l\wedge\cdot$ for all $l$. This implies (\ref{eqn:xi-contraction_adjoint}) at $x$.

 In our case, $\omega$ is given by (\ref{eqn:bundle-like-metric_def}), so $\alpha_{\frac{\partial}{\partial z_l}} = \frac{\partial}{\partial z_l}\lrcorner\omega$ is given by

$$\sum\limits_{k=1}^r\omega_{l\bar{k}}\,id\bar{z}_k \hspace{2ex} \mbox{if}\hspace{2ex} l\in\{1,\dots , r\} \hspace{2ex} \mbox{and by} \hspace{2ex} \sum\limits_{k=r+1}^n\omega_{l\bar{k}}\,id\bar{z}_k \hspace{2ex} \mbox{if}\hspace{2ex} l\in\{r+1,\dots , n\}.$$

\noindent Applying (\ref{eqn:xi-contraction_adjoint}) with $\xi= \frac{\partial}{\partial z_l}$, we get

$$\bigg(\frac{\partial}{\partial z_l}\lrcorner\cdot\bigg)^{\star} = \sum\limits_{k=1}^r\omega_{k\bar{l}}\,dz_k\wedge\cdot \hspace{1ex} \mbox{if}\hspace{2ex} l\in\{1,\dots , r\} \hspace{2ex} \mbox{and} \hspace{2ex} \bigg(\frac{\partial}{\partial z_l}\lrcorner\cdot\bigg)^{\star} = \sum\limits_{k=r+1}^n\omega_{k\bar{l}}\,dz_k\wedge\cdot \hspace{1ex} \mbox{if}\hspace{2ex} l\in\{r+1,\dots , n\}.$$

\noindent Thus, $\frac{\partial}{\partial z_l}\lrcorner\cdot = \sum\limits_{k=1}^r\omega_{l\bar{k}}\,(dz_k\wedge\cdot)^{\star}$ if $l\in\{1,\dots , r\}$ and  $\frac{\partial}{\partial z_l}\lrcorner\cdot = \sum\limits_{k=r+1}^n\omega_{l\bar{k}}\,(dz_k\wedge\cdot)^{\star}$ if $l\in\{r+1,\dots , n\}$.

\noindent These identities are equivalent to (\ref{eqn:dz_l-adjoint}).  \hfill $\Box$

\begin{Lem}\label{Lem:Lambda_N-del_F} Let $\Lambda_N:=(\omega_N\wedge\cdot)^{\star}$ and $\Lambda_F:=(\omega_F\wedge\cdot)^{\star}$.

 If $\omega$ is a bundle-like metric in the sense of Definition \ref{Def:bundle-like-metric}, then $[\Lambda_N,\,\partial_F]=0$.

\end{Lem}

\noindent {\it Proof.} Using formula (\ref{eqn:dz_l-adjoint}), we get

$$\Lambda_N = (\sum\limits_{j,\,k=1}^r\omega_{j\bar{k}}\,idz_j\wedge d\bar{z}_k\wedge\cdot)^{\star} = -i\,\sum\limits_{j,\,k=1}^r \omega_{k\bar{j}} \sum\limits_{s=1}^r\omega^{s\bar{k}}\,\frac{\partial}{\partial\bar{z}_s}\lrcorner\bigg(\sum\limits_{l=1}^r\omega^{j\bar{l}}\,\frac{\partial}{\partial z_l}\lrcorner\cdot\bigg) = -i\,\sum\limits_{l,\,s=1}^ra_{s\bar{l}}\,\frac{\partial}{\partial\bar{z}_s}\lrcorner\frac{\partial}{\partial z_l}\lrcorner\cdot,$$

\noindent where we have put $a_{s\bar{l}} := \sum\limits_{j,\,k=1}^r\omega_{k\bar{j}}\,\omega^{s\bar{k}}\,\omega^{j\bar{l}}$. From this, we get for every form $u$:

\begin{eqnarray}\nonumber [\Lambda_N,\,\partial_F]\,u & = & -i\,\sum\limits_{l,\,s=1}^ra_{s\bar{l}}\,\frac{\partial}{\partial\bar{z}_s}\lrcorner\bigg[\frac{\partial}{\partial z_l}\lrcorner \sum\limits_{\tau=r+1}^ndz_{\tau}\wedge\frac{\partial u}{\partial z_{\tau}}\bigg] + i\, \sum\limits_{\tau=r+1}^ndz_{\tau}\wedge\frac{\partial}{\partial z_{\tau}}\bigg(\sum\limits_{l,\,s=1}^ra_{s\bar{l}}\,\frac{\partial}{\partial\bar{z}_s}\lrcorner\frac{\partial}{\partial z_l}\lrcorner u\bigg).\end{eqnarray}

\noindent Since $\frac{\partial}{\partial z_l}\lrcorner \sum\limits_{\tau=r+1}^ndz_{\tau}\wedge\frac{\partial u}{\partial z_{\tau}} = \sum\limits_{\tau=r+1}^n\bigg[\delta_{l\tau}\,\frac{\partial u}{\partial z_{\tau}} - dz_{\tau}\wedge\bigg(\frac{\partial}{\partial z_l}\lrcorner\frac{\partial u}{\partial z_{\tau}}\bigg)\bigg]$ and since $\delta_{l\tau} = 0$ for all $l\in\{1,\dots , r\}$ and all $\tau\in\{r+1, \dots , n\}$, the above identity reads

\begin{eqnarray}\label{eqn:Lambda_N-del_F_1}\nonumber [\Lambda_N,\,\partial_F]\,u & = & -i\,\sum\limits_{l,\,s=1}^ra_{s\bar{l}}\,\sum\limits_{\tau=r+1}^ndz_{\tau}\wedge\bigg(\frac{\partial}{\partial\bar{z}_s}\lrcorner\frac{\partial}{\partial z_l}\lrcorner\frac{\partial u}{\partial z_{\tau}}\bigg) + i\,\sum\limits_{l,\,s=1}^r\sum\limits_{\tau=r+1}^na_{s\bar{l}}\,dz_{\tau}\wedge\bigg(\frac{\partial}{\partial\bar{z}_s}\lrcorner\frac{\partial}{\partial z_l}\lrcorner\frac{\partial u}{\partial z_{\tau}}\bigg) \\
  & + & i\,\sum\limits_{l,\,s=1}^r\sum\limits_{\tau=r+1}^n \frac{\partial a_{s\bar{l}}}{\partial z_{\tau}}\,dz_{\tau}\wedge\bigg(\frac{\partial}{\partial\bar{z}_s}\lrcorner\frac{\partial}{\partial z_l}\lrcorner u\bigg).\end{eqnarray}

\noindent Hence $[\Lambda_N,\,\partial_F]\,u = i\,\sum\limits_{l,\,s=1}^r\sum\limits_{\tau=r+1}^n \frac{\partial a_{s\bar{l}}}{\partial z_{\tau}}\,dz_{\tau}\wedge\bigg(\frac{\partial}{\partial\bar{z}_s}\lrcorner\frac{\partial}{\partial z_l}\lrcorner u\bigg)$ since the top line on the r.h.s. of (\ref{eqn:Lambda_N-del_F_1}) clearly vanishes. 

 On the other hand, if $\omega$ is a bundle-like metric in the sense of Definition \ref{Def:bundle-like-metric}, the coefficients $\omega_{j\bar{k}}$ with $j,k\in\{1, \dots , r\}$ depend only on the variables $z'=(z_1,\dots , z_r)$, so the derivatives of all the quantities $\omega_{k\bar{j}},\,\omega^{s\bar{k}},\,\omega^{j\bar{l}}$ with indices $j,k,s,l\in\{1,\dots , r\}$ w.r.t. any of the variables $z_{\tau}$ with $\tau\in\{r+1, \dots , n\}$ vanish. This implies that $\frac{\partial a_{s\bar{l}}}{\partial z_{\tau}} = 0$ for all $l, s\in\{1,\dots , r\}$ and all $\tau\in\{r+1, \dots , n\}$, which further implies that $[\Lambda_N,\,\partial_F] = 0$ thanks to the above formula.  \hfill $\Box$

\vspace{2ex}

 The next observation is that the standard Hermitian commutatin relations (\ref{eqn:standard-comm-rel}) are still valid for $\partial_N$ and $\partial_F$. 

\begin{Lem}\label{Lem:N-F_comm-rel} Let $\tau_N:=[\Lambda_N,\,\partial_N\omega_N\wedge\cdot]$ denote the $N$-directional torsion operator of $\omega$. If $\omega$ is a bundle-like metric in the sense of Definition \ref{Def:bundle-like-metric}, the following identity holds:

\begin{equation}\label{eqn:N-F_comm-rel} [\Lambda_N,\,\partial_N] = i\,(\bar\partial_N^{\star} + \bar\tau_N^{\star}).\end{equation}

\end{Lem}

\noindent {\it Proof.} We shall only prove this identity in the case when $\partial_N\omega_N=0$ (hence $\tau_N=0$) since this is the only case we need in this paper. Fix an arbitrary point $x_0\in X$ and choose local holomorphic product coordinates about $x_0$ such that $\omega_{j\bar{k}}(x_0) = \delta_{jk}$ for all $j.k\in\{1,\dots ,\,r\}$. Thanks to the assumption $\partial_N\omega_N=0$, we can choose these coordinates to also satisfy the property:

\begin{equation}\label{eqn:coordinates_1st-derivatives}\frac{\partial\omega_{j\bar{k}}}{\partial z_l}(x_0) = 0  \hspace{3ex} \mbox{for all}\hspace{1ex} j,k,l\in\{1,\dots , r\}. \end{equation} 

 We now use the notation and the computations in the proof of Lemma \ref{Lem:Lambda_N-del_F} to deduce the following identity (which is the analogue of (\ref{eqn:Lambda_N-del_F_1}) but with an extra term) for every form $u$:

\begin{eqnarray}\label{eqn:Lambda_N-del_N_1} [\Lambda_N,\,\partial_N]\,u = -i\, \sum\limits_{l,\,s=1}^r a_{s\bar{l}} \frac{\partial}{\partial\bar{z}_s}\lrcorner\frac{\partial u}{\partial z_l} + i\,\sum\limits_{l,\,s=1}^r\sum\limits_{\tau=1}^r \frac{\partial a_{s\bar{l}}}{\partial z_{\tau}}\,dz_{\tau}\wedge\bigg(\frac{\partial}{\partial\bar{z}_s}\lrcorner\frac{\partial}{\partial z_l}\lrcorner u\bigg).\end{eqnarray}

\noindent The extra term comes from the fact that $\tau$ ranges now from $1$ to $r$, so $(\partial/\partial z_l)\lrcorner dz_{\tau} = \delta_{l\tau}$ and $\delta_{l\tau} = 1$ when $l=\tau$ (a situation which, unlike in the case of the proof of Lemma \ref{Lem:Lambda_N-del_F}, can occur).

\noindent With our choice of coordinates, we have $a_{s\bar{l}}(x_0) = \delta_{sl}$ and also, thanks to (\ref{eqn:coordinates_1st-derivatives}), $(\partial a_{s\bar{l}}/\partial z_{\tau})(x_0) = 0$ for all $l,s,\tau\in\{1, \dots , r\}$. We conclude from (\ref{eqn:Lambda_N-del_N_1}) that

\begin{eqnarray}\label{eqn:Lambda_N-del_N_2} [\Lambda_N,\,\partial_N]\,u = -i\, \sum\limits_{l=1}^r \frac{\partial}{\partial\bar{z}_l}\lrcorner\frac{\partial u}{\partial z_l} \hspace{2ex} \mbox{at} \hspace{1ex} x_0  \hspace{2ex} \mbox{if} \hspace{1ex} \partial_N\omega_N = 0.\end{eqnarray}

 Meanwhile, a simple integration by parts shows that the adjoint of $\partial/\partial z_l$ acting on functions is

\begin{equation}\label{eqn:adjoint_d-dz_l} \bigg(\frac{\partial}{\partial z_l}\bigg)^{\star}g = -(\det\omega)\,\frac{\partial g}{\partial\bar{z}_l} - \frac{\partial(\det\omega)}{\partial\bar{z}_l}\,g  \hspace{2ex} \mbox{for every function} \hspace{1ex} g \hspace{2ex} \mbox{and every} \hspace{1ex} l\in\{1,\dots , n\},\end{equation}

\noindent where $\det\omega$ stands for the determinant of the matrix $(\omega_{j\bar{k}})_{j,k=1,\dots , n}$. Using (\ref{eqn:dz_l-adjoint}) and (\ref{eqn:adjoint_d-dz_l}), we get

\begin{eqnarray}\label{eqn:d-bar_N-star1} \bar\partial_N^{\star} & = & \bigg(\sum\limits_{l=1}^r d\bar{z}_l\wedge\frac{\partial}{\partial\bar{z}_l}\bigg)^{\star} = \sum\limits_{l=1}^r \bigg(\frac{\partial}{\partial\bar{z}_l}\bigg)^{\star}\bigg(\sum\limits_{k=1}^r \omega^{k\bar{l}}\,\frac{\partial}{\partial\bar{z}_k}\lrcorner\cdot\bigg) \\
 \nonumber   & = & - (\det\omega)\,\sum\limits_{k,l=1}^r \frac{\partial\omega^{k\bar{l}}}{\partial z_l}\bigg(\frac{\partial}{\partial\bar{z}_k}\lrcorner\cdot\bigg) - (\det\omega)\,\sum\limits_{k,l=1}^r \omega^{k\bar{l}}\,\frac{\partial}{\partial\bar{z}_k}\lrcorner\frac{\partial}{\partial z_l}  - \sum\limits_{k,l=1}^r \frac{\partial(\det\omega)}{\partial z_l}\,\omega^{k\bar{l}}\,\frac{\partial}{\partial\bar{z}_k}\lrcorner\cdot \end{eqnarray}

\noindent Now, $\det\omega = \det(\omega_{j\bar{k}})_{j,\,k=1,\dots , r}\,\det(\omega_{j\bar{k}})_{j,\,k=r+1,\dots , n}$ since $\omega$ is a product metric. Meanwhile, $\partial\omega_{j\bar{k}}/\partial z_l = 0$ for all $j,k\in\{r+1,\dots , n\}$ and all $l\in\{1,\dots , r\}$ by the {\it bundle-like} assumption on $\omega$. Combined with (\ref{eqn:coordinates_1st-derivatives}), these properties imply that $\partial(\det\omega)/\partial z_l (x_0) =0$ for all $l\in\{1,\dots , r\}$. Since $\partial\omega^{k\bar{l}}/\partial z_l (x_0) = 0$ for all $k,l\in\{1,\dots , r\}$ by  (\ref{eqn:coordinates_1st-derivatives}), (\ref{eqn:d-bar_N-star1}) translates to the following identity at $x_0$ for every form $u$:

\begin{eqnarray}\label{eqn:d-bar_N-star2} \bar\partial_N^{\star}u = - \sum\limits_{l=1}^r \frac{\partial}{\partial\bar{z}_l}\lrcorner\frac{\partial u}{\partial z_l}   \hspace{2ex} \mbox{at} \hspace{1ex} x_0  \hspace{2ex} \mbox{if} \hspace{1ex} \partial_N\omega_N = 0 \hspace{2ex} \mbox{and if} \hspace{1ex} \omega \hspace{1ex} \mbox{is a bundle-like metric}.\end{eqnarray}

 The conjunction of (\ref{eqn:Lambda_N-del_N_2}) and (\ref{eqn:d-bar_N-star2}) proves (\ref{eqn:N-F_comm-rel}) when $\partial_N\omega_N = 0$ and $\omega$ is a bundle-like metric. \hfill  $\Box$

\begin{Lem}\label{Lem:Lambda_N-del_F-delbar_N} We always have $[\partial_F,\,\bar\partial_N]=0$.

\end{Lem}

\noindent {\it Proof.} $(i)$\, If $u$ is an arbitrary differential form, we have
 \begin{eqnarray}\nonumber [\partial_F,\,\bar\partial_N]\,u & = & \sum\limits_{l=r+1}^ndz_l\wedge\frac{\partial}{\partial z_l}\bigg(\sum\limits_{j=1}^rd\bar{z}_j\wedge\frac{\partial u}{\partial\bar{z}_j}\bigg) + \sum\limits_{j=1}^rd\bar{z}_j\wedge\frac{\partial}{\partial\bar{z}_j}\bigg(\sum\limits_{l=r+1}^ndz_l\wedge\frac{\partial u}{\partial z_l}\bigg) \\
\nonumber & = & \sum\limits_{j=1}^r\sum\limits_{l=r+1}^n dz_l\wedge d\bar{z}_j\wedge \frac{\partial^2u}{\partial z_l \partial\bar{z}_j} + \sum\limits_{j=1}^r\sum\limits_{l=r+1}^n d\bar{z}_j\wedge dz_l\wedge \frac{\partial^2u}{\partial z_l \partial\bar{z}_j} = 0.\end{eqnarray}

\hfill $\Box$

\vspace{3ex}

 \noindent {\it End of proof of Theorem \ref{The:delN-delFstar_anticomm}.} We shall actually prove the identity

\begin{equation}\label{eqn:delF-delNstar_anticomm} [\partial_F,\,\partial_N^{\star}]=0\end{equation}

\noindent under the hypotheses of Theorem \ref{The:delN-delFstar_anticomm} (which imply, in particular, that $\tau_N = 0$). Taking adjoints, (\ref{eqn:delF-delNstar_anticomm}) is seen to be equivalent to (\ref{eqn:delN-delFstar_anticomm}).

 Since $\tau_N = 0$, the conjugate of (\ref{eqn:N-F_comm-rel}) reads $\partial_N^{\star} = i\,[\Lambda_N,\,\bar\partial_N]$, so we get

$$ [\partial_F,\,\partial_N^{\star}] = i\,[\partial_F,\,[\Lambda_N,\,\bar\partial_N]] = i\,[\Lambda_N,\,[\bar\partial_N,\,\partial_F]] +  i\,[\bar\partial_N,\,[\partial_F,\,\Lambda_N]].$$

 \noindent Since $[\bar\partial_N,\,\partial_F] = 0$ by Lemma \ref{Lem:Lambda_N-del_F-delbar_N} and $[\partial_F,\,\Lambda_N] = 0$ by Lemma \ref{Lem:Lambda_N-del_F}, we get $[\partial_F,\,\partial_N^{\star}] = 0$.  \hfill $\Box$

\section{Appendix: Commutation relations}\label{Appendix}

Let $(X,\,\omega)$ be a compact complex Hermitian manifold. Recall the following standard Hermitian commutation relations ([Dem84], see also [Dem97, VII, $\S.1$]):

\begin{eqnarray}\label{eqn:standard-comm-rel}\nonumber &  & (i)\,\,(\partial + \tau)^{\star} = i\,[\Lambda,\,\bar\partial];  \hspace{3ex} (ii)\,\,(\bar\partial + \bar\tau)^{\star} = - i\,[\Lambda,\,\partial]; \\
&  & (iii)\,\, \partial + \tau = -i\,[\bar\partial^{\star},\,L]; \hspace{3ex} (iv)\,\,
\bar\partial + \bar\tau = i\,[\partial^{\star},\,L],\end{eqnarray}

\noindent where the upper symbol $\star$ stands for the formal adjoint w.r.t. the $L^2$ inner product induced by $\omega$, $L=L_{\omega}:=\omega\wedge\cdot$ is the Lefschetz operator of multiplication by $\omega$, $\Lambda=\Lambda_{\omega}:=L^{\star}$ and $\tau:=[\Lambda,\,\partial\omega\wedge\cdot]$ is the torsion operator (of order zero and type $(1,\,0)$) associated with the metric $\omega$.

 Again following [Dem97, VII, $\S.1$], recall that the commutation relations $(1)$ immediately induce via the Jacobi identity the Bochner-Kodaira-Nakano-type identity

\begin{equation}\label{eqn:BKN_demailly1}\Delta'' = \Delta' + [\partial,\,\tau^{\star}] - [\bar\partial,\,\bar{\tau}^{\star}]\end{equation} 

\noindent relating the $\bar\partial$-Laplacian $\Delta''=[\bar\partial,\,\bar\partial^{\star}]=\bar\partial\bar\partial^{\star} + \bar\partial^{\star} \bar\partial$ and the $\partial$-Laplacian $\Delta'=[\partial,\,\partial^{\star}]=\partial\partial^{\star} + \partial^{\star}\partial$. This, in turn, induces the following Bochner-Kodaira-Nakano-type identity (cf. [Dem84]) in which the first-order terms have been absorbed in the twisted Laplace-type operator $\Delta'_{\tau}:=[\partial+\tau,\, (\partial+\tau)^{\star}]$:

 \begin{equation}\label{eqn:BKN_demailly2}\Delta'' = \Delta'_{\tau} + T_{\omega},\end{equation} 

\noindent where $T_{\omega}:=\bigg[\Lambda,\,[\Lambda,\,\frac{i}{2}\,\partial\bar\partial\omega]\bigg] - [\partial\omega\wedge\cdot,\,(\partial\omega\wedge\cdot)^{\star}]$ is a zeroth order operator of type $(0,\,0)$ associated with the torsion of $\omega$. Formula (\ref{eqn:BKN_demailly2}) is obtained from (\ref{eqn:BKN_demailly1}) via the following identities (cf. [Dem84] or [Dem97, VII, $\S.1$]) which have an interest of their own:

\begin{eqnarray}\label{eqn:BKN_demailly_auxiliary}\nonumber &  & (i)\,\,[L,\,\tau] = 3\,\partial\omega\wedge\cdot,  \hspace{3ex} (ii)\,\, [\Lambda,\,\tau] = -2i\,\bar{\tau}^{\star},\\
 &  & (iii)\,\, [\partial,\,\bar{\tau}^{\star}] = - [\partial,\,\bar\partial^{\star}] = [\tau,\,\bar\partial^{\star}],  \hspace{3ex}  (iv)\,\, -[\bar\partial,\,\bar\tau^{\star}] = [\tau,\, (\partial+\tau)^{\star}] + T_{\omega}.\end{eqnarray}

\noindent Note that $(iii)$ yields, in particular, that $\partial$ and $\bar\partial^{\star} + \bar\tau^{\star}$ anti-commute, hence by conjugation, $\bar\partial$ and $\partial^{\star} + \tau^{\star}$ anti-commute, i.e.

\begin{equation}\label{eqn:basic-anticommutation}[\partial,\,\bar\partial^{\star} + \bar\tau^{\star}] = 0 \hspace{2ex} \mbox{and} \hspace{2ex} [\bar\partial,\,\partial^{\star} + \tau^{\star}] = 0.\end{equation}

\subsection{Bochner-Kodaira-Nakano formula for the Witten Laplacians}\label{subsection:BKN_Witten}

 The computations in this subsection are motivated by an attempt at widening the spectral gap of $\Delta' + \Delta''$ to make it exceed the triple of the torsion upper bound so that hypothesis (\ref{eqn:torsion-smallness_spectral-gap_introd}) in Theorem \ref{The:E_2_small-torsionSKT_introd} may be satisfied. If this were possible, we would get the degeneration at $E_2$ of the Fr\"olicher spectral sequence under the mere assumption that an SKT metric exist on $X$. The idea was to use the Witten twisting of the connection $\partial + \bar\partial$ which produces an isomorphic cohomology (hence only changes the spectral sequence by an isomorphism) without changing the metric (and thus keeping the same torsion). Although our original goal could not be achieved by this method, we believe that Theorem \ref{The:BKN_Witten2} is of independent interest. For this reason, we spell out these calculations here.  

 Let $(X,\,\omega)$ still stand for a compact Hermitian manifold with $\mbox{dim}_{\C}X=n$. We start by recalling a few standard facts. For a fixed $C^{\infty}$ function $\varphi:X\rightarrow\R$, consider the twisted $\partial$ and $\bar\partial$ operators in the sense of Witten (cf. [Wit82]):

\begin{equation}\label{def:Witten-twisting_def}\partial_{\varphi}\alpha:=e^{\varphi}\partial(e^{-\varphi}\alpha)  \hspace{3ex} \mbox{and} \hspace{3ex} \bar\partial_{\varphi}\alpha:=e^{\varphi}\bar\partial(e^{-\varphi}\alpha),\end{equation}

\noindent where $\alpha$ is an arbitrary smooth form on $X$. The operators $\partial_{\varphi}$ and $\bar\partial_{\varphi}$ are clearly integrable and anti-commute with each other, i.e.

\begin{equation}\label{eqn:integrability_anti-commutation}\partial_{\varphi}^2=0, \hspace{2ex} \bar\partial_{\varphi}^2=0 \hspace{2ex} \mbox{and} \hspace{2ex} \partial_{\varphi}\bar\partial_{\varphi} + \bar\partial_{\varphi}\partial_{\varphi} = 0.\end{equation}

\noindent It is clear that the multiplication by $e^{-\varphi}$ defines isomorphisms between $\ker\partial_{\varphi}$ and $\ker\partial$, between $\mbox{Im}\partial_{\varphi}$ and $\mbox{Im}\partial$, between $\ker\bar\partial_{\varphi}$ and $\ker\bar\partial$, as well as between $\mbox{Im}\bar\partial_{\varphi}$ and $\mbox{Im}\bar\partial$. Therefore, we get isomorphisms in cohomology for every $p,q\in\{0, \dots , n\}$:

\begin{eqnarray}\label{eqn:Witten_cohomology_isomorphisms} H^{p,\,q}_{\partial_{\varphi}}(X,\,\C)\longrightarrow H^{p,\,q}_{\partial}(X,\,\C) & \mbox{and} &  H^{p,\,q}_{\bar\partial_{\varphi}}(X,\,\C)\longrightarrow H^{p,\,q}_{\bar\partial}(X,\,\C)\\
 \nonumber [\alpha]_{\partial_{\varphi}}\longmapsto [e^{-\varphi}\alpha]_{\partial}  & \mbox{and} &  [\alpha]_{\bar\partial_{\varphi}}\longmapsto [e^{-\varphi}\alpha]_{\bar\partial},\end{eqnarray}

\noindent where $H^{p,\,q}_{\partial_{\varphi}}(X,\,\C):=\ker\partial_{\varphi}/\mbox{Im}\partial_{\varphi}$ at the level of $(p,\,q)$-forms, while $[\alpha]_{\partial_{\varphi}}$ stands for the $\partial_{\varphi}$-cohomology class of $\alpha$ and the analogous notation is used for the $\bar\partial_{\varphi}$, $\partial$ and $\bar\partial$-cohomologies. This illustrates one of the main ideas in [Wit82]: the twisting (\ref{def:Witten-twisting_def}) of the connection does not change the cohomology.

 On the other hand, we have

 \begin{equation}\label{def:Witten-twisting_expanded}\partial_{\varphi} = \partial - \partial\varphi\wedge\cdot \hspace{3ex} \mbox{and} \hspace{3ex} \bar\partial_{\varphi} = \bar\partial - \bar\partial\varphi\wedge\cdot.\end{equation}

\noindent Taking adjoints w.r.t. the $L^2$ inner product defined by the metric $\omega$, we get

\begin{equation}\label{def:Witten-twisting_adjoints}\partial_{\varphi}^{\star} = \partial^{\star} - i\,\xi_{\varphi}\lrcorner\cdot \hspace{3ex} \mbox{and} \hspace{3ex} \bar\partial_{\varphi}^{\star} = \bar\partial^{\star} + i\,\bar\xi_{\varphi}\lrcorner\cdot,\end{equation}

\noindent where $\xi_{\varphi}\in C^{\infty}(X,\,T^{1,\,0}X)$ (resp. its conjugate $\bar\xi_{\varphi}\in C^{\infty}(X,\,T^{0,\,1}X)$) is the unique vector field of type $(1,\,0)$ (resp. $(0,\,1)$) satisfying

\begin{equation}\label{eqn:xi_varphi_def}\xi_{\varphi}\lrcorner\omega = \bar\partial\varphi  \hspace{3ex} \mbox{resp.} \hspace{3ex} \bar\xi_{\varphi}\lrcorner\omega = \partial\varphi.\end{equation}

\noindent Indeed, the following formulae (the first two of which have been used to infer (\ref{def:Witten-twisting_adjoints}) from (\ref{def:Witten-twisting_expanded})), can be easily checked to hold for any smooth differential forms $u$ and $v$ of any degrees. (The use of complex Lie derivatives of the shape (\ref{eqn:Lie-deriv_def}) was suggested in e.g. [Sar78] and variants have been recently used in e.g. [LRY15].)

\begin{Formulae}\label{Formulae:adjoints}

 $(a)$\,\, $(\partial\varphi\wedge\cdot)^{\star} = i\,\xi_{\varphi}\lrcorner\cdot,$  \hspace{3ex} $(b)$\,\, $(\bar\partial\varphi\wedge\cdot)^{\star} = -i\,\bar\xi_{\varphi}\lrcorner\cdot,$  

\vspace{1ex}

$(c)$\,\,$\partial^{\star}(\varphi\,u) = \varphi\,\partial^{\star}u - (\partial\varphi\wedge\cdot)^{\star}u = \varphi\,\partial^{\star}u - i\,\xi_{\varphi}\lrcorner u,$

\vspace{1ex}

$(d)$\,\,$\bar\partial^{\star}(\varphi\,u) = \varphi\, \bar\partial^{\star}u - (\bar\partial\varphi\wedge\cdot)^{\star}u = \varphi\, \bar\partial^{\star}u + i\,\bar\xi_{\varphi}\lrcorner u,$

 \vspace{1ex}

$(e)$\,\,$\partial(\xi_{\varphi}\lrcorner v) = L^{1,\,0}_{\xi_{\varphi}}v - \xi_{\varphi}\lrcorner\partial v$,  \hspace{3ex}  $(f)$\,\,$\bar\partial(\bar\xi_{\varphi}\lrcorner v) = L^{0,\,1}_{\bar\xi_{\varphi}}v - \bar\xi_{\varphi}\lrcorner\bar\partial v$,

\vspace{1ex}

\noindent where, by analogy with the standard Lie derivative w.r.t. a vector field in the real context, we define the $(1,\,0)$-Lie derivative w.r.t. $\xi_{\varphi}$, resp. the $(0,\,1)$-Lie derivative w.r.t. $\bar\xi_{\varphi}$, by

\begin{equation}\label{eqn:Lie-deriv_def}L^{1,\,0}_{\xi_{\varphi}}:=[\xi_{\varphi}\lrcorner\cdot,\,\partial],   \hspace{2ex} \mbox{resp.} \hspace{2ex} L^{0,\,1}_{\bar\xi_{\varphi}}:=[\bar\xi_{\varphi}\lrcorner\cdot,\,\bar\partial].\end{equation}

\end{Formulae}

 We now consider the Laplace-Beltrami operators associated with $\partial_{\varphi}$, resp. $\bar\partial_{\varphi}$, and the Hermitian metric $\omega$, i.e. the holomorphic, resp. anti-holomorphic Witten Laplacians induced by $\varphi$ and $\omega$:

\begin{equation}\label{eqn:Witten-Laplacians_def}\Delta'_{\varphi}:=[\partial_{\varphi},\,\partial_{\varphi}^{\star}] = \partial_{\varphi}\partial_{\varphi}^{\star} + \partial_{\varphi}^{\star}\partial_{\varphi} \hspace{2ex} \mbox{and} \hspace{2ex} \Delta''_{\varphi}:=[\bar\partial_{\varphi},\,\bar\partial_{\varphi}^{\star}] = \bar\partial_{\varphi}\bar\partial_{\varphi}^{\star} + \bar\partial_{\varphi}^{\star}\bar\partial_{\varphi}.\end{equation}

\begin{Lem}\label{Lem:Witten-Laplacians_formulae} The Witten Laplacians are given in terms of the standard Laplacians by the formulae

\begin{equation}\label{eqn:Witten-Laplacians_formulae}\Delta'_{\varphi} = \Delta' -i\,L^{1,\,0}_{\xi_{\varphi}} - (i\,L^{1,\,0}_{\xi_{\varphi}})^{\star} + |\partial\varphi|_{\omega}^2  \hspace{2ex} \mbox{and} \hspace{2ex} \Delta''_{\varphi} = \Delta'' + i\,L^{0,\,1}_{\bar\xi_{\varphi}} + (i\,L^{0,\,1}_{\bar\xi_{\varphi}})^{\star} + |\partial\varphi|_{\omega}^2,\end{equation}

\noindent where $|\partial\varphi|_{\omega}^2 = |\bar\partial\varphi|_{\omega}^2 = \Lambda_{\omega}(i\partial\varphi\wedge\bar\partial\varphi) = |\xi_{\varphi}|_{\omega}^2 = |\bar\xi_{\varphi}|_{\omega}^2$ is the pointwise squared norm w.r.t. $\omega$ of the $(1,\,0)$-form $\partial\varphi$ and also of the $(1,\,0)$-vector field $\xi_{\varphi}$.

\end{Lem}

\noindent {\it Proof.} The formula for $\Delta''_{\varphi}$ in (\ref{eqn:Witten-Laplacians_formulae}) is the conjugate of the one for $\Delta'_{\varphi}$, so it suffices to prove the former. Using (\ref{def:Witten-twisting_expanded}) and (\ref{def:Witten-twisting_adjoints}), for an arbitrary differential form $\alpha$ we get

\begin{eqnarray}\label{eqn:Delta'_varphi_computation1}\nonumber \Delta'_{\varphi}\alpha & = & (\partial - \partial\varphi\wedge\cdot)(\partial^{\star}\alpha - i\xi_{\varphi}\lrcorner\alpha) + (\partial^{\star} - i\xi_{\varphi}\lrcorner\cdot)(\partial\alpha - \partial\varphi\wedge\alpha)\\
    & = & \Delta' -i\,(\partial(\xi_{\varphi}\lrcorner\alpha) + \xi_{\varphi}\lrcorner\partial\alpha) - (\partial\varphi\wedge\partial^{\star}\alpha + \partial^{\star}(\partial\varphi\wedge\alpha)) + (i\,\xi_{\varphi}\lrcorner\partial\varphi)\,\alpha.\end{eqnarray}

\noindent Now, $\partial(\xi_{\varphi}\lrcorner\alpha) + \xi_{\varphi}\lrcorner\partial\alpha = [\xi_{\varphi}\lrcorner\cdot,\,\partial]\,(\alpha) = L^{1,\,0}_{\xi_{\varphi}}\,\alpha$ thanks to definition (\ref{eqn:Lie-deriv_def}). Moreover, $- (\partial\varphi\wedge\partial^{\star}\alpha + \partial^{\star}(\partial\varphi\wedge\alpha)) = -[\partial^{\star},\,\partial\varphi\wedge\cdot](\alpha) = - [i\,\xi_{\varphi}\lrcorner\cdot,\,\partial]^{\star}\,(\alpha) = -(i\,L^{1,\,0}_{\xi_{\varphi}})^{\star}\alpha$. Combined with (\ref{eqn:Delta'_varphi_computation1}), these identities prove (\ref{eqn:Witten-Laplacians_formulae}) if we can show that

\begin{equation}\label{eqn:i-xi_varphi-del_varphi}i\,\xi_{\varphi}\lrcorner\partial\varphi = |\partial\varphi|_{\omega}^2.\end{equation}

\noindent This is easily proved in local coordinates $z_1,\dots, z_n$. If we denote

$$\omega = \sum\limits_{\lambda,\,\mu=1}^n\omega_{\lambda\bar{\mu}}\,idz_{\lambda}\wedge d\bar{z}_{\mu}  \hspace{2ex} \mbox{and} \hspace{2ex} \xi_{\varphi} = \sum\limits_{l=1}^n\xi_l\,\frac{\partial}{\partial z_l},$$

\noindent the identity $\xi_{\varphi}\lrcorner\omega = \bar\partial\varphi$ defining $\xi_{\varphi}$ (cf. (\ref{eqn:xi_varphi_def})) is equivalent to 

\begin{equation}\label{eqn:xi_varphi_def_coordinates}\sum\limits_{\mu=1}^n\bigg(\sum\limits_{l=1}^ni\,\xi_l\,\omega_{l\bar{\mu}}\bigg)\,d\bar{z}_{\mu} = \sum\limits_{\mu=1}^n\frac{\partial\varphi}{\partial\bar{z}_{\mu}}\,d\bar{z}_{\mu}, \hspace{2ex} \mbox{i.e. to} \hspace{2ex} i\,\sum\limits_{l=1}^n\xi_l\,\omega_{l\bar{\mu}} = \frac{\partial\varphi}{\partial\bar{z}_{\mu}} \hspace{2ex} \mbox{for all} \hspace{2ex} \mu=1,\dots , n.\end{equation}

\noindent Now, $i\,\xi_{\varphi}\lrcorner\partial\varphi = i\,\sum\limits_{l=1}^n\xi_l\,\frac{\partial\varphi}{\partial z_l}$, so using the identity $\frac{\partial\varphi}{\partial z_l} = -i\,\sum\limits_{r=1}^n\bar\xi_r\,\omega_{l\bar{r}}$ (the conjugate of (\ref{eqn:xi_varphi_def_coordinates})), we get

\begin{equation}\label{eqn:i-xi_varphi-del_varphi_norm}i\,\xi_{\varphi}\lrcorner\partial\varphi = \sum\limits_{l,\,r=1}^n\xi_l\,\bar\xi_r\,\omega_{l\bar{r}} = |\xi_{\varphi}|_{\omega}^2 = -i\,\sum\limits_{r=1}^n \bar\xi_r\,\frac{\partial\varphi}{\partial\bar{z}_r},\end{equation}

\noindent where for the last identity we used again (\ref{eqn:xi_varphi_def_coordinates}). On the other hand, writing (\ref{eqn:xi_varphi_def_coordinates}) in matrix form, we see that (\ref{eqn:xi_varphi_def_coordinates}) is equivalent to

\begin{equation}\label{eqn:xi_varphi_del-bar_varphi_coordinates}\xi_r = -i\,\sum\limits_{l=1}^n\frac{\partial\varphi}{\partial\bar{z}_l}\,\omega^{r\bar{l}} \hspace{2ex} \mbox{for all} \hspace{2ex} r=1,\dots , n.\end{equation}

\noindent Plugging the conjugate of (\ref{eqn:xi_varphi_del-bar_varphi_coordinates}) into (\ref{eqn:i-xi_varphi-del_varphi_norm}), we get

\begin{equation}\label{eqn:i-xi_varphi-del_varphi_norm_final}i\,\xi_{\varphi}\lrcorner\partial\varphi = \sum\limits_{l,\,r=1}^n\frac{\partial\varphi}{\partial z_l}\,\frac{\partial\varphi}{\partial\bar{z}_r}\,\omega^{l\bar{r}} = |\partial\varphi|_{\omega}^2.\end{equation}

\noindent This proves (\ref{eqn:i-xi_varphi-del_varphi}) and (\ref{eqn:Witten-Laplacians_formulae}).   \hfill $\Box$

\begin{Cor}\label{Cor:3-space-decomp_Witten} The Witten Laplacians $\Delta'_{\varphi}$ and $\Delta''_{\varphi}$ are elliptic and induce $3$-space orthogonal decompositions of the space of smooth $(p,\,q)$-forms for all $p,q\in\{0,\dots , n\}$:

\begin{equation}\label{eqn:3-space-decomp_Witten}C^{\infty}_{p,\,q}(X,\,\C) = \ker\Delta'_{\varphi} \oplus \mbox{Im}\,\partial_{\varphi} \oplus \mbox{Im}\,\partial_{\varphi}^{\star} \hspace{2ex} \mbox{and} \hspace{2ex} C^{\infty}_{p,\,q}(X,\,\C) = \ker\Delta''_{\varphi} \oplus \mbox{Im}\,\bar\partial_{\varphi} \oplus \mbox{Im}\,\bar\partial_{\varphi}^{\star}\end{equation}

\noindent in which $\ker\Delta'_{\varphi}$, $\ker\Delta''_{\varphi}$ are finite-dimensional, $\ker\partial_{\varphi} = \ker\Delta'_{\varphi} \oplus \mbox{Im}\,\partial_{\varphi}$ and $\ker\bar\partial_{\varphi} = \ker\Delta''_{\varphi} \oplus \mbox{Im}\,\bar\partial_{\varphi}$.

\end{Cor}

\noindent {\it Proof.} Formulae (\ref{eqn:Witten-Laplacians_formulae}) show that $\Delta'_{\varphi}$ and $\Delta''_{\varphi}$ have the same principal parts as $\Delta'$ resp. $\Delta''$ which are known to be elliptic, hence they are themselves elliptic. Since $X$ is compact and $\Delta'_{\varphi}, \Delta''_{\varphi}$ are self-adjoint, standard elliptic theory ensures that $\ker\Delta'_{\varphi}$ and $\ker\Delta''_{\varphi}$ are finite-dimensional, $\mbox{Im}\Delta'_{\varphi}$ and $\mbox{Im}\Delta''_{\varphi}$ are closed and finite-codimensional and we have orthogonal decompositions

\begin{equation}\label{eqn:2-space-decomp_Witten}C^{\infty}_{p,\,q}(X,\,\C) = \ker\Delta'_{\varphi} \oplus \mbox{Im}\,\Delta'_{\varphi} \hspace{2ex} \mbox{and} \hspace{2ex} C^{\infty}_{p,\,q}(X,\,\C) = \ker\Delta''_{\varphi} \oplus \mbox{Im}\,\Delta''_{\varphi}.\end{equation}

\noindent Using these splittings and the integrability properties (\ref{eqn:integrability_anti-commutation}), we easily infer that $\mbox{Im}\,\Delta'_{\varphi}$ and $\mbox{Im}\,\Delta''_{\varphi}$ further split orthogonally as 

$$\mbox{Im}\,\Delta'_{\varphi} = \mbox{Im}\,\partial_{\varphi} \oplus \mbox{Im}\,\partial_{\varphi}^{\star} \hspace{2ex} \mbox{and} \hspace{2ex} \mbox{Im}\,\Delta''_{\varphi} = \mbox{Im}\,\bar\partial_{\varphi} \oplus \mbox{Im}\,\bar\partial_{\varphi}^{\star}$$

\noindent and we get (\ref{eqn:3-space-decomp_Witten}).   \hfill $\Box$

\begin{Prop}\label{Prop:Witten-comm-rel} For every compact Hermitian manifold $(X,\,\omega)$ and every smooth function $\varphi:X\rightarrow\R$, the following {\bf twisted commutation relations} hold:

\begin{eqnarray}\label{eqn:Witten-comm-rel}\nonumber (i)\,\,[\Lambda,\,\bar\partial_{\varphi}] & = & -i\,(\partial_{\varphi}^{\star} + \tau^{\star} + 2i\,\xi_{\varphi}\lrcorner\cdot) = -i\,(\partial^{\star} + \tau^{\star} + i\,\xi_{\varphi}\lrcorner\cdot); \\
\nonumber (ii)\,\,[\Lambda,\,\partial_{\varphi}] & = & i\,(\bar\partial_{\varphi}^{\star} + \bar\tau^{\star} - 2i\,\bar\xi_{\varphi}\lrcorner\cdot) = i\,(\bar\partial^{\star} + \bar\tau^{\star} - i\,\bar\xi_{\varphi}\lrcorner\cdot); \\
\nonumber (iii)\,\, [\bar\partial^{\star}_{\varphi},\,L] & = & i\,(\partial_{\varphi} + \tau + 2\partial\varphi\wedge\cdot) = i\,(\partial + \tau + \partial\varphi\wedge\cdot); \\
\nonumber (iv)\,\,[\partial^{\star}_{\varphi},\,L] & = & -i\,(\bar\partial_{\varphi} + \bar\tau + 2\bar\partial\varphi\wedge\cdot) = -i\,(\bar\partial + \bar\tau + \bar\partial\varphi\wedge\cdot).\end{eqnarray}

\noindent They induce the following Bochner-Kodaira-Nakano-type identity

\begin{equation}\label{eqn:BKN_Witten1}\Delta''_{\varphi} = \Delta'_{\varphi} + [\partial_{\varphi},\,\tau^{\star}] - [\bar\partial_{\varphi},\,\bar{\tau}^{\star}] + 2\,[\partial_{\varphi},\,i\,\xi_{\varphi}\lrcorner\cdot] + 2\,[\bar\partial_{\varphi},\,i\,\bar\xi_{\varphi}\lrcorner\cdot]\end{equation} 

\noindent which can be seen as the twisted version of (\ref{eqn:BKN_demailly1}).

\end{Prop}

\noindent {\it Proof.} As usual, it suffices to prove one of the identities $(i)-(iv)$ as the others follow by conjugation and adjunction. We shall prove $(iii)$. Using (\ref{def:Witten-twisting_adjoints}), we get

\vspace{1ex}

\hspace{15ex} $[\bar\partial^{\star}_{\varphi},\,L] = [\bar\partial^{\star},\,L] + [i\,\bar\xi_{\varphi}\lrcorner,\,L] =i\,(\partial + \tau) + i\,(\bar\xi_{\varphi}\lrcorner\omega)\wedge\cdot,$

\noindent where we have used $(iii)$ of (\ref{eqn:standard-comm-rel}). This proves $(iii)$ of since $\bar\xi_{\varphi}\lrcorner\omega = \partial\varphi$ by (\ref{eqn:xi_varphi_def}).

 From $(ii)$, we get $\bar\partial_{\varphi}^{\star} = -i\,[\Lambda,\,\partial_{\varphi}] - \bar\tau^{\star} + 2i\,\bar\xi_{\varphi}\lrcorner\cdot$, hence

\begin{equation}\label{eqn:Delta_varphi''_comm-rel}\Delta''_{\varphi} = [\bar\partial_{\varphi},\,\bar\partial_{\varphi}^{\star}] = -i\,[\bar\partial_{\varphi},\,[\Lambda,\,\partial_{\varphi}]] - [\bar\partial_{\varphi},\,\bar\tau^{\star}] + 2\,[\bar\partial_{\varphi},\,i\,\bar\xi_{\varphi}\lrcorner\cdot].\end{equation}

\noindent Jacobi's identity reads

\vspace{1ex}

$-[\bar\partial_{\varphi},\,[\Lambda,\,\partial_{\varphi}]] + [\Lambda,\,[\partial_{\varphi},\,\bar\partial_{\varphi}]] + [\partial_{\varphi},\,[\bar\partial_{\varphi},\,\Lambda]] = 0, \hspace{2ex} \mbox{hence} \hspace{2ex} -[\bar\partial_{\varphi},\,[\Lambda,\,\partial_{\varphi}]] =  [\partial_{\varphi},\,[\Lambda,\,\bar\partial_{\varphi}]]$

\vspace{1ex}

\noindent because $[\partial_{\varphi},\,\bar\partial_{\varphi}] = 0$ by (\ref{eqn:integrability_anti-commutation}), hence $-i\,[\bar\partial_{\varphi},\,[\Lambda,\,\partial_{\varphi}]] = [\partial_{\varphi},\,\partial_{\varphi}^{\star} + \tau^{\star} + 2i\,\xi_{\varphi}\lrcorner\cdot]$ thanks to $(i)$. Plugging this into (\ref {eqn:Delta_varphi''_comm-rel}) and using $[\partial_{\varphi},\,\partial_{\varphi}^{\star}] = \Delta'_{\varphi}$, we get (\ref{eqn:BKN_Witten1}).  \hfill  $\Box$

\vspace{2ex}

 We shall now absorb the first-order terms on the r.h.s. of (\ref{eqn:BKN_Witten1}) into a new twisted Laplacian to obtain a twisted version of Demailly's formula (\ref{eqn:BKN_demailly2}).

\begin{The}\label{The:BKN_Witten2} For every compact Hermitian manifold $(X,\,\omega)$ and every smooth function $\varphi:X\rightarrow\R$, the following {\bf twisted Bochner-Kodaira-Nakano-type identity} holds:

\begin{equation}\label{eqn:BKN_Witten2}\Delta''_{\varphi} = \Delta'_{\varphi,\,\tau} -2\,[i\partial\bar\partial\varphi\wedge\cdot,\,\Lambda] + T_{\omega},\end{equation}

\noindent where $\Delta'_{\varphi,\,\tau}:=[\partial_{\varphi} + \tau + 2\partial\varphi\wedge\cdot,\, (\partial_{\varphi} + \tau + 2\partial\varphi\wedge\cdot)^{\star}] = [\partial_{\varphi} + \tau + 2\partial\varphi\wedge\cdot,\, \partial_{\varphi}^{\star} + \tau^{\star} + 2i\,\xi_{\varphi}\lrcorner\cdot]$ and $T_{\omega}:=\bigg[\Lambda,\,[\Lambda,\,\frac{i}{2}\,\partial\bar\partial\omega]\bigg] - [\partial\omega\wedge\cdot,\,(\partial\omega\wedge\cdot)^{\star}]$. Clearly, $\Delta'_{\varphi,\,\tau}$ is a non-negative, self-adjoint, elliptic operator of order $2$ having the same principal part as $\Delta'$, while $T_{\omega}$ is the same zeroth-order operator depending only on the torsion of $\omega$ as in Demailly's formula (\ref{eqn:BKN_demailly2}).

\end{The}

\noindent {\it Proof.} Starting from (\ref{eqn:BKN_Witten1}), we notice that

$$-[\bar\partial_{\varphi},\,\bar{\tau}^{\star}] = -[\bar\partial,\,\bar{\tau}^{\star}] + [\bar\partial\varphi\wedge\cdot,\,\bar{\tau}^{\star}] = [\tau,\,\partial^{\star} + \tau^{\star}] + T_{\omega} + [\bar\partial\varphi\wedge\cdot,\,\bar{\tau}^{\star}],$$

\noindent where the last identity followed from $(iv)$ of (\ref{eqn:BKN_demailly_auxiliary}). Plugging this into (\ref{eqn:BKN_Witten1}), we get

\begin{equation}\label{eqn:BKN_Witten2_prelim1}\Delta''_{\varphi} = [\partial_{\varphi} + \tau,\, \partial_{\varphi}^{\star} + \tau^{\star} + 2i\,\xi_{\varphi}\lrcorner\cdot] - [\tau,\,i\,\xi_{\varphi}\lrcorner\cdot] + [\bar\partial\varphi\wedge\cdot,\,\bar\tau^{\star}] + 2\,[\bar\partial_{\varphi},\,i\,\bar\xi_{\varphi}\lrcorner\cdot] + T_{\omega}.\end{equation}

\begin{Claim}\label{Claim1} We have: $- [\tau,\,i\,\xi_{\varphi}\lrcorner\cdot] + [\bar\partial\varphi\wedge\cdot,\,\bar\tau^{\star}] = 0$.\end{Claim}

 To prove this claim, we start by noticing that

\begin{equation}\label{eqn:proof-claim1_1}[\bar\partial\varphi\wedge\cdot,\,\bar\tau^{\star}] = [\bar\partial\varphi\wedge\cdot,\,\frac{i}{2}\,[\Lambda,\,\tau]] = \frac{i}{2}\,[\Lambda,\,[\tau,\,\bar\partial\varphi\wedge\cdot]] + \frac{i}{2}\,[\tau,\,[\bar\partial\varphi\wedge\cdot,\,\Lambda]],\end{equation}

\noindent where the first identity above followed from $(ii)$ of (\ref{eqn:BKN_demailly_auxiliary}) and the second identity is the Jacobi identity.

 Now, $[\tau,\,\bar\partial\varphi\wedge\cdot] = [\bar\partial\varphi\wedge\cdot,\,[\Lambda,\,\partial\omega\wedge\cdot]] = [\Lambda,\,[\partial\omega\wedge\cdot,\, \bar\partial\varphi\wedge\cdot]] + [\partial\omega\wedge\cdot,\,[\bar\partial\varphi\wedge\cdot,\, \Lambda]]$, hence 

\begin{equation}\label{eqn:comm_tau-delbar-varphi}[\tau,\,\bar\partial\varphi\wedge\cdot] = [\partial\omega\wedge\cdot,\,[\bar\partial\varphi\wedge\cdot,\, \Lambda]]\end{equation}

\noindent since $[\partial\omega\wedge\cdot,\, \bar\partial\varphi\wedge\cdot] = 0$. Now, $[\bar\partial\varphi\wedge\cdot,\, \Lambda]u = \bar\partial\varphi\wedge\Lambda u - \Lambda( \bar\partial\varphi\wedge u)$ for any form $u$. Thanks to Claim \ref{Claim2} below, this means that

\begin{equation}\label{eqn:comm_delbar_varphi-Lambda}[\bar\partial\varphi\wedge\cdot,\, \Lambda] = \xi_{\varphi}\lrcorner\cdot, \hspace{2ex}  \mbox{hence (\ref{eqn:comm_tau-delbar-varphi}) becomes} \hspace{2ex} [\tau,\,\bar\partial\varphi\wedge\cdot] = [\partial\omega\wedge\cdot,\,  \xi_{\varphi}\lrcorner\cdot] = (\xi_{\varphi}\lrcorner\partial\omega)\wedge\cdot.\end{equation}

\begin{Claim}\label{Claim2} For any differential form $u$, we have: $\Lambda(\bar\partial\varphi\wedge u) = -\xi_{\varphi}\lrcorner u + \bar\partial\varphi\wedge\Lambda u$.\end{Claim}

\noindent To prove Claim \ref{Claim2}, let $v$ be any smooth form of the same bidegree as $\Lambda(\bar\partial\varphi\wedge u)$. Then

\begin{eqnarray}\nonumber \langle\langle\Lambda(\bar\partial\varphi\wedge u),\, v\rangle\rangle & = & \langle\langle\bar\partial\varphi\wedge u,\, \omega\wedge v\rangle\rangle = \langle\langle u,\, -i\,\bar\xi_{\varphi}\lrcorner(\omega\wedge v)\rangle\rangle \\
\nonumber & = & \langle\langle u,\, -i\,\partial\varphi\wedge v\rangle\rangle + \langle\langle u,\, -\omega\wedge i\,\bar\xi_{\varphi}\lrcorner v\rangle\rangle =  \langle\langle -\xi_{\varphi}\lrcorner u,\,v\rangle\rangle +  \langle\langle \bar\partial\varphi\wedge\Lambda u,\,v\rangle\rangle,\end{eqnarray}

\noindent where (\ref{eqn:xi_varphi_def}) and Formulae \ref{Formulae:adjoints} have been used. Thus, Claim \ref{Claim2} is proved.

 Putting (\ref{eqn:proof-claim1_1}) and (\ref{eqn:comm_delbar_varphi-Lambda}) together, we get

\begin{equation}\label{eqn:proof-claim1_2}[\bar\partial\varphi\wedge\cdot,\,\bar\tau^{\star}] = \frac{i}{2}\,[\Lambda,\,(\xi_{\varphi}\lrcorner\partial\omega)\wedge\cdot] + \frac{1}{2}\,[\tau,\,i\,\xi_{\varphi}\lrcorner\cdot],\end{equation}

\noindent To compute $[\Lambda,\,(\xi_{\varphi}\lrcorner\partial\omega)\wedge\cdot]$, notice that $[\Lambda,\,(\xi_{\varphi}\lrcorner\partial\omega)\wedge\cdot]\,u = \Lambda((\xi_{\varphi}\lrcorner\partial\omega)\wedge u) - (\xi_{\varphi}\lrcorner\partial\omega)\wedge\Lambda u$ for every form $u$ and that for every form $v$ of the same bidegree as $\Lambda((\xi_{\varphi}\lrcorner\partial\omega)\wedge u)$, we have

\begin{eqnarray}\nonumber\langle\langle\Lambda((\xi_{\varphi}\lrcorner\partial\omega)\wedge u),\, v\rangle\rangle & = & \langle\langle\xi_{\varphi}\lrcorner(\partial\omega\wedge u) + \partial\omega\wedge(\xi_{\varphi}\lrcorner u),\,\omega\wedge v\rangle\rangle \\
\nonumber & = & \langle\langle\partial\omega\wedge u,\,i\partial\varphi\wedge\omega\wedge v\rangle\rangle +  \langle\langle\Lambda(\partial\omega\wedge(\xi_{\varphi}\lrcorner u)),\,v\rangle\rangle \\
\nonumber & = & \langle\langle\xi_{\varphi}\lrcorner\Lambda(\partial\omega\wedge u),\, v\rangle\rangle  + \langle\langle\Lambda(\partial\omega\wedge(\xi_{\varphi}\lrcorner u)),\,v\rangle\rangle.\end{eqnarray}

\noindent Hence, $\Lambda((\xi_{\varphi}\lrcorner\partial\omega)\wedge u) = \xi_{\varphi}\lrcorner\Lambda(\partial\omega\wedge u) + \Lambda(\partial\omega\wedge(\xi_{\varphi}\lrcorner u))$, so (\ref{eqn:proof-claim1_2}) gives

\begin{equation}\label{eqn:proof-claim1_3}[\bar\partial\varphi\wedge\cdot,\,\bar\tau^{\star}]\,u = \frac{1}{2}\,\bigg(i\,\xi_{\varphi}\lrcorner\Lambda(\partial\omega\wedge u) + \Lambda(\partial\omega\wedge(i\,\xi_{\varphi}\lrcorner u)) - (i\,\xi_{\varphi}\lrcorner\partial\omega)\wedge\Lambda u\bigg) + \frac{1}{2}\,[\tau,\,i\,\xi_{\varphi}\lrcorner\cdot]\,u.\end{equation}

 Meanwhile, we have

\begin{eqnarray}\label{eqn:proof-claim1_4}\nonumber [\tau,\,i\,\xi_{\varphi}\lrcorner\cdot]\,u & = & \tau(i\,\xi_{\varphi}\lrcorner u) + i\,\xi_{\varphi}\lrcorner\tau u \\
\nonumber & = & \Lambda(\partial\omega\wedge(i\,\xi_{\varphi}\lrcorner u)) - \partial\omega\wedge\Lambda(i\,\xi_{\varphi}\lrcorner u) + i\,\xi_{\varphi}\lrcorner\Lambda(\partial\omega\wedge u) - i\,\xi_{\varphi}\lrcorner(\partial\omega\wedge\Lambda u)\\
\nonumber    & = & \bigg(i\,\xi_{\varphi}\lrcorner\Lambda(\partial\omega\wedge u) + \Lambda(\partial\omega\wedge(i\,\xi_{\varphi}\lrcorner u)) - (i\,\xi_{\varphi}\lrcorner\partial\omega)\wedge\Lambda u\bigg) + \partial\omega\wedge[i\,\xi_{\varphi}\lrcorner\cdot,\,\Lambda]\,u\\
    & = & \bigg(i\,\xi_{\varphi}\lrcorner\Lambda(\partial\omega\wedge u) + \Lambda(\partial\omega\wedge(i\,\xi_{\varphi}\lrcorner u)) - (i\,\xi_{\varphi}\lrcorner\partial\omega)\wedge\Lambda u\bigg),\end{eqnarray}

\noindent where the last identity followed from $[i\,\xi_{\varphi}\lrcorner\cdot,\,\Lambda] = [L,\,\partial\varphi\wedge\cdot]^{\star} = 0$.

 Using (\ref{eqn:proof-claim1_4}), (\ref{eqn:proof-claim1_3}) reads $[\bar\partial\varphi\wedge\cdot,\,\bar\tau^{\star}]\,u = [\tau,\,i\,\xi_{\varphi}\lrcorner\cdot]\,u$ for every form $u$. This proves Claim \ref{Claim1}. Using Claim \ref{Claim1}, (\ref{eqn:BKN_Witten2_prelim1}) reduces to

\begin{equation}\label{eqn:BKN_Witten2_prelim2}\Delta''_{\varphi} = [\partial_{\varphi} + \tau,\, \partial_{\varphi}^{\star} + \tau^{\star} + 2i\,\xi_{\varphi}\lrcorner\cdot] + 2\,[\bar\partial_{\varphi},\,i\,\bar\xi_{\varphi}\lrcorner\cdot] + T_{\omega}.\end{equation}

 We shall now compute the term by which the r.h.s. of (\ref{eqn:BKN_Witten2_prelim2}) still falls short of the expected $\Delta'_{\varphi,\,\tau}$. This is the quantity

\begin{eqnarray}\label{eqn:missing-quantity0}[2\partial\varphi\wedge\cdot,\,\partial_{\varphi}^{\star} + \tau^{\star} + 2i\,\xi_{\varphi}\lrcorner\cdot] = [2\partial\varphi\wedge\cdot,\,\partial^{\star} + \tau^{\star}] + [2\partial\varphi\wedge\cdot,\,i\,\xi_{\varphi}\lrcorner\cdot].\end{eqnarray}

\noindent For every form $u$, we have $[\partial\varphi\wedge\cdot,\,i\,\xi_{\varphi}\lrcorner\cdot]\,u = \partial\varphi\wedge(i\,\xi_{\varphi}\lrcorner u) + i\,\xi_{\varphi}\lrcorner(\partial\varphi\wedge u) = (i\,\xi_{\varphi}\lrcorner\partial\varphi)\, u = |\partial\varphi|^2_{\omega}\, u$ (cf. (\ref{eqn:i-xi_varphi-del_varphi}) for the last identity). This means that

\begin{equation}\label{eqn:del-varphi-xi-varphi}[2\,\partial\varphi\wedge\cdot,\,i\,\xi_{\varphi}\lrcorner\cdot] = 2|\partial\varphi|^2_{\omega}.\end{equation}

\noindent Meanwhile, we know that $\partial^{\star} + \tau^{\star} = i\,[\Lambda,\,\bar\partial]$ by the standard commutation relation $(i)$ of (\ref{eqn:standard-comm-rel}). Thus

\begin{eqnarray}\label{eqn:del-varphi-delstar-taustar_1}[\partial\varphi\wedge\cdot,\,\partial^{\star} + \tau^{\star}] = [\partial\varphi\wedge\cdot,\,i\,[\Lambda,\,\bar\partial]] = [\bar\partial,\,i\,[\partial\varphi\wedge\cdot,\,\Lambda]] + [\Lambda,\,i\,[\bar\partial,\,\partial\varphi\wedge\cdot]],\end{eqnarray}

\noindent where the last identity is the Jacobi identity. Now, $[\bar\partial,\,\partial\varphi\wedge\cdot]\,u = \bar\partial(\partial\varphi\wedge u) -  \partial\varphi\wedge\bar\partial u = \bar\partial\partial\varphi\wedge u$ for every form $u$, hence $[\bar\partial,\,\partial\varphi\wedge\cdot] = \bar\partial\partial\varphi\wedge\cdot$. We also have $[\partial\varphi\wedge\cdot,\,\Lambda] = \bar\xi_{\varphi}\lrcorner\cdot$ by taking conjugates in (\ref{eqn:comm_delbar_varphi-Lambda}). Thus (\ref{eqn:del-varphi-delstar-taustar_1}) becomes

\begin{eqnarray}\label{eqn:del-varphi-delstar-taustar_2}[\partial\varphi\wedge\cdot,\,\partial^{\star} + \tau^{\star}] = [\bar\partial,\,i\,\bar\xi_{\varphi}\lrcorner\cdot] - [\Lambda,\,i\,\partial\bar\partial\varphi\wedge\cdot].\end{eqnarray}

\noindent Putting together (\ref{eqn:missing-quantity0}), (\ref{eqn:del-varphi-xi-varphi}) and (\ref{eqn:del-varphi-delstar-taustar_1}), we get

\begin{eqnarray}\label{eqn:missing-quantity1}[2\partial\varphi\wedge\cdot,\,\partial_{\varphi}^{\star} + \tau^{\star} + 2i\,\xi_{\varphi}\lrcorner\cdot] = 2\,[\bar\partial,\,i\,\bar\xi_{\varphi}\lrcorner\cdot] - 2\,[\Lambda,\,i\,\partial\bar\partial\varphi\wedge\cdot] + 2|\partial\varphi|^2_{\omega}.\end{eqnarray}

\noindent On the other hand, the second term on the r.h.s. of (\ref{eqn:BKN_Witten2_prelim2}) reads $2\,[\bar\partial_{\varphi},\,i\,\bar\xi_{\varphi}\lrcorner\cdot] = 2\,[\bar\partial,\,i\,\bar\xi_{\varphi}\lrcorner\cdot] -2\,[\bar\partial\varphi\wedge\cdot,\,i\,\bar\xi_{\varphi}\lrcorner\cdot]$, so using the expression for  $2\,[\bar\partial,\,i\,\bar\xi_{\varphi}\lrcorner\cdot]$ given by (\ref{eqn:missing-quantity0}), we get

\begin{eqnarray}\label{eqn:missing-quantity2}\nonumber 2\,[\bar\partial_{\varphi},\,i\,\bar\xi_{\varphi}\lrcorner\cdot] & = & [2\partial\varphi\wedge\cdot,\,\partial_{\varphi}^{\star} + \tau^{\star} + 2i\,\xi_{\varphi}\lrcorner\cdot] + 2\,[\Lambda,\,i\,\partial\bar\partial\varphi\wedge\cdot] - 2|\partial\varphi|^2_{\omega} -2\,[\bar\partial\varphi\wedge\cdot,\,i\,\bar\xi_{\varphi}\lrcorner\cdot] \\
   & = & [2\partial\varphi\wedge\cdot,\,\partial_{\varphi}^{\star} + \tau^{\star} + 2i\,\xi_{\varphi}\lrcorner\cdot] + 2\,[\Lambda,\,i\,\partial\bar\partial\varphi\wedge\cdot],\end{eqnarray}

\noindent where the last identity is due to the fact that $-2\,[\bar\partial\varphi\wedge\cdot,\,i\,\bar\xi_{\varphi}\lrcorner\cdot] = 2|\partial\varphi|^2_{\omega}$ which follows, in turn, by taking conjugates in (\ref{eqn:del-varphi-xi-varphi}).

 It is now clear that (\ref{eqn:BKN_Witten2_prelim2}) and (\ref{eqn:missing-quantity2}) prove between them the desired identity (\ref{eqn:BKN_Witten2}). The proof of Theorem \ref{The:BKN_Witten2} is complete.  \hfill $\Box$

\vspace{3ex}

\subsection{Torsion and spectral gap}\label{subsection:torsion-spectral-gap}

 We have seen in Lemma \ref{Lem:reduction-ineq-starstar} that inequality (\ref{eqn:reduction-ineq-starstar}) suffices to guarantee the degeneration at $E_2$ of the Fr\"olicher spectral sequence of $X$. We spell out in this section a few calculations aimed at bounding below the term $(1-\varepsilon-\delta)\,(||p''_{\perp}\partial u||^2 + ||p''_{\perp}\partial^{\star} u||^2)$ featuring on the l.h.s. of (\ref{eqn:reduction-ineq-starstar}). The arguments resemble those in $\S.$ \ref{subsection:use-of-b}.

 Since $X$ is compact, the spectres of the non-negative elliptic operators $\Delta'=\Delta^{'p,\,q},\Delta''= \Delta^{''p,\,q}:C^{\infty}_{p,\,q}(X,\,\C)\longrightarrow C^{\infty}_{p,\,q}(X,\,\C)$ are discrete, contained in $[0,\,+\infty)$ and have no finite accumulation points, so each of them has a smallest positive eigenvalue. We denote these eigenvalues by

\begin{equation}\label{eqn:lambda_0mu_0_def}\lambda_0^{p,\,q}:=\mbox{min}\,\bigg(\mbox{Spec}\,\Delta^{'p,\,q}\cap(0,\,+\infty)\bigg)>0 \hspace{2ex} \mbox{and} \hspace{2ex}  \mu_0^{p,\,q}:=\mbox{min}\,\bigg(\mbox{Spec}\,\Delta^{''p,\,q}\cap(0,\,+\infty)\bigg)>0.\end{equation}

\begin{Lem}\label{Lem:p''_perp_lbound} Suppose the metric $\omega$ is SKT and the following condition is satisfied on $X$:

\begin{equation}\label{eqn:p''_perp_lbound_hyp}\lambda_0^{p,\,q} > ||[\bar\tau,\,\bar\tau^{\star}] - [\bar\partial\omega\wedge\cdot,\,(\bar\partial\omega\wedge\cdot)^{\star}]||_{p,\,q}:=\sup\limits_{u\in C^{\infty}_{p,\,q}(X,\,\C),\,\,||u||=1}\langle\langle([\bar\tau,\,\bar\tau^{\star}] - [\bar\partial\omega\wedge\cdot,\,(\bar\partial\omega\wedge\cdot)^{\star}])\,u,\,u\rangle\rangle.\end{equation}

\noindent Then each of the following equivalent statements holds. 

\vspace{1ex}

 $(i)$\, There exits a constant $0<\varepsilon_1<1$ such that 

\begin{equation}\label{eqn:p''_perp_lbound}||p''_{\perp}v||^2 \geq (1-\varepsilon_1)\,||v||^2 \hspace{2ex} \mbox{for all}\hspace{1ex} v\in C^{\infty}_{p,\,q}(X,\,\C) \hspace{1ex} \mbox{such that} \hspace{1ex} v\in(\ker\Delta')^{\perp}; \end{equation}

\vspace{1ex}

 $(ii)$\, The map $p''_{\perp}:(\ker\Delta')^{\perp} \longrightarrow (\ker\Delta'')^{\perp}$ is injective;

\vspace{1ex}

 $(iii)$\, $\ker\Delta''\cap(\ker\Delta')^{\perp} = \{0\}.$

\end{Lem}

\noindent {\it Proof.} The equivalences $(i)\Leftrightarrow (ii)\Leftrightarrow (iii)$ are clear (using the Open Mapping Theorem in Fr\'echet spaces), while the fact that the positive constant in (\ref{eqn:p''_perp_lbound}) is smaller than $1$ follows from the obvious inequality $||p''_{\perp}v||\leq||v||$ for all $v$. Note that the conjugate of Demailly's non-K\"ahler Bochner-Kodaira-Nakano identity $\Delta'' = \Delta'_{\tau} + T_{\omega}$ (cf. (\ref{eqn:BKN_demailly2})) is

\begin{equation}\label{eqn:BKN_conjugate_demailly2_appendix}\Delta' = \Delta''_{\tau} + \overline{T}_{\omega},\end{equation}

\noindent where $\Delta''_{\tau}:=[\bar\partial + \bar\tau,\, \bar\partial^{\star} + \bar\tau^{\star}]$ and $\overline{T}_{\omega} = [\Lambda,\,[\Lambda,\,\frac{i}{2}\,\partial\bar\partial\omega\wedge\cdot]] - [\bar\partial\omega\wedge\cdot,\,(\bar\partial\omega\wedge\cdot)^{\star}]$.

 To prove $(iii)$, let $v\in\ker\Delta''\cap(\ker\Delta')^{\perp}$. Then, using (\ref{eqn:BKN_conjugate_demailly2_appendix}), we get the identity below:

\begin{eqnarray}\label{eqn:Delta'-Delta''_comparison_ineq}\nonumber \lambda_0^{p,\,q}\,||v||^2 \leq \langle\langle\Delta'v,\,v\rangle\rangle & = & \langle\langle\Delta''v,\,v\rangle\rangle + \langle\langle[\bar\partial,\,\bar\tau^{\star}]\,v,\,v\rangle\rangle + \langle\langle[\bar\tau,\,\bar\partial^{\star}]\,v,\,v\rangle\rangle \\
  & + & \langle\langle[\bar\tau,\,\bar\tau^{\star}]\,v,\,v\rangle\rangle - \langle\langle[\bar\partial\omega\wedge\cdot,\,(\bar\partial\omega\wedge\cdot)^{\star}]\,v,\,v\rangle\rangle,\end{eqnarray}

\noindent where the SKT assumption on $\omega$ has been used to have the signless term in $\overline{T}_{\omega}$ vanish. Since $v\in\ker\Delta'' = \ker\bar\partial\cap\ker\bar\partial^{\star}$, the first three terms on the r.h.s. of (\ref{eqn:Delta'-Delta''_comparison_ineq}) vanish. Indeed, \begin{eqnarray}\label{eqn:Delta''-torsion-terms_vanishing_appendix}\nonumber \langle\langle[\bar\partial,\,\bar\tau^{\star}]\,v,\,v\rangle\rangle & = & \langle\langle\bar\tau^{\star}\,v,\,\bar\partial^{\star}v\rangle\rangle + \langle\langle\bar\partial v,\,\bar\tau v\rangle\rangle = 0 + 0 = 0, \\ 
\nonumber \langle\langle[\bar\tau,\,\bar\partial^{\star}]\,v,\,v\rangle\rangle & = & \langle\langle\bar\partial^{\star}\,v,\,\bar\tau^{\star}v\rangle\rangle + \langle\langle\bar\tau v,\,\bar\partial v\rangle\rangle = 0 + 0 = 0.\end{eqnarray}

\noindent Thus, for every $(p,\,q)$-form $v\in\ker\Delta''\cap(\ker\Delta')^{\perp}$, (\ref{eqn:Delta'-Delta''_comparison_ineq}) reduces to

\begin{eqnarray}\label{eqn:Delta'-Delta''_comparison_ineq1}\nonumber \lambda_0^{p,\,q}\,||v||^2 \leq \langle\langle\Delta'v,\,v\rangle\rangle & = & \langle\langle([\bar\tau,\,\bar\tau^{\star}] - [\bar\partial\omega\wedge\cdot,\,(\bar\partial\omega\wedge\cdot)^{\star}])\,v,\,v\rangle\rangle\\
 & \leq & ||[\bar\tau,\,\bar\tau^{\star}] - [\bar\partial\omega\wedge\cdot,\,(\bar\partial\omega\wedge\cdot)^{\star}]||_{p,\,q}\,||v||^2.\end{eqnarray}

\noindent If $v\neq 0$, this is clearly impossible under the assumption (\ref{eqn:p''_perp_lbound_hyp}). So $v=0$ and $(iii)$ is proved.  \hfill $\Box$

\begin{Cor}\label{Cor:p''_perp_lbound} If $\omega$ is SKT and if condition (\ref{eqn:p''_perp_lbound_hyp}) is satisfied for all $(p,\,q)$, then there exists a constant $0<\varepsilon_1<1$ such that

\begin{equation}\label{eqn:p''_perp-del_p''_perp-del-star} ||p''_{\perp}\partial u||^2 + ||p''_{\perp}\partial^{\star} u||^2 \geq (1-\varepsilon_1)\,\langle\langle\Delta'u,\,u\rangle\rangle\end{equation}

\noindent for all $(p,\,q)$ and all $u\in C^{\infty}_{p,\,q}(X,\,\C)$. In other words, $\Delta'_{p''_{\perp}} \geq (1-\varepsilon_1)\,\Delta'$.

\end{Cor}

\noindent {\it Proof.} Thanks to the first of the orthogonal $3$-space decompositions (\ref{eqn:3-space-decomp}), both choices $v=\partial u$ and $v=\partial^{\star}u$ satisfy the requirement $v\in(\ker\Delta')^{\perp}$, so (\ref{eqn:p''_perp_lbound}) applies. Then (\ref{eqn:p''_perp-del_p''_perp-del-star}) follows from (\ref{eqn:p''_perp_lbound}) and from $||\partial u||^2 + ||\partial^{\star}u||^2 = \langle\langle\Delta'u,\,u\rangle\rangle$.  \hfill $\Box$

\newpage

\noindent {\bf References} \\

\noindent [BS76]\, L. Boutet de Monvel, J. Sj\"ostrand --- {\it Sur la singularit\'e des noyaux de Bergman et de Szeg\"o} --- Soc. Math. de France, Ast\'erisque {\bf 34-35} (1976), 123-164.

\vspace{1ex}

\noindent [Del68]\, P. Deligne --- {\it Th\'eor\`eme de Lefschetz et crit\`eres de d\'eg\'en\'erescence de suites spectrales} --- Publ. Math. IHES {\bf 35} (1968), 107-126.

\vspace{1ex}

\noindent [Dem 84]\, J.-P. Demailly --- {\it Sur l'identit\'e de Bochner-Kodaira-Nakano en g\'eom\'etrie hermitienne} --- S\'eminaire d'analyse P. Lelong, P. Dolbeault, H. Skoda (editors) 1983/1984, Lecture Notes in Math., no. {\bf 1198}, Springer Verlag (1986), 88-97.

\vspace{1ex}

\noindent [Dem 96]\, J.-P. Demailly --- {\it Th\'eorie de Hodge $L^2$ et th\'eor\`emes d'annulation} --- in ``Introduction \`a la th\'eorie de Hodge'', J. Bertin, J.-P. Demailly, L. Illusie, C. Peters, Panoramas et Synth\`eses {\bf 3}, SMF (1996).

\vspace{1ex}

\noindent [Dem 97]\, J.-P. Demailly --- {\it Complex Analytic and Algebraic Geometry}---http://www-fourier.ujf-grenoble.fr/~demailly/books.html

\vspace{1ex}

\noindent [LRY15]\, K. Liu, S. Rao, X. Yang --- {\it Quasi-isometry and Deformations of Calabi-Yau Manifolds} --- Invent. Math. {\bf 199} (2015), no. 2, 423-453.

\vspace{1ex}

\noindent [MS74]\, A. Melin, J. Sj\"ostrand --- {\it Fourier Integral Operators with Complex-Valued Phase Functions} --- in ``Fourier Integral Operators and Partial Differential Equations'', Lecture Notes in Mathematics, Springer, Vol. {\bf 459} (1974), 120- 223.

\vspace{1ex}

\noindent [Pop13]\, D. Popovici --- {\it Deformation Limits of Projective Manifolds\!\!: Hodge Numbers and Strongly Gauduchon Metrics} --- Invent. Math. {\bf 194} (2013), 515-534.

\vspace{1ex}

\noindent [Raw77]\, J.H. Rawnsley --- {\it On the Cohomology Groups of a Polarisation and Diagonal Quantisation} --- Trans. Amer. Math. Soc. {\bf 230} (1977) 235- 255.

\vspace{1ex}

\noindent [Rei58]\, B. Reinhart --- {\it Harmonic Integrals on Almost Product Manifolds} --- Trans. Amer. Math. Soc. {\bf 88} (1958) 243-276.

\vspace{1ex}

\noindent [Rei59]\, B. Reinhart --- {\it Foliated Manifolds with Bundle-Like Metrics} --- Ann. of Math. (2) {\bf 69} (1959) 119-132.

\vspace{1ex}

\noindent [Sar78]\, K.S. Sarkaria --- {\it A Finiteness Theorem for Foliated Manifolds} --- J. Math. Soc. Japan {\bf 30}, No.4 (1978) 687-696.

\vspace{1ex}

\noindent [Voi02]\, C. Voisin --- {\it Hodge Theory and Complex Algebraic Geometry, $I$} --- {\bf 76} (2002), {\it Cambridge Studies in Advanced Mathematics}, Cambridge University Press, Cambridge. (Translated from the French original by L. Schneps.)

\vspace{1ex}

\noindent [Wit82]\, E. Witten --- {\it Supersymmetry and Morse Theory} --- J. Diff. Geom., {\bf 17} (1982), 661-692.

\vspace{10ex}

\noindent Institut de Math\'ematiques de Toulouse, Universit\'e Paul Sabatier,

\noindent 118 route de Narbonne, 31062 Toulouse, France

\noindent Email: popovici@math.univ-toulouse.fr

\end{document}